A Mathematical Monograph

# Substitution drawing rules on the Fibonacci word

A Research MSc from The Open University (UK)

## Martin Hansen

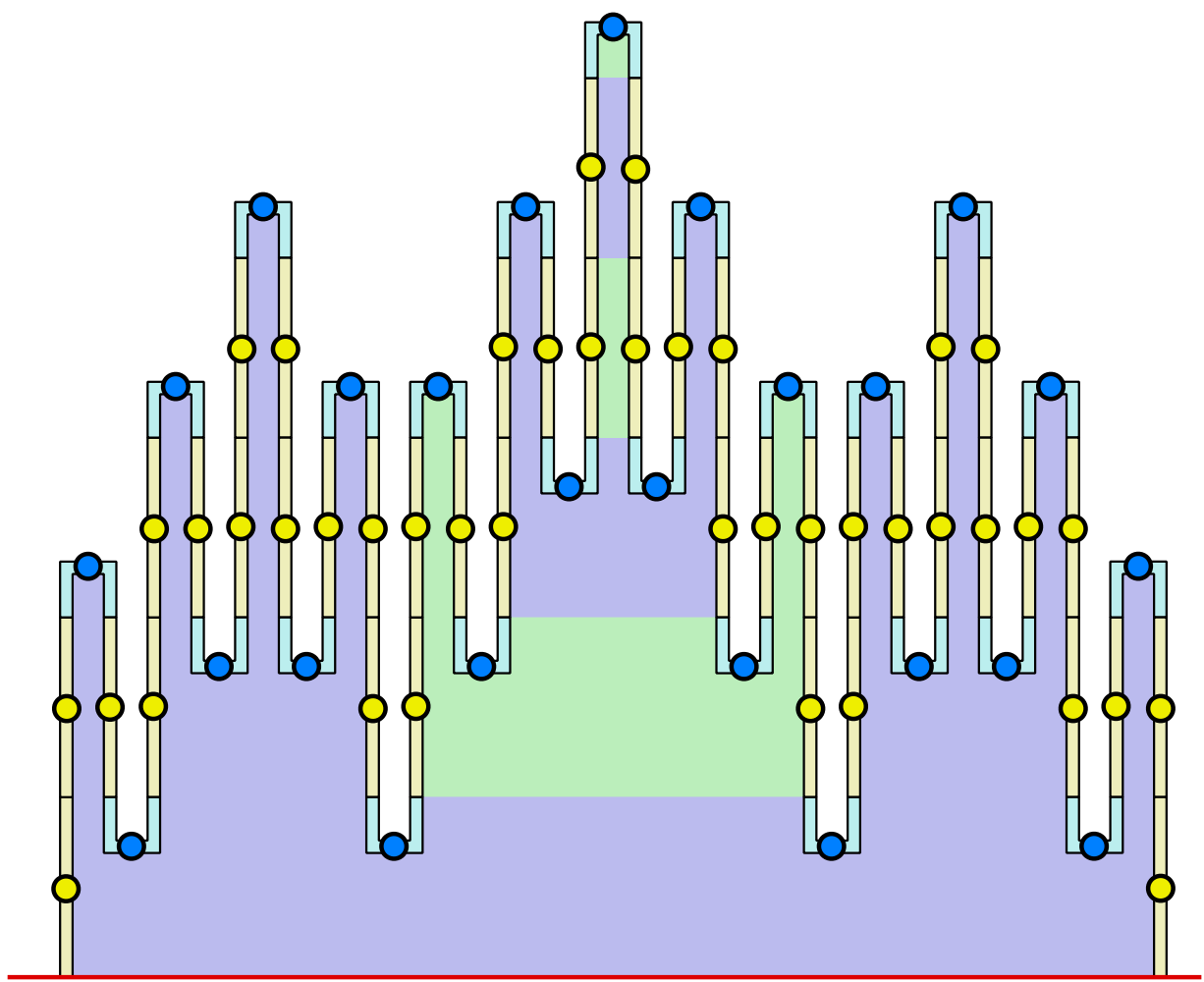

## MSc Dissertation

"Substitution drawing rules on the Fibonacci word" by Martin Hansen.

In fulfilment of the requirements for the degree of Master of Science under the supervision of Dr Reem Yassawi and Dr Dan Rust at the School of Mathematics and Statistics at The Open University, Milton Keynes, England.

## Dedication

To Helen for her unwavering support for the enterprise over the past year.

## Reviews & Comments

"Martin not only received a distinction for his dissertation but he achieved a mark that placed him in the top five students in his cohort of approximately eighty. He demonstrated a proficiency with several advanced combinatorial and computational methods in his work, to the point where he was able to prove modest statements that are, to my knowledge, new in the literature; an incredible feat for a student in their master's studies."

"Martin drew upon an impressive array of mathematical methods and subjects in this work, including but not limited to combinatorics on words, fractal geometry. linear algebra, number theory, mathematical induction and computer visualisation. This dissertation is full of detailed and well-reasoned proofs and calculations that comprise a true tour de force."

"Martin is an adept writer and conveys intricate mathematical arguments in a rigorous and thought-provoking way, with a flare for injecting a certain amount of fun and whimsy, to which I took great pleasure when reading his dissertation."

"Martin's dissertation was excellent and a joy to read, and really very impressive. Quite apart from the excellent mathematical content, which was at times intricate and fiddly (but precise and rigorous), and the original creation of a fractal out of the Fibonacci word (including the use of excellent computer generated illustrations to elucidate), what I noticed the most was the clear sense that Martin was enjoying himself; Martin conveys to the reader the joy and excitement of doing mathematics."

"Many of the proofs in this dissertation were of Martin's own results, exploring two extremely interesting geometric constructions. I particularly enjoyed that Martin was able to reuse Lemmas and Propositions from earlier sections throughout the dissertation, reinforcing that each section of the work is linked and builds from one section to the next. This helped to shape an engaging and convincing narrative."

# Substitution drawing rules on the Fibonacci word

A Research MSc from The Open University (UK)

# Martin Hansen

Martin Hansen
Shrewsbury
SY2 5DT
mhh@shrewsbury.org.uk
07841 288198 (text only)

# Abstract

This paper is a sharp and focussed exploration of the Fibonacci substitution and the mathematical entity it gives rise to, the Fibonacci word. Our investigations are both of an algebraic and a geometric nature. Indeed, it is the combination of the two that gives this paper its overall character. The work is in four parts. Chapter 1 is a brisk tour of necessary basics; definitions, key theorems, and a number of techniques subsequently used extensively. A simple one dimensional drawing rule is investigated in chapter 2 with the aid of what is thought to be an original geometric figure that we will call a "deviation from zero diagram". A highlight of the chapter is its concluding elementary proof of a non-trivial result. Chapter 3 presents a two dimensional drawing rule. Although selected because it is amongst the simplest form possible, this time the object derived from the rule is a fractal. That this is so is proven and its fractal (Hausdorff) dimension calculated. This fractal is a (possibly previously unexplored) variant of that known in the literature as "The Fibonacci Fractal". By way of an overall conclusion, the last chapter, the fourth, suggests a few aspects of those preceding it worthy of further analysis.

<div style="text-align: right">

Martin Hansen
The Open University
May 2022

</div>

# Symbols Used

In addition to standard mathematical notation and following symbols have a particular meaning in this work.

| Symbol | Meaning |
|---|---|
| $\mathcal{A}$ | Alphabet : A set of letters |
| $\mathcal{D}_n$ | The $n^{th}$ deviation from zero structure ($n$ is a catalogue entry). |
| $dim_{SYM}$ | Similarity dimension of a fractal |
| $\mathcal{E}_n$ | A piece of tiling associated with an embedded word |
| $\mathcal{F}_n$ | The Fibonacci word, $\mathcal{F}_n = \theta^n(a)$ |
| $\mathcal{F}_n^*$ | $\mathcal{F}_n$ with the leading $(aba)$ removed; $\mathcal{F}_n = (aba)\,\mathcal{F}_n^*$ |
| $|\mathcal{F}_n|$ | The number of letters in the word $\mathcal{F}_n$ |
| $\overrightarrow{\mathcal{F}_n}$ | Displacement vector between the start and finish of the associated path |
| $\mathcal{F}$ | The infinite Fibonacci word |
| $\mathcal{G}_n$ | A piece of tiling in a growth chart ($n$ is a catalogue entry) |
| $L_A$ | The length of a tile (tile $A$, for example) |
| $\mathbf{M}_\Phi$ | The incidence matrix (for the substitution $\Phi$, for example) |
| $\mathcal{T}_n$ | $\mathcal{F}_n$ with the last two letters exchanged, each with the other |
| $\mathcal{W}_n$ | The Fibonacci word $\mathcal{F}_n$ with the two rightmost letters removed |
| $\overrightarrow{\mathcal{W}_n}$ | Displacement vector between the start and finish of the associated path |
| $\phi$ | The golden ratio, $\phi = \frac{1+\sqrt{5}}{2}$ |
| $\lambda$ | An eigenvalue ($\lambda_{PF}$ is the Perron-Frobenius eigenvalue). |
| $\theta$ | The Fibonacci substitution ($a \to ab, b \to a$) |
| $\Theta$ | An inflation mapping on tile lengths |
| $\square$ | Denotes the end of a proof |

# Contents

**Word Count : Total is approximately 11853 words.**

| Chapter 1 | | | Chapter 2 | | | Chapter 3 | | | Chapter 4 | | |
|---|---|---|---|---|---|---|---|---|---|---|---|
| Page | N° | Total | Page | N° | Total | Page | N° | Total | Page | N° | Total |
| 07 | 467 | 00467 | 17 | 413 | 04110 | 39 | 256 | 07927 | 55 | 249 | 11228 |
| 08 | 270 | 00737 | 18 | 188 | 04298 | 40 | 195 | 08122 | 56 | 000 | 11228 |
| 09 | 346 | 01083 | 19 | 178 | 04476 | 41 | 315 | 08437 | 57 | 170 | 11398 |
| 10 | 475 | 01558 | 20 | 319 | 04795 | 42 | 312 | 08749 | 58 | 295 | 11693 |
| 11 | 533 | 02091 | 21 | 000 | 04795 | 43 | 257 | 09006 | 59 | 160 | 11853 |
| 12 | 299 | 02390 | 22 | 000 | 04795 | 44 | 103 | 09109 | | | |
| 13 | 384 | 02774 | 23 | 138 | 04933 | 45 | 171 | 09280 | | | |
| 14 | 383 | 03157 | 24 | 205 | 05138 | 46 | 293 | 09573 | | | |
| 15 | 232 | 03389 | 25 | 192 | 05330 | 47 | 182 | 09755 | | | |
| 16 | 308 | 03697 | 26 | 172 | 05502 | 48 | 063 | 09818 | | | |
| | | | 27 | 258 | 05760 | 49 | 000 | 09818 | | | |
| | | | 28 | 108 | 05868 | 50 | 111 | 09929 | | | |
| | | | 29 | 000 | 05868 | 51 | 227 | 10156 | | | |
| | | | 30 | 143 | 06011 | 52 | 329 | 10485 | | | |
| | | | 31 | 088 | 06099 | 53 | 233 | 10718 | | | |
| | | | 32 | 291 | 06390 | 54 | 261 | 10979 | | | |
| | | | 33 | 000 | 06390 | | | | | | |
| | | | 34 | 314 | 06704 | | | | | | |
| | | | 35 | 185 | 06889 | | | | | | |
| | | | 36 | 264 | 07153 | | | | | | |
| | | | 37 | 245 | 07398 | | | | | | |
| | | | 38 | 271 | 07669 | | | | | | |

# Chapter 1

# The Fibonacci Word

**1.1 Mathematical Beauty**

Mathematics is often described by its practitioners as being beautiful. "Beauty", some say, "is in the eye of the beholder". To me, mathematical beauty revolves around simple ideas giving rise to interesting entities with visualisations that can be sketched on paper and yet have extraordinary properties and depth. Over the following pages a fascinating object, the Fibonacci word, will be explained and explored. To set the scene and make this work self-contained, this first chapter is a minimalist introduction to the Fibonacci word. This initial material is well known. Of necessity, it begins with definitions of the objects of interest. These definitions are deliberately paired down and elementary to make this work as accessible as possible. Readers wanting a more technical treatment of the foundations and to see how they are embedded in a more general context, are recommended to read, for example, the classic textbook by Pytheas Fogg [Fogg]. The main contribution of this chapter to the established literature is to string together snippets found in the standard works (for example, [BGr13], and [ASh03]), and in the the use of a uniform notation and terminology. It strives to enhance established explanations, to untangle dependencies, and to present a narrative that covers the key steps, stands by itself, and flows. This chapter is sharply focussed on the subset of material needed to set up the Fibonacci word ready for the visualisations to be explored in the subsequent chapters. That is where the true beauty of the Fibonacci word starts to be revealed.

**1.2 Definitions**

By definition, a *word* is a finite or infinite sequence of elements, termed *letters*, all of which are taken from a finite set called an *alphabet*. In this dissertation we work exclusively with the two letter alphabet, $\mathcal{A} = \{\,a, b\,\}$. From this alphabet, for example, the letters *a*, *b*, *a*, *a* and *b* could be selected, in that order. Then, by *concatenation*, the word *abaab* formed. Concatenation is a simple placing of the selected letters one after the other, working left to right, to form a word. Words can themselves be concatenated. For example, *aba* concatenated with *ab* again forms the word *abaab*. The interest here is not in forming words from a random selection of letters from the alphabet but rather in starting with the simplest of words, a single isolated letter *a*, and repeatedly applying a *substitution*. In general a substitution can be thought of as "a replacement rule". Our interest is in the Fibonacci substitution, denoted $\theta$. Given a word, the Fibonacci substitution replaces each occurrences of the letter *a* with the letter pair *ab* and each occurrence of the letter *b* with the letter *a*. Mathematically $\theta$ is given by,

$$a \to ab$$
$$b \to a$$

When applied iteratively thrice to an initial letter *a*, the word *abaab* is obtained via the following steps,

$$a \to ab \to aba \to abaab$$

It could be written that $\theta^3(a) = abaab$, but we will write $\mathcal{F}_3 = abaab$ because the focus throughout this dissertation is entirely on the specific simple case of the Fibonacci substitution being applied iteratively to an initial isolated letter $a$. By definition, the $n^{th}$ Fibonacci word is given by,

$$\mathcal{F}_n = \theta^n(a)$$

Table 1.1 shows an initial letter $a$ having the Fibonacci substitution applied repeatedly, giving rise to the first few Fibonacci words.

| $n$ | $\mathcal{F}_n = \theta^n(a)$ | $\lvert a \rvert$ | $\lvert b \rvert$ | $\lvert \mathcal{F}_n \rvert$ |
|---|---|---|---|---|
| 0 | $\mathcal{F}_0 = a$ | 1 | 0 | 1 |
| 1 | $\mathcal{F}_1 = ab$ | 1 | 1 | 2 |
| 2 | $\mathcal{F}_2 = aba$ | 2 | 1 | 3 |
| 3 | $\mathcal{F}_3 = abaab$ | 3 | 2 | 5 |
| 4 | $\mathcal{F}_4 = abaababa$ | 5 | 3 | 8 |
| 5 | $\mathcal{F}_5 = abaababaabaab$ | 8 | 5 | 13 |
| 6 | $\mathcal{F}_6 = abaababaabaababaababa$ | 13 | 8 | 21 |
| 7 | $\mathcal{F}_7 = abaababaabaababaababaabaababaabaab$ | 21 | 13 | 34 |

**Table 1.1** : The Fibonacci substitution applied iteratively to an initial letter $a$.

With any iterative process, any object in the nature of a fixed point is of interest. With that in mind, notice that the word $\mathcal{F}_n$, corresponding to any given value of $n$, occurs at the leftmost end of all subsequent words. The infinite fixed word that results as $n$ tends to infinity is known as the *Fibonacci word, $\mathcal{F}$*, and has been described as "one of the most studied examples in the combinatorial theory of infinite words" [CRR14, page 40].

As table 1.1 suggests, successive Fibonacci words, $\mathcal{F}_n$, have the number of occurrences of the letter $a$, which we will denote $\lvert a \rvert$, the number of occurrences of the letter $b$, denoted $\lvert b \rvert$, and the overall word length, $\lvert \mathcal{F}_n \rvert$, all working their way through the world famous Fibonacci number sequence which is presented in table 1.2. With a little thought it can be seen that the number patterns of table 1.1 arise directly from the Fibonacci substitution ($a \to ab$, $b \to a$) where each letter in a previous word yields one letter $a$ in the next, and, in addition, each letter $a$ in a previous word yields one letter $b$ in the next.

| Position | 0 | 1 | 2 | 3 | 4 | 5 | 6 | 7 | 8 | 9 | 10 | ... |
|---|---|---|---|---|---|---|---|---|---|---|---|---|
| Term | 0 | 1 | 1 | 2 | 3 | 5 | 8 | 13 | 21 | 34 | 55 | ... |

**Table 1.2** : The Fibonacci number sequence.

There is a remarkable relationship where any given Fibonacci word is the concatenation of the two previous Fibonacci words;

$$\mathcal{F}_n = \mathcal{F}_{n-1}\mathcal{F}_{n-2} \text{ for } n \in \mathbb{Z}, n \geqslant 3,$$

$$\text{with } \mathcal{F}_1 = ab \text{ and } \mathcal{F}_2 = aba.$$

In table 1.1, for example, notice that $\mathcal{F}_3$ = *abaab*, $\mathcal{F}_4$ = *abaababa* and that $\mathcal{F}_5$ is the concatenation $\mathcal{F}_4\,\mathcal{F}_3$ which is $\mathcal{F}_5$ = *abaababaabaab*. Respected authors, such as Lothaire, take this relationship as the definition of the Fibonacci words [Lot83, page 10]. As an aside, notice that in principle this relationship provides a computationally efficient way to generate the finite Fibonacci words up to any given *n* by repeated concatenation of two previous words.

A subword of the infinite Fibonacci word, $\mathcal{F}$, is defined to be any word that can be found as a piece of $\mathcal{F}$. Clearly, the Fibonacci words, $\mathcal{F}_n$, are subwords of $\mathcal{F}$. From studying $\mathcal{F}$ = *abaababaabaababaabaababaabaab...* it can be seen that other subwords of $\mathcal{F}$ include, for example, *bab* and *baab*.

We now define the infinite Fibonacci word's language to simply be all possible subwords of $\mathcal{F}$. Table 1.3 lists the language's subwords up to six letters length.

| length | 1 | 2 | 3 | 4 | 5 | 6 | ... |
|---|---|---|---|---|---|---|---|
| $\mathcal{F}$ words | a | aa | aab | aaba | aabaa | aabaab | ... |
| | b | ab | aba | abaa | aabab | aababa | ... |
| | | ba | baa | abab | abaab | abaaba | ... |
| | | | bab | baab | ababa | ababaa | ... |
| | | | | baba | baaba | baabaa | ... |
| | | | | | babaa | baabab | ... |
| | | | | | | babaab | ... |
| | | | | | | | ... |
| N° of words | 2 | 3 | 4 | 5 | 6 | 7 | ... |

**Table 1.3** : All subwords up to six letters in length in the Fibonacci word's language.

### 1.3  Two Key Language Exclusions

Intentionally, table 1.3 has been presented somewhat "out of the blue" and this is because, whilst it gives illuminating hints of structure within the Fibonacci word, a detailed derivation of it is not relevant to the aspects of the Fibonacci word of interest in chapters 2 and 3. However, the fact that neither of the words *bb* or *aaa* occur is of crucial importance and we now attend to a proof of this.

**Lemma 1.1 : Exclusion from $\mathcal{F}$ of *bb***

Of the four possible two letter subwords that can be formed from the letters *a* and *b* all but *bb* are subwords of the infinite Fibonacci word, $\mathcal{F}$.

*Proof*

Note that *bb* does not occur in $\mathcal{F}_0$, $\mathcal{F}_1$, $\mathcal{F}_2$ or $\mathcal{F}_3$. Inspection of $\mathcal{F}_3$ = *abaab* reveals that the subwords of length two letters in the Fibonacci word's language include *ab, ba* and *aa*. For the letter *b* to occur in a word subsequent to $\mathcal{F}_3$ it must come from the part of the Fibonacci substitution $a \rightarrow ab$. Whenever *b* occurs in $\mathcal{F}_n$ with $n \geqslant 3$ it must be preceded by an *a*. Thus, in the Fibonacci word's language, *bb* can not appear as a subword in any word of any length. □

**Lemma 1.2 : Word ending alternation of *ab* with *ba***

For $n \geq 1$, the ending of the Fibonacci words alternates between *ab* and *ba*. In fact, the Fibonacci word $\mathcal{F}_n$ ends in *ab* when $n$ is odd and ends in *ba* when $n$ is even.

*Proof*

In general, if an $n^{th}$ word ends in *ab*, then the $(n + 1)^{th}$ word must end in *aba* where the *a* on the end of the $(n + 1)^{th}$ word is generated by the substitution acting on the end *b* of the $n^{th}$ word, and the $(n + 1)^{th}$ word's penultimate *b* is generated by the substitution acting on the penultimate *a* of the $n^{th}$ word. On the other hand, if an $n^{th}$ word ends in *ba* then the $(n + 1)^{th}$ word will end in *ab*, this being generated by the substitution acting on the end *a* of the $n^{th}$ word. This proves the first part of the Lemma, that for $n \geq 1$, the endings alternate between *ab* and *ba*. Now observe that $\mathcal{F}_1 = ab$. That is, $\mathcal{F}_1$ ends in *ab*, from which it is deduced that $\mathcal{F}_n$ ends in *ab* when $n$ is odd. By further deduction, because of the alternating nature of the endings, $\mathcal{F}_n$ ends in *ba* when $n$ is even.    □

**Lemma 1.3 : Exclusion from $\mathcal{F}$ of *aaa***

Of the eight possible three letter subwords that can be formed from the letters *a* and *b* only *aab, aba, baa* and *bab* are subwords of the infinite Fibonacci word. In particular the subword *aaa* can not occur.

*Proof*

The eight possible three letter subwords referred to are, *aaa, aab, aba, baa, abb, bab, bba* and *bbb*. Of these, Lemma 1.1 excludes those with *bb* as a subword, leaving, *aaa, aab, aba, baa* and *bab*.

We will now prove that *aaa* must be removed from this shortened list.

Clearly, *aaa* does not occur in $\mathcal{F}_0 = a$ or in $\mathcal{F}_1 = ab$.

The basis for an inductive proof is established by noting that neither $\mathcal{F}_2 = aba$ nor $\mathcal{F}_3 = abaab$ contain the subword *aaa*. Recall the property that any given Fibonacci word is the concatenation of the two previous words. Our proof by induction requires we assume that for any given $\mathcal{F}_n$ with $n \geq 4$, the two previous words $\mathcal{F}_{n-1}$ and $\mathcal{F}_{n-2}$ do not contain *aaa*. In forming $\mathcal{F}_n$, the *aaa* can only occur at the seam of the concatenation of the two previous words.

Let the symbol : denote that seam.

Case 1 : The concatenation is ... *a* : *aa* ...
   This cannot happen as no word begins with *aa*, they all begin with *ab*.

Case 2 : The concatenation is ... *aa* : *a* ...
   This cannot happen as no word ends with *aa*, they all end alternatively with *ab* and *ba*, by Lemma 1.2.

By induction, the subword *aaa* cannot occur in any Fibonacci word.

Furthermore, with *aaa* now also removed from the list, what remains is,

$$aab, \ aba, \ baa, \ bab.$$

Inspection of $\mathcal{F}_6 = abaababaabaababaababa$ demonstrates that no further words should be removed from the list.    □

## 1.4 Desubstitution

Lemma 1.1 and Lemma 1.3 provide simple examples of what is more generally termed *pattern avoidance*. For the Fibonacci word it is said that *bb* and *aaa* are *illegal* subwords. In fact there are infinitely many other illegal subwords of $\mathcal{F}$. However, thankfully, the Fibonacci subwords have a desirable property that means the other illegals are not of particular interest here. This property results in it being easier to work out the legal status of a subword when the need arises rather than have long lists of which subwords are, and which are not, legal. The property is that of being readily amenable to *desubstitution*. To explain desubstitution it first needs to be understood that we will end up wanting to analyse alleged pieces of the infinite Fibonacci word. Typically, it will not be known exactly where in the word the piece to be analysed is located nor, indeed, if it is actually a legal piece of the Fibonacci word at all. Furthermore, if legal, the piece could be in multiple possible locations, maybe even in an infinite number of locations! Desubstitution provides a method of determining the legal status of the piece. It is the act of inverse substitution; working out for a given word, what previous word (or words) it could have arisen from.

Recall that the Fibonacci substitution $\theta$ is given by,

$$a \rightarrow ab$$
$$b \rightarrow a$$

where some colour has been added to assist with understanding a forthcoming example. We can now reason that when desubstituting, any letter *b* along with the letter *a* that must be to its left, must have come from a letter *a* in a preceding word. Once all of the occurrences of *ab* have been dealt with, any remaining occurrences of the letter *a* must have come from a letter *b* in the preceding word.

For example, suppose it is wishes to determine the legal status of the following alleged piece of the Fibonacci word,

$$...abaabaababaabababa...$$

This would initially be scanned for any occurrence of *bb* or *aaa* to check if it is obviously illegal. As the word passes this initial scan the next step is to bracket each occurrence of the letter *b* along with the *a* immediately to its left, like this,

$$...(ab)\,a\,(ab)\,a\,(ab)\,(ab)\,a\,(ab)\,(ab)\,(ab)\,a\,...$$

and now the desubstitution can be made which reveals the previous word to be,

$$...ababaabaaab...$$

The illegal subword *aaa* is now in the desubstituted word and so it is deduced that the alleged piece of the Fibonacci word being analysed is also illegal. More examples of desubstitution will be given in chapter 2. For now, note that if the word obtained from an initial desubstitution does not contain *aaa* the desubstitution process is repeated until either *aaa* does occur, or a sufficiently short enough word that is clearly legal is obtained. The process will, of course, eventually terminate in one of these two situations because the length of the word reduces each time a desubstitution is applied. The astute reader may be wondering why, other than in the initial scan, a lookout for the illegal *bb* is not being kept. This is because ...*aaa*... desubstitutes to either ...*bbb*... or ...*bba*... both of which contain the illegal *bb* and so, when desubstituting, the illegal subword *aaa* will always be encountered in the desubstitution before the illegal subword *bb* can occur.

### 1.5 The Incidence Matrix

In general, a substitution, $\Phi$, has associated with it a matrix termed its *incidence matrix*, $\mathbf{M}_\Phi$. For the Fibonacci substitution $\theta$, $\mathbf{M}_\theta$ is of great utility in determining information about the relative frequencies of the letters in the infinite Fibonacci word and also in suggesting a natural tiling geometry that flows from that word. We'll look first at determining relative frequencies.

Here is a standard definition of what an incidence matrix is;

For a substitution $\Phi$ acting on an alphabet $\mathcal{A} = \{ a_1, a_2, ..., a_n \}$ of cardinality $n$, the incidence matrix is defined to be the $n \times n$ matrix $\mathbf{M}_\Phi = (m_{rc})$ where $m_{rc}$ is equal to the number of occurrences of $a_r$ in $\Phi(a_c)$.         (Taken from [GY21])

For the Fibonacci substitution, $\theta$,

$$\begin{array}{l} a \to ab \\ b \to a \end{array} \quad \text{has incidence matrix} \quad \mathbf{M}_\theta = \begin{pmatrix} 1 & 1 \\ 1 & 0 \end{pmatrix}$$

In $\mathbf{M}_\theta$ the upper left $1$ indicates a single letter $a$ in $\theta(a)$ and the lower left $1$ indicates a single letter $b$ in $\theta(a)$. Similarly, the upper right 1 shows there is one letter $a$ in $\theta(b)$ and the lower right 0 that there was no letter $b$ in $\theta(b)$.

To better understand how this carries information about the letter frequencies of the Fibonacci words, consider the word $\mathcal{F}_6 = abaababaabaababaababa$ and observe that, using the notation defined previously, $|a| = 13$ and $|b| = 8$.

The connection with $\begin{pmatrix} 1 & 1 \\ 1 & 0 \end{pmatrix}^6 = \begin{pmatrix} 13 & 8 \\ 8 & 5 \end{pmatrix}$ is suggestive of the following,

---

**Proposition 1.1 : Powers of an Incidence Matrix**

If $\mathbf{M}_\Phi$ is the incidence matrix for $\Phi$ then $\mathbf{M}_\Phi^n$ is the incidence matrix for $\Phi^n$.

---

*Proof*

In general, the matrix $\mathbf{M}_\Phi$, which is of necessity square, has the form,

$$\mathbf{M}_\Phi = \begin{array}{c} a_1 \\ a_2 \\ \vdots \\ a_r \\ \vdots \\ a_n \end{array} \begin{pmatrix} m_{11} & m_{12} & \dots & m_{1c} & \dots & m_{1n} \\ m_{21} & m_{22} & \dots & m_{2c} & \dots & m_{2n} \\ \dots & \dots & \dots & \dots & \dots & \dots \\ m_{r1} & m_{r2} & \dots & m_{rc} & \dots & m_{rn} \\ \dots & \dots & \dots & \dots & \dots & \dots \\ m_{n1} & m_{n2} & \dots & m_{nc} & \dots & m_{nn} \end{pmatrix}$$

$$\Phi(a_1) \quad \Phi(a_2) \quad \dots \quad \Phi(a_c) \quad \dots \quad \Phi(a_n)$$

where $\Phi$ is the substitution on an alphabet $\mathcal{A} = \{a_1, a_2, ..., a_r, ..., a_n\}$. The scaffolding around the matrix (not routinely shown) is helping to illustrate that the $rc^{th}$ entry gives the number of times the letter $a_r$ occurs in $\Phi(a_c)$ where $a_c$ is the $c^{th}$ letter in $\mathcal{A}$.

In the square of this matrix the $rc^{th}$ entry is given by the scalar product of the $r^{th}$ row with the $c^{th}$ column;

$$\mathbf{M}_\Phi^2 (m_{rc}) = m_{r1}m_{1c} + m_{r2}m_{2c} + \dots + m_{rc}m_{rc} + \dots + m_{rn}m_{nc}$$

which will be the number of times the letter $a_r$ occurs in $\Phi^2(a_c)$.

This establishes a basis for a proof by induction for the case $n = 2$ in $\mathbf{M}_\Phi^n$. For the inductive step, assume that for $n = k$ for some integer $k \geq 3$, the $rc^{th}$ entry of $\mathbf{M}_\Phi^k$ gives the number of times the letter $a_r$ occurs in $\Phi^k(a_c)$. Taking the scalar product of the $r^{th}$ row of $\mathbf{M}_\Phi$ with the $c^{th}$ column of $\mathbf{M}_\Phi^k$ gives how often the letter $a_r$ appears in $\Phi^{k+1}(a_c)$. So, the $rc^{th}$ entry of $\mathbf{M}_\Phi^{k+1}$ represents the occurrences of the $a_r^{th}$ letter in $\Phi^{k+1}(a_c)$. By induction, $\mathbf{M}_\Phi^n$ is the incidence matrix for $\Phi^n$. □

---

The major implication of linking an incidence matrix to the Fibonacci word is that results from Linear Algebra can be utilised. In particular, the incidence matrix of the Fibonacci substitution has the desirable property of being *primitive*. By definition, a real matrix $\mathbf{M}$ is primitive if it is non-negative and its $m^{th}$ power is positive for some natural number $m$; That is, all entries of $\mathbf{M}^m$ are strictly positive for some $m \in \mathbb{Z}^+$. As it stands, the incidence matrix for the Fibonacci substitution is non-negative as required but the zero in its incidence matrix is an initial cause for concern! However, it is primitive because, for example, its square is strictly positive;

$$\begin{pmatrix} 1 & 1 \\ 1 & 0 \end{pmatrix}^2 = \begin{pmatrix} 2 & 1 \\ 1 & 1 \end{pmatrix}$$

A substitution with a primitive incidence matrix is itself termed primitive. The incidence matrix, $\mathbf{M}$, of a primitive substitution has the marvellous property of having a strictly positive simple eigenvalue, $\lambda_{PF}$. In absolute value, $\lambda_{PF}$ is strictly the largest eigenvalue of $\mathbf{M}$, a fact that is a consequence of the Perron-Frobenius theorem. See [Q95, page 132] for a proof. Additionally, $\lambda_{PF}$ has an associated eigenvector with strictly positive entries. For an infinite fixed point of a substitution, provided it is a primitive substitution, the relative frequencies with which the various letters occur exist. Furthermore, $\lambda_{PF}$ is the key to determining their values. So, what are the relative frequencies with which the letters $a$ and $b$ occur in the infinite Fibonacci word, $\mathcal{F}$?

To answer this question, note that the Fibonacci substitution's incidence matrix has characteristic polynomial $\lambda^2 - \lambda - 1 = 0$, yielding $\lambda_{PF} = \phi$, where $\phi$ denotes the golden ratio, an irrational number that has an exact value of,

$$\phi = \frac{1 + \sqrt{5}}{2}.$$

The relative letter frequencies are given by a right eigenvector $(x, y)^T$ for $\lambda_{PF}$, with sum of entries $x + y = 1$.
Solving under these constraints yields,

$$x = relative\,frequency\,(a) = \frac{1}{\lambda_{PF}} = \phi - 1 \quad \text{(About 61.8\%)}$$

$$y = relative\,frequency\,(b) = 1 - \frac{1}{\lambda_{PF}} = 2 - \phi \quad \text{(About 38.2\%)}$$

Usefully, the ratio of $|a| : |b|$ can be written as being precisely $\phi : 1$ and from this follows a straight forward proof that the Fibonacci word is non-periodic. This is presented next as Proposition 1.2.

**Proposition 1.2 : Non-Periodic Nature of the Fibonacci Word**

The Fibonacci word with letter ratio $|a|:|b|$ of $\phi:1$ is non-periodic.

The golden ratio, $\phi$, is an irrational number with an exact value of $\dfrac{1+\sqrt{5}}{2}$.

*Proof*

Assume that the Fibonacci word is strictly periodic in which case there would be some partitioning of the word along the lines of, for example,

$$(aba...ab)(aba...ab)(aba...ab)(aba...ab)\,...$$

Within a partition, let the number occurrences of the letter *a* be *p* and the number of occurrences of the letter *b* be *q* and note that from this construction both *p* and *q* must be positive integers.

The ratio of the occurrences of the two letters is thus given by $p:q$.

More revealingly, this ratio can be expressed as $\dfrac{p}{q}:1$ with $p, q \in \mathbb{Z}^+$.

A direct comparison is now made with the letter ratio $\phi:1$ of the Fibonacci word. $\dfrac{p}{q}$, is, from the definition, a rational number, and so cannot equal $\phi$.

Therefore the Fibonacci word, as claimed, must be non-periodic.   □

Proving the Fibonacci word is non-periodic is not enough to also claim that it has an aperiodic structure. To be aperiodic, we must rule out the possibility that the Fibonacci word contain an arbitrarily large periodic part. By way of explaining the need for this exclusion, consider the word that starts with five *a*, then alternates between *b* then *a* ad infinitum;

*aaaaabababababab...*

This is non-periodic but no more interesting than if it were. It is non-periodic via an annoying technicality! To show that the Fibonacci word is not of this nature it needs to be demonstrated that it does not contain arbitrarily large periodic parts, which is what would make it, by definition, *aperiodic* [Wik21a]. In 1983 Juhani Karhumäki, as part of a more general result, proved that in the Fibonacci word no subword can occur more than three times in succession [Kar83]. This is the property that makes the Fibonacci word so mathematically intriguing.

A further subtlety, that will be encountered in chapter 2, is that any given subword of $\mathcal{F}$ will occur an infinite number of times. This makes the Fibonacci word an example of a *recurrent* word [Lot02 page 31]. Table 1.4 provides examples of a recurrence, a square and a cube in the Fibonacci word.

$\mathcal{F}$ = aba abab aaba abab a abab aaba abab aaba...

$\mathcal{F}$ = a baababaa baababaa babaabaababaaba...

$\mathcal{F}$ = abaababaa baaba baaba baaba ababaaba...

**Table 1.4 :** Upper : Any subword of the Fibonacci word, such as *abab*, occurs an infinite number of times but not more than three times in succession.
Middle : The square; *baababaa* occurring twice in succession in the Fibonacci word.
Lower : The cube; *baaba* occurring thrice in succession in the Fibonacci word.

## 1.6 An Associated Tiling Geometry

The Fibonacci word can be visualised as a one-dimensional tiling. In striving to find a tiling that is faithful to the word, two lengths of tile would seem logical, tile *A* of length $L_A$ representing a letter *a* and tile *B* of length $L_B$ representing a letter *b*. The inflation mapping, $\Theta$, needs to inflate tile *A* to have a length that is the sum of the lengths of a tile *A* plus a tile *B* and also inflate tile *B* to have the length of a tile *A*. The tile length inflation geometry of $\Theta$ can be summarised in a manner that makes a correspondence with $\theta$ obvious, namely;

$$L_A \rightarrow L_A + L_B$$

$$L_B \rightarrow L_A$$

Previously, in our work of relative frequencies, it was the right eigenvector of the incidence matrix that yielded the frequency information. In general, it turns out, however, that when an associated geometry is sought, it is the left eigenvector that is required. The eigenvalue required is again $\lambda_{PF} = \phi$ because this is the only eigenvalue with a modulus greater than unity, and so able to yield the desired expansive mapping. If we arbitrarily assign a length of 1 to tile *B*, then tile *A* will be of length $\phi$. Figure 1.1 gives a summary of the resulting inflation mapping.

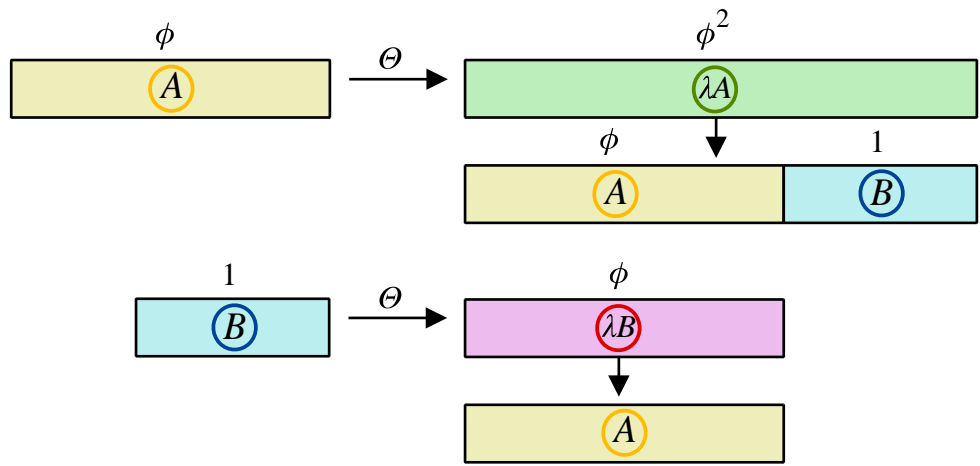

**Figure 1.1** : A summary of the inflation mapping, $\Theta$, associated with the Fibonacci word.

To conclude this chapter, figure 1.2 presents the first few Fibonacci words as one-dimensional geometric tilings of $\mathbb{R}^+$. The lengths are drawn to scale and they have been given an arbitrary width in order they be seen.

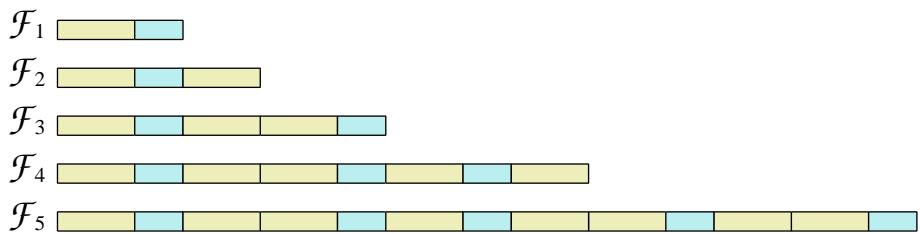

**Figure 1.2** : The first few Fibonacci words represented as tilings on a line segment.

### 1.7 Author's Notes on Chapter 1

The definition of mathematical beauty that launches the chapter is intended to be provocative, and is somewhat tongue in cheek. I'm very aware of the joy that mathematicians experience when working only with symbols, effortlessly moving around a topic area with which they are highly familiar, bending it to solve a problem or construct a proof, and without any recourse to a physical pen and paper diagram. That this beauty is hidden to those outside of mathematics prompted Samuel Eilenberg at Columbia University to once remark that, "Mathematics is a performance art, but one whose only audience is fellow performers" [PM16]. Possibly, this is a little less so since the advent of personal computers in the 1980s and the spectacular visuals associated with some branches of mathematics, in particular fractal geometry and the Mandelbrot set.

An earlier version of this chapter contained several proofs concerning the structure of the Fibonacci word; that it is Sturmian, for example. However, I felt that the overall effect was a loss of a clear direction, and so they were ruthlessly removed to sharpen the focus upon what was necessary to prepare the groundwork for the remainder of this dissertation. In a similar vein there was felt to be no need to introduce left, two-sided, or periodic fixed points.

Doorways to mathematically deep waters abound in this material. For example, when analysing words that arise as the fixed points of substitutions under iteration, just how effective desubstitution can or can't be is termed *recognizability*, and is an active area of current research. See, for example [Kyr19, Chapter 3]. Pattern avoidance as a general phenomenon is another example of a vast subject area in which active research is ongoing. See, for example [Ram04], [Ram12].

My initial understanding of the material covered in this chapter stems from The Open University's topic guide (not available outside of the University) to their course *Aperiodic Tilings and Symbolic Dynamics*, written by Reem Yassawi and the late Uwe Grimm [GY21].

I would like to thank Dan Rust of The Open University for his valuable comments on an early draft of this chapter and in particular his recommendation that I wield an editor's knife to remove that which didn't need to be there.

# Chapter 2

# Iterated Function System Drawing Rules

### 2.1  Mathematical DNA

One of the wonders of the modern age is the discovery that each living creature has the blueprint of its construction, development, function and reproduction within itself, in each of its cells in a molecule called DNA. In essence this long molecule, twisted with an antiparallel copy of itself as a double helix, can be represented as a one-dimensional string of letters. The letters represent the four different basic units of DNA, the nucleotides with a base of either thymine, **T**, adenine, **A**, guanine, **G** or cytosine, **C**. In figure 2.1 an example of a DNA segment is given. It is the sequencing of these four letters, over three billion of them in human DNA, that is the biological code at the centre of what a creature is and can become.

   **ATCCAAGCGCCCGCTAATTCTGTTCTGTTAATGTTCATACCAAGAACCGGC**
**Figure 2.1** : An example of a segment of DNA

With the foregoing in mind, it's reasonable to wonder if the Fibonacci word can be used as the code to produce another mathematical entity. Indeed, that is exactly that was achieved in chapter 1 with the Fibonacci words being viewed as tilings of line segments. Rather than that being the end of the matter, this chapter will take the idea forward with some further exploration of the idea.

In biology the concept of a string rewriting system is attributed to Aristid Lindenmayer. This Hungarian biologist sought a method to describe elementary plant development, publishing his thoughts in 1968 [Lin68]. Today referred to as *L-systems*, their grammar, *G*, is a collection of the following three parts;

- An alphabet, *V,* of symbols.
- A string of symbols *ω* defining the system's initial state.
- A set of production rules, *P*, specifying how a next string is produced from a previous string.

By repeatedly applying the rules, *P,* an *L-system* generates a sequence of strings in a similar fashion to the generation of the Fibonacci words where,

$$G = (V, \omega, P) = (\mathcal{A}, a, \theta)$$

and where, $\mathcal{A} = \{ a, b \}$, $\omega = $ "*a*" and *θ* is the substitution $a \to ab, b \to a$.

In 1986, the Polish computer scientist Przemyslaw Prusinkiewicz formalised the mathematics of *L-systems* and the interpretation of the strings as drawing instructions for a "turtle", the on-screen drawing pen controlled by a computer running the 1980s popular programming language LOGO [Prz86].

The enthusiasm at the time was to produce *fractals* following the publication in 1982 of Benoit Mandelbrot's "The Fractal Geometry of Nature" [Man82]. In this chapter the Fibonacci words are interpreted via a one dimensional drawing rule. The resulting visualisations are not thought to have been seen or explored previously. They give intriguing methods of revealing aspects of, and symmetries within, the Fibonacci words.

## 2.2 The To and Fro Drawing Rule

Essentially, the *L-system* iterative generation of the Fibonacci word was taken care of in chapter 1. Attention here can focus on drawing rules. Table 2.1 gives the first rule to be considered. It produces, from a new perspective, the tiling of section 1.6 along $\mathbb{R}^+$. This could be regarded as an identity drawing rule.

| Symbol | Action |
|:------:|:------:|
| *a* | forward $\phi$ |
| *b* | forward 1 |

**Table 2.1** : Identity drawing rule

By adding a simple about turn to the midpoint of each *b* tile the to and fro drawing rule of table 2.2 is created. This rule is believed to be previously unexplored and is what we shall focus on throughout this chapter.

| Symbol | Action |
|:------:|:------:|
| *a* | forward $\phi$ |
| *b* | forward 0.5, turn 180°, forward 0.5 |

**Table 2.2** : To and fro drawing rule.

When the to and fro drawing rule is applied to the Fibonacci word a tiled *path* results that repeatedly overwrites itself on a one-dimensional line. To see what is going on, figure 2.2 makes use of the otherwise vacant second dimension, giving a *deviation from zero diagram*. It shows the rule applied to $\mathcal{F}_6$ giving a tiled path running from *A* to *B*. Each tile has a *control point* at its midpoint, added to enhance the ease with which a reader can see what is going on. As the path descends in a zigzag fashion artificial regions are created, coloured purple. These are bounded by the path and an artificial *y*-axis (shown red) and give an impression of how much each zig and zag swings away from the origin.

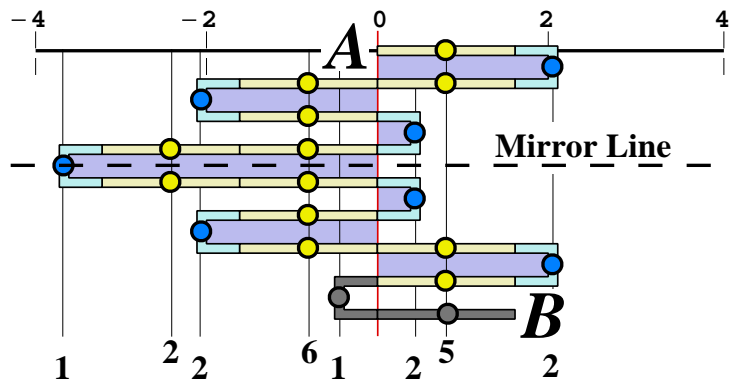

**Figure 2.2** : The $\mathcal{F}_6$ deviation from zero diagram from the to and fro drawing rule. To see the back and forth motion along a one-dimensional line the path is presented as a two-dimensional zigzag. The numbers at the foot of the diagram tally the number of tile control points at each possible value of *x* that such points occur. Originally, it was intended to explore the mathematics of these tallies but that become sidelined as interest shifted to the swings back and forth about the artificial *y*-axis. The control point tallies have been left in for a possible future analysis.

A striking feature of figure 2.2 is that, with the last two tiles greyed out, the deviation from zero diagram has reflection symmetry in the mirror shown as a broken line. Such a line must pass through the central *b* tile of a *palindromic* word. Table 2.3, below, suggests that palindromic words feature in all Fibonacci words, some with a central letter *b*, others with a central letter *a*, and some with an empty central letter. The words are coloured to highlight these observations.

| $\mathcal{F}_n = \theta^n(a)$ |
| --- |
| $\mathcal{F}_1 = ab$ |
| $\mathcal{F}_2 = aba$ |
| $\mathcal{F}_3 = abaab$ |
| $\mathcal{F}_4 = abaababa$ |
| $\mathcal{F}_5 = abaababaabaab$ |
| $\mathcal{F}_6 = abaababaabaababaababa$ |
| $\mathcal{F}_7 = abaababaabaababaababaabaababaabaab$ |

**Table 2.3** : The last two letters of each Fibonacci word shown are greyed out. What remains is a palindrome in one of three types in a repeated modulo 3 sequence of an empty central letter, then a central letter *a*, then a central letter *b*.

---

**Lemma 2.1 : Palindromic Nature of The Fibonacci Word**

Fibonacci words $\mathcal{F}_n$ with the last two letters removed are palindromic, $n \geq 1$.

---

*Proof*

Let $\mathcal{W}_n$ be the Fibonacci word $\mathcal{F}_n$ with the last two letters removed for $n \geq 1$.
Then $\mathcal{W}_1$ is the empty word, $\mathcal{W}_2 = a$ and $\mathcal{W}_3 = aba$ which are palindromes.
We thus have a basis for a proof by strong induction.
Assume when $n \leq k$ for some constant $k \geq 3$, $\mathcal{W}_k$ are palindromic.
Consider the finite Fibonacci word $\mathcal{F}_{k+1}$ and recall that,

$$\mathcal{F}_{k+1} = \mathcal{F}_k \, \mathcal{F}_{k-1}$$
$$= \mathcal{F}_{k-1} \, \mathcal{F}_{k-2} \, \mathcal{F}_{k-1}$$

The Fibonacci words, as observed previously in Lemma 1.2, end alternatively with the letters pairs *ab* and *ba*.
Case 1 : $\mathcal{F}_{k-1}$ ends in *ab* in which case $\mathcal{F}_{k-2}$ ends in *ba* and so,

$$\mathcal{F}_{k+1} = \mathcal{W}_{k-1}(ab) \, \mathcal{W}_{k-2}(ba) \, \mathcal{W}_{k-1}(ab)$$
$$\mathcal{W}_{k+1} = \mathcal{W}_{k-1}(ab) \, \mathcal{W}_{k-2}(ba) \, \mathcal{W}_{k-1} \quad \text{which is a palindrome.}$$

Case 2 : $\mathcal{F}_{k-1}$ ends in *ba* in which case $\mathcal{F}_{k-2}$ ends in *ab* and so,

$$\mathcal{F}_{k+1} = \mathcal{W}_{k-1}(ba) \, \mathcal{W}_{k-2}(ab) \, \mathcal{W}_{k-1}(ba)$$
$$\mathcal{W}_{k+1} = \mathcal{W}_{k-1}(ba) \, \mathcal{W}_{k-2}(ab) \, \mathcal{W}_{k-1} \quad \text{which is a palindrome.}$$

In both cases, $\mathcal{W}_{k+1}$ is a palindrome.
By strong induction the proof is complete.    □

**Lemma 2.2 : The Central Letter of a Palindromic Word**

Let $\mathcal{W}_n$ be the Fibonacci word $\mathcal{F}_n$ with the last two letters removed, $n \geqslant 1$.
By Lemma 2.1, $\mathcal{W}_n$ is a palindrome.
Furthermore, the central letter of this palindrome is,
- empty    if $n \equiv 1 \pmod 3$
- $a$          if $n \equiv 2 \pmod 3$
- $b$          if $n \equiv 0 \pmod 3$

*Proof*

For $\mathcal{F}_1 = ab$, $\mathcal{W}_1$ is empty and so has an empty central letter with $n \equiv 1 \pmod 3$.
For $\mathcal{F}_2 = aba$, $\mathcal{W}_2 = a$ and so has central letter $a$ with $n \equiv 2 \pmod 3$.
For $\mathcal{F}_3 = abaab$, $\mathcal{W}_3 = aba$ and so has central letter $b$ with $n \equiv 0 \pmod 3$.
The above three statements establish a basis for an inductive proof.
Assume when $n \leqslant k$ for some constant $k \geqslant 3$, $\mathcal{W}_k$ has an empty central letter when $k \equiv 1 \pmod 3$, a central letter $a$ when $k \equiv 2 \pmod 3$ and a central letter $b$ when $k \equiv 0 \pmod 3$.
Consider the construction used in Lemma 2.1, namely,

$$\begin{aligned}
\mathcal{F}_{k+1} &= \mathcal{F}_k \, \mathcal{F}_{k-1} \\
&= \mathcal{F}_{k-1} \, \mathcal{F}_{k-2} \, \mathcal{F}_{k-1} \\
&= \mathcal{W}_{k-1}(uv) \, \mathcal{W}_{k-2}(vu) \, \mathcal{W}_{k-3}(ab) \\
\mathcal{W}_{k+1} &= \mathcal{W}_{k-1}(uv) \, \mathcal{W}_{k-2}(vu) \, \mathcal{W}_{k-1}
\end{aligned}$$

where $uv$ is a letter for letter replacement for one of $ab$ or $ba$.
This shows that the palindrome $\mathcal{W}_{k+1}$ has the same central letter as $\mathcal{W}_{k-2}$
Strong induction now completes the proof.  □

From a broader perspective, Lemma 2.2 can be viewed as a consequence of long established theorems from the number theory of Fibonacci numbers. One relevant result is that the Fibonacci number sequence modulo $n$ is periodic. The length of the period modulo $n$ is denoted $\pi(n)$. It is termed the $n^{\text{th}}$ *Pisano period* [Wik21c]. The fact proven in Lemma 2.2, that the central letter of $\mathcal{W}_n$ repeatedly cycles through three possibilities, stems from $\pi(2) = 3$.
An implication of Lemma 2.2 is that, after $\mathcal{F}_6$, the next Fibonacci word with mirror symmetry through a *B* tile is $\mathcal{F}_9$ and figure 2.3 suggests this is so, subject to checking the "in between" deviation from zero diagrams for $\mathcal{F}_7$ and $\mathcal{F}_8$.
A little thought, as captured by figure 2.4, indicates that in order to be,
- palindromic
- without the illegal double *B* or triple *A* tiling combination

they will have to have half turn rotational symmetry about either the empty central tile or the central *A* tile (and thus, in both cases, not have mirror symmetry). Figure 2.5 shows $\mathcal{F}_7$ where the half turn rotational symmetry is centred on the *zero deviation* line. This is the line of the artificial *y*-axis representing the origin of $\mathbb{R}^+$. Figure 2.6 shows $\mathcal{F}_8$ and again, as expected, it has half turn rotational symmetry but this time not centred on the zero deviation line.

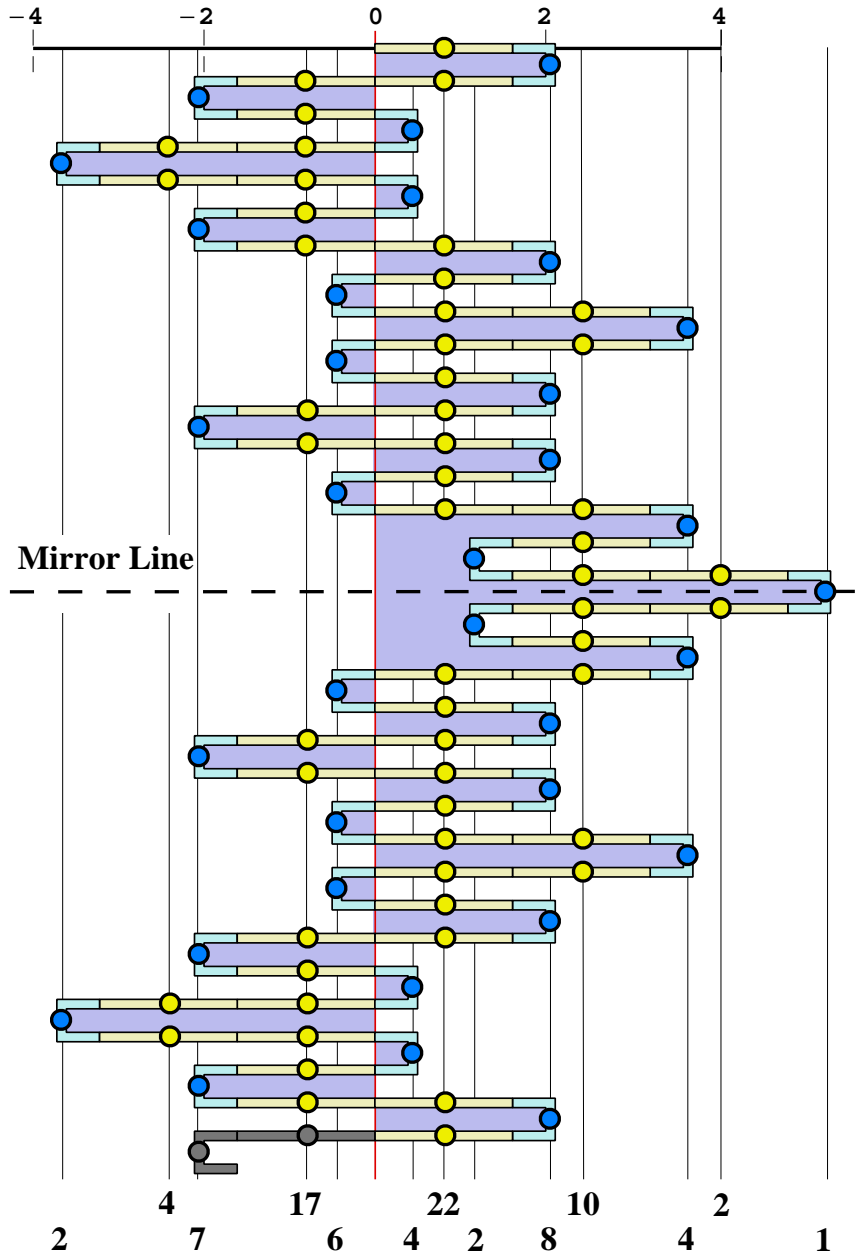

**Figure 2.3** : The deviation from zero diagram for $\mathcal{F}_9$, under the to and fro drawing rule. Like that for $\mathcal{F}_6$ (figure 2.2) it has a line of mirror symmetry through a central letter $b$ of its corresponding palindromic word ($\mathcal{W}_9$).

| $n \;(mod\; 3) \equiv 0$ | $n \;(mod\; 3) \equiv 1$ | $n \;(mod\; 3) \equiv 2$ |
|---|---|---|

**Figure 2.4** : The possible configurations at the centre of the palindromic words, $\mathcal{W}_n$, associated with the Fibonacci words, $\mathcal{F}_n$.

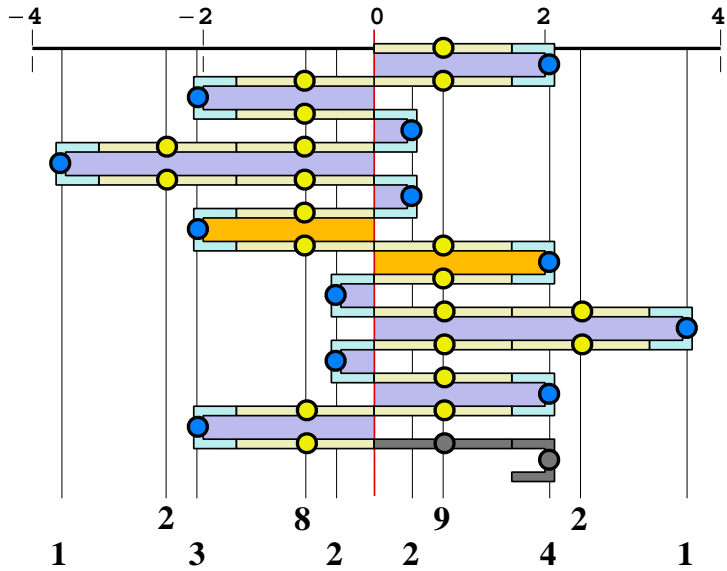

**Figure 2.5** : The deviation from zero diagram for $\mathcal{F}_7$ under the to and fro drawing rule with the $\mathcal{W}_7$ centre, corresponding to an $n \,(mod\, 3) \equiv 1$ case (figure 2.4) highlighted in orange. $\mathcal{W}_7$ has half-turn rotational symmetry about its centre which is on the zero deviation line (shown red).

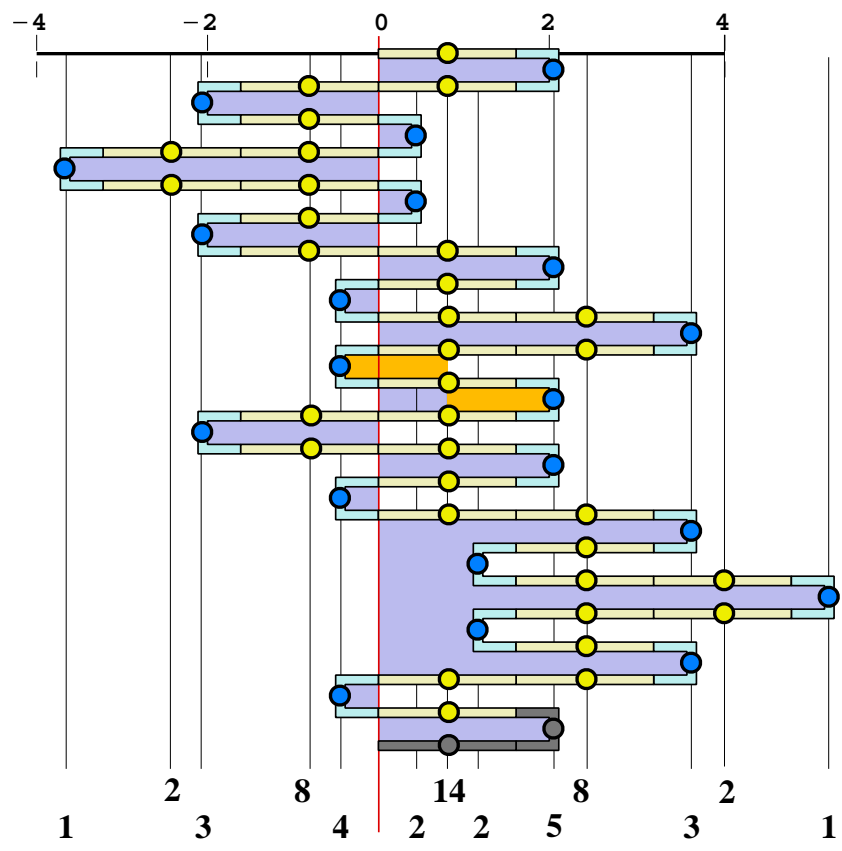

**Figure 2.6** : The deviation from zero diagram for $\mathcal{F}_8$ under the to and fro drawing rule with the $\mathcal{W}_8$ centre, corresponding to an $n \,(mod\, 3) \equiv 2$ case, (figure 2.4) highlighted in orange. $\mathcal{W}_8$ has half-turn rotational symmetry about its centre but, unlike figure 2.5 for $\mathcal{F}_7$, this centre is not on the zero deviation line (shown red).

## 2.3  About The Zero Deviation Line

The to and fro drawing rule was so named because of the manner in which the one-dimensional path swings repeatedly from one side of the origin to the other then back along $\mathbb{R}$. It's natural to wonder how far away from the origin it can swing and what are the characteristics of the possible deviations away from the origin. To get a feel for such "bigger picture" questions, we need to look at longer Fibonacci words as, with iteration in general, emergent behaviours are not necessarily obvious from the first few iterations of the system. The deviation from zero diagram for $\mathcal{F}_{11}$ (a "longer word") is given in figure 2.7 where the individual tiles and control points are no longer shown, only the areas enclosed by them and the zero deviation line.

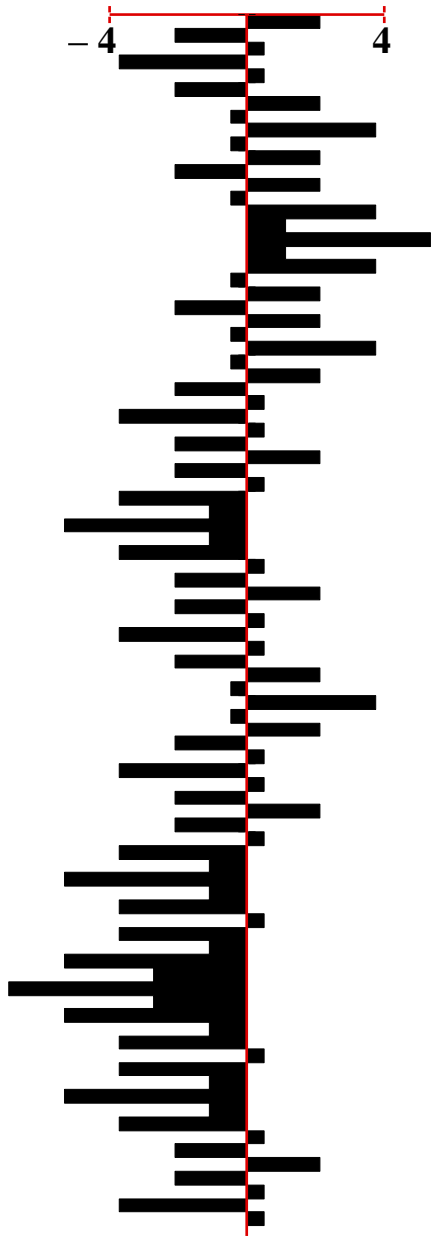

**Figure 2.7** : The deviation from zero diagram for $\mathcal{F}_{11}$ under the to and fro drawing rule. This was drawn by a computer running the LOGO programming language. Samples of the code written are presented in Appendix A.

At first sight, figure 2.7 looks complicated but once the seemingly different structures within it are isolated, there are remarkably few. These are given in our next diagram, figure 2.8, and each assigned a catalogue identifier.

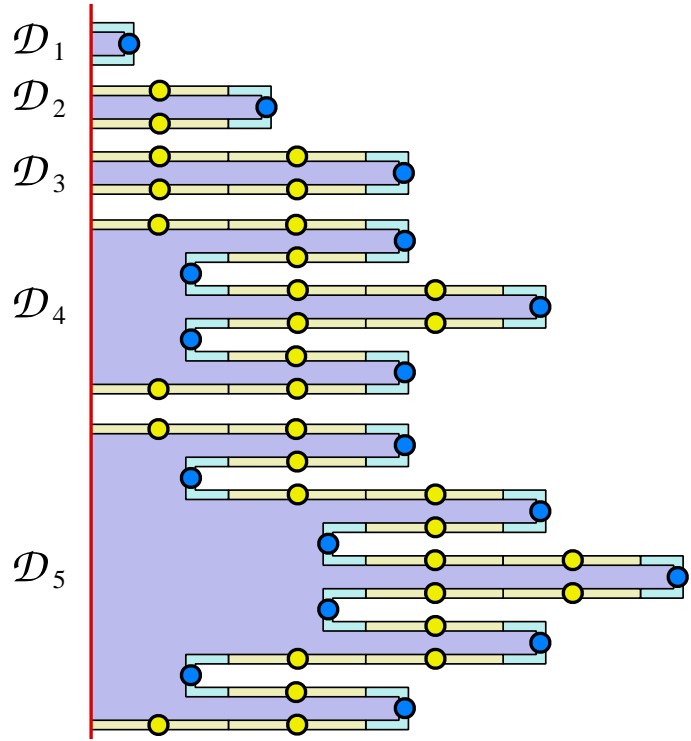

**Figure 2.8** : All five of the seemingly different structures evident in the deviation from zero diagram for $\mathcal{F}_{11}$ under the to and fro drawing rule.

One could be forgiven for expecting more complicated structures to start to occur as the length of the Fibonacci word increases. Pondering figure 2.8, however, suggests an alternative viewpoint; what is emerging is an ever larger piece of a single structure.

The realisation there may only be a single fundamental structure underpinning all of the Fibonacci words rather begs the question of why some other structures have not appeared and, indeed, why they cannot do so as the word length is increased further. Figure 2.9 gives an example of a simple structure that has not occurred in any of the deviation from zero diagrams up to and including $\mathcal{F}_{11}$. So, will this specific structure eventually pop up in a Fibonacci word, one not yet investigated, or can it never occur?

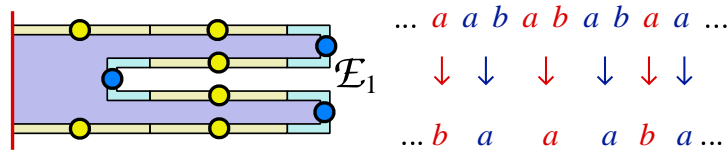

**Figure 2.9** : Desubstitution shows that $\mathcal{E}_1$ is an excluded structure that can never occur.

The key to answering the question is to desubstitute the part of the word associated with the structure. This desubstitution process was described earlier in section 1.4. Here it gives *...baaaba...* as the previous word. The triple *a* is illegal and the inescapable conclusion is that this structure can never occur.

A piece of word associated with a structure, or a proposed structure, will be referred to as an *embedded word*. This is because, when desubstitution is applied, it may be necessary to speculate upon what the letters immediately before or after the embedded word might be. It does not necessarily have a fixed start or finish. A desubstitution may need to be applied repeatedly to determine the legal status of an embedded word. For example the embedded word ...*aababaabaabaaba*... after two successive desubstitutions is shown to be illegal because it would have had to come from a word with the illegal triple *a* in it. Figure 2.10 provides the details.

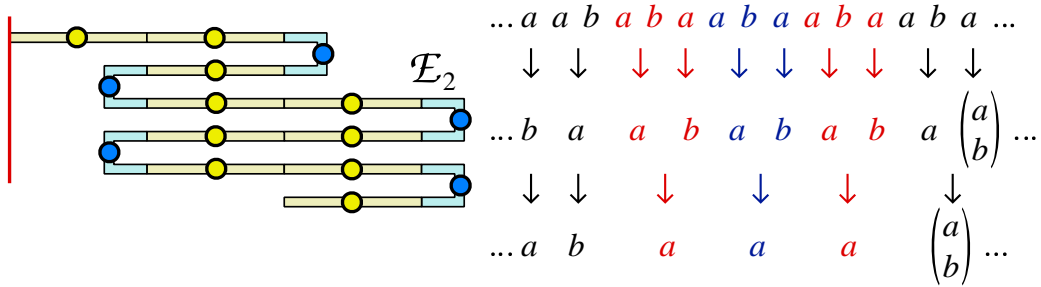

**Figure 2.10** : It takes two successive desubstitutions to show this structure can never occur. The brackets $\begin{pmatrix} a \\ b \end{pmatrix}$ are used to show that the letter immediately to the right of the embedded word could alter the letter at that end of the desubstituted word; one of an *a* or *b* could be in that location.

Figure 2.11 gives an example where it takes three successive desubstitutions to show that the embedded word ...*aababaababaababaaba*... is illegal.

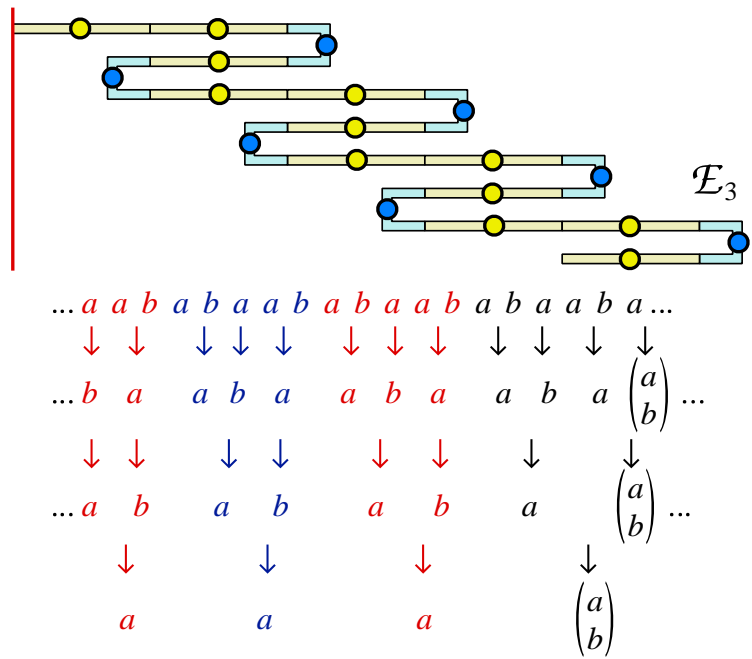

**Figure 2.11** : It takes three successive desubstitutions to show this structure can never occur.

Desubstitution provides a mechanism to test proposed structures for their legality. When applied, its effectiveness when dealing with the Fibonacci word stems from the fact that the possible desubstitutions are often unique or only uncertain at one end. This may not be the case with other substitutions.

## 2.4 Charting Growth

The desubstitution method described in section 2.3 can be thought of as a "top down" method of analysing a proposed structure's legality. However, it's a "stabbing in the dark" method of trying to find structures. They need to be guessed first, then their validity checked via desubstitution. A more systematic approach to producing a catalogue of structures is to work "bottom up" and produce a *growth chart* of all the ways in which structures can legally be constructed. Figure 2.12a shows the start of a methodical, thorough and robust search for all possible deviation from zero structures.

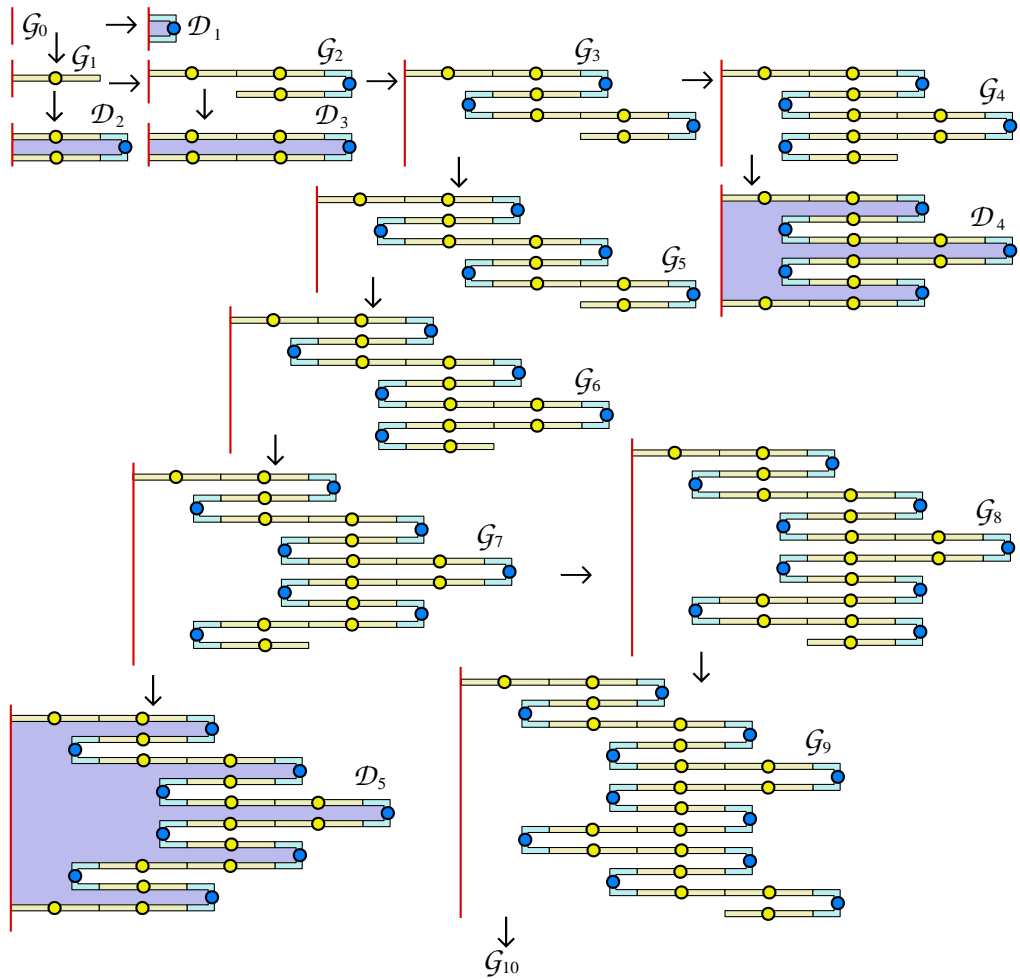

**Figure 2.12a** : The start of the growth chart, mapping all possible deviation from zero structures of less than thirty tiles.

Rather than add one letter and its associated tile each time, the growth chart speeds the process up by adding either three or five. This relies on the fact that the Fibonacci word can be factorised in the following manner;

$$\mathcal{F} = (aba)(aba)(baaba)(aba)(baaba)(baaba)(aba)\ldots$$

where the only factors involved are *aba* or *baaba*. Observe that in the growth chart when moving from one diagram to the next either the tiles *ABA* or *BAABA* has been added to a previous tiling path (once past the initial terms). The factorisation guarantees that one of *ABA* or *BAABA* can always be added. Occasionally both are valid which gives branching, at $\mathcal{G}_3$ for example.

The growth chart shows that the structures already observed in figure 2.8 (that is, $\mathcal{D}_1$, $\mathcal{D}_2$, $\mathcal{D}_3$, $\mathcal{D}_4$ and $\mathcal{D}_5$) are the only structures of less than thirty tiles that can exist. As the arrow at the foot of figure 2.12a indicates, the growth chart continues. Before presenting the extension, it was felt to be prudent to prove the claimed factorisation of the Fibonacci word. There is considerable literature on the topic of factorising words similar to the Fibonacci word. For example, Amy Glen's 2006 PhD thesis is devoted entirely to various decompositions of words (with similar properties to the Fibonacci word) into factors, with an emphasis on palindromic factorisations [Gle06]. However, for our modest needs, at this stage, a less general and more targeted proof will suffice.

---

**Theorem 2.1 : The *aba*, *baaba* Factorisation of $\mathcal{F}$**

The Fibonacci word can be factorised into the form,

$$\mathcal{F} = (aba)(aba)(baaba)(aba)(baaba)(baaba)(aba)\ldots$$

where the only factors involved are *aba* or *baaba*.

---

*Proof*

First observe that all occurrences of (*baaba*) can be further factorised as (*ba*)(*aba*) and that this is an unambiguous reversible piece of algebraic manipulation within the context of the Fibonacci word where the leftmost factor is an *aba*. Thus the claim is equivalent to proving that,

$$\mathcal{F} = (aba)(aba)(ba)(aba)(aba)(ba)(aba)(ba)(aba)(aba)\ldots$$

where the only factors involved are *aba* or *ba*.

All Fibonacci words, $\mathcal{F}_n$ begin with *aba* for $n \geq 2$.

Let $\mathcal{F}_n^*$ be $\mathcal{F}_n$ with that leading *aba* removed for $n \geq 2$.

Inductive base: $\mathcal{F}_1 = (ab)$, $\mathcal{F}_2 = (aba)$ and $\mathcal{F}_4 = (aba)(aba)(ba)$.

Assume when $n \leq k$ for some even constant $k \geq 4$, $\mathcal{F}_k$ (and hence $\mathcal{F}_k^*$) factorise as claimed.

Case 1 : $\mathcal{F}_{EVEN}$

$$\begin{aligned}
\mathcal{F}_{k+2} &= \mathcal{F}_{k+1}\,\mathcal{F}_k && k \text{ is even} \\
&= \mathcal{F}_k\,\mathcal{F}_{k-1}\,\mathcal{F}_k \\
&= \mathcal{F}_k\,\mathcal{F}_{k-2}\ldots\mathcal{F}_6\,\mathcal{F}_4\,\mathcal{F}_2\,\mathcal{F}_1\,\mathcal{F}_k \\
&= \mathcal{F}_k\,\mathcal{F}_{k-2}\ldots\mathcal{F}_2\,(ab)(aba)\,\mathcal{F}_k^* \\
&= \mathcal{F}_k\,\mathcal{F}_{k-2}\ldots\mathcal{F}_2\,(aba)(ba)\,\mathcal{F}_k^*
\end{aligned}$$

Case 2 : $\mathcal{F}_{ODD}$

$$\begin{aligned}
\mathcal{F}_{k+1} &= \mathcal{F}_k\,\mathcal{F}_{k-1} && k \text{ is still even} \\
&= \mathcal{F}_k\,\mathcal{F}_{k-2}\,\mathcal{F}_{k-4}\ldots\mathcal{F}_6\,\mathcal{F}_4\,\mathcal{F}_2\,\mathcal{F}_1 \\
&= \mathcal{F}_k\,\mathcal{F}_{k-2}\,\mathcal{F}_{k-4}\ldots\mathcal{F}_2\,(ab)
\end{aligned}$$

The result for $\mathcal{F}_{EVEN}$ now follows by strong induction and for $\mathcal{F}_{ODD}$ the words factorise as required except for a factor of (*ab*) at the end of each word. The infinite Fibonacci word is the word formed as $n \to \infty$. and so, in the limit, all words factorise as claimed. The proof of Theorem 2.1 is complete.    □

When working with an embedded word, there is ambiguity over whether the extension to the right of the word, which will always end with *aba*, is followed with an *aba* or *baaba* and so desubstitution is deployed each time the word is extended to determine which is legal. In fact, on occasion both *aba* and *baaba* can be legal extensions to the word, and hence give two valid extensions to the associated tiling path. This gives rise to the branching phenomenon observed in the growth chart. However, by Theorem 2.1. at least one will always be legal. Any illegal case always leads, via desubstitution, to the illegal *aaa*.

Everything is now in place to extend the growth chart started in figure 2.12a, and this extension is presented over the next few pages.

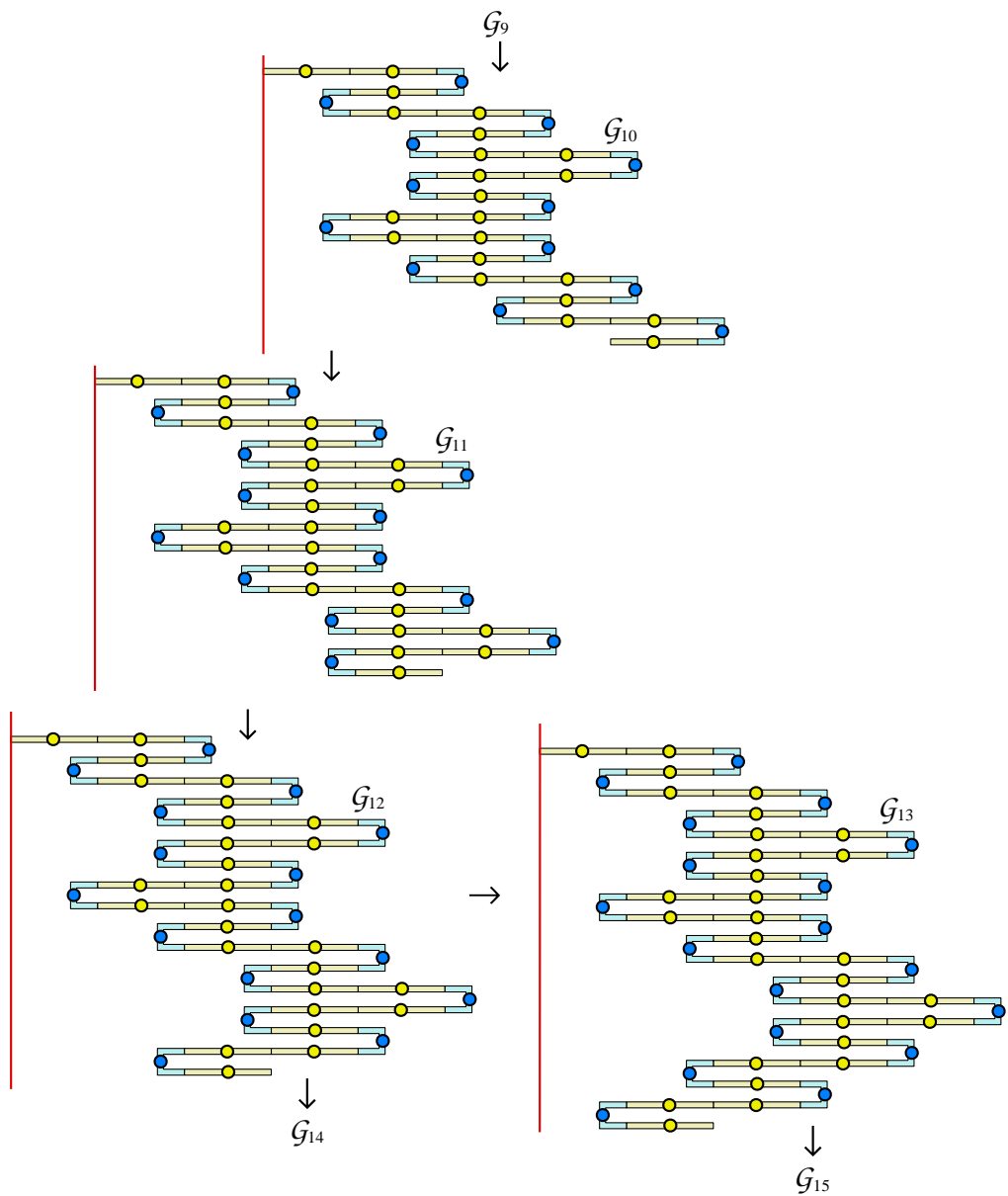

**Figure 2.12b** : The continuation of the growth chart started in figure 2.12a.
        A further continuation is in figure 2.12c.

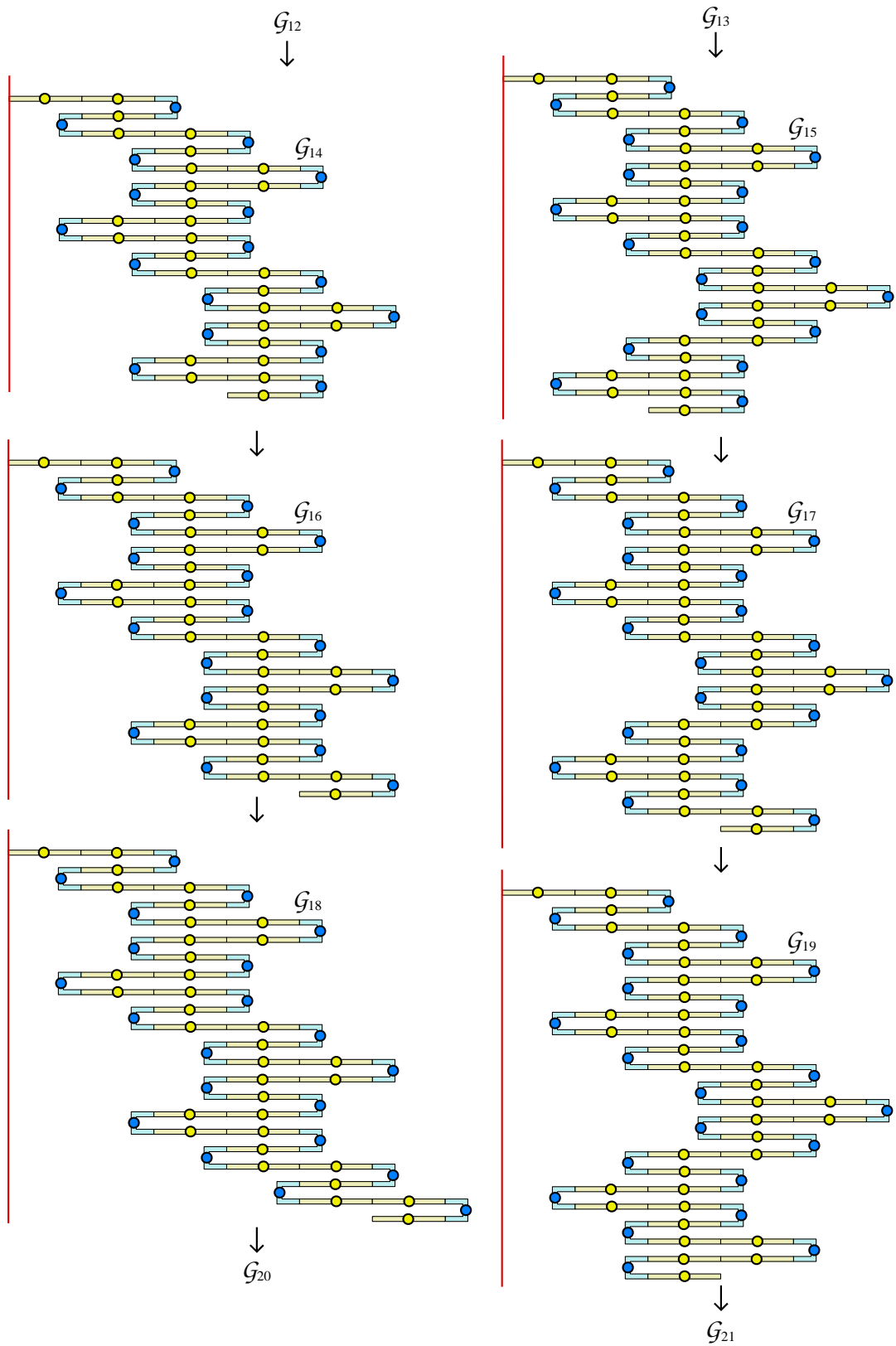

**Figure 2.12c** : The continuation of the growth chart from figure 2.12b.
The next continuation is in figure 2.12d.

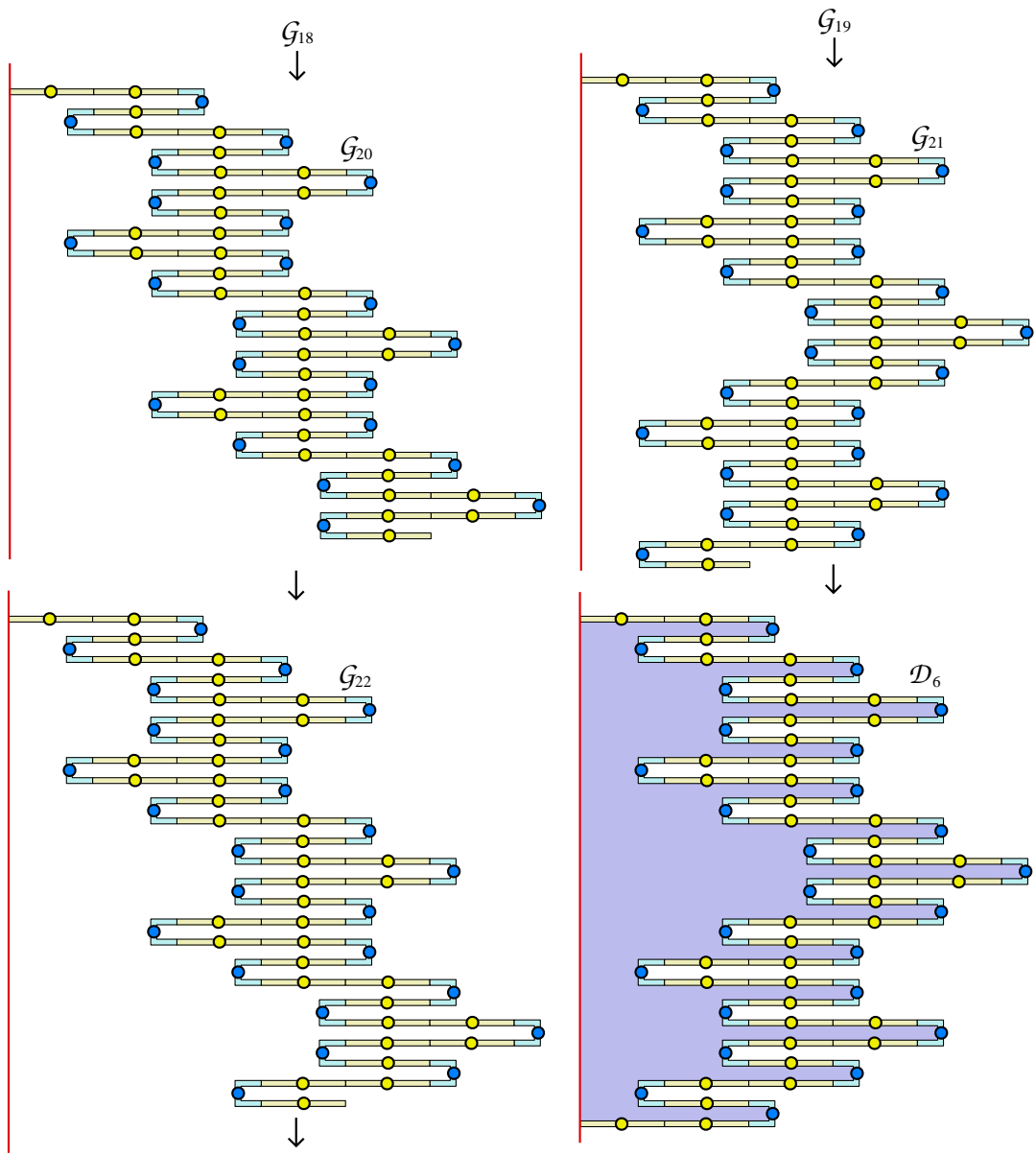

**Figure 2.12d** : The continuation of the growth chart from figure 2.12c.
It shows the deviation from zero structure, $\mathcal{D}_6$.

The growth chart presented across the four figures 2.12a, 2.12b, 2,12c and 2.12d has determined all possible deviation from zero structures up to the 67 tiles of $\mathcal{D}_6$ and the 64 tiles of $\mathcal{G}_{22}$. The growth chart could be continued further. In fact, it would seem likely that the aperiodic nature of the Fibonacci word guarantees it will continue indefinitely. All the structures found continue to suggest that they are simply increasingly large chunks of an underpinning single structure.

## 2.5 Exploring the $n \equiv 0 \pmod{3}$ Case

Two facts that give rise to some speculative thoughts are, firstly, that for the Fibonacci words $\mathcal{F}_n$ with $n \equiv 0 \pmod{3}$, the associated palindromic words $\mathcal{W}_n$ have mirror symmetry about a central $B$ tile. Secondly, as all words start on the zero deviation line, the mirror symmetry guarantees that $\mathcal{W}_n$ with $n \equiv 0 \pmod{3}$ will end on it too. In consequence, with the tiling path thus pegged to the zero deviation line at either end, it is reasonable to postulate that the maximum deviation away from the zero deviation line is predisposed to occur somewhere towards the centre of the word. As noted in figure 2.4, $\mathcal{W}_n$ with $n \equiv 0 \pmod{3}$ has a central tiling configuration of *ABA*. Theorem 2.1, the *aba, baaba* factorisation of $\mathcal{F}$ theorem, can now be applied. Some care is needed, as the *ABA* matches up to an *aba* in $\mathcal{F}$ that is at the right hand end of a factor (*ba**aba*). The aim now is to produce a growth chart, moving outward from the central *aba* and figure 2.13 is the result.

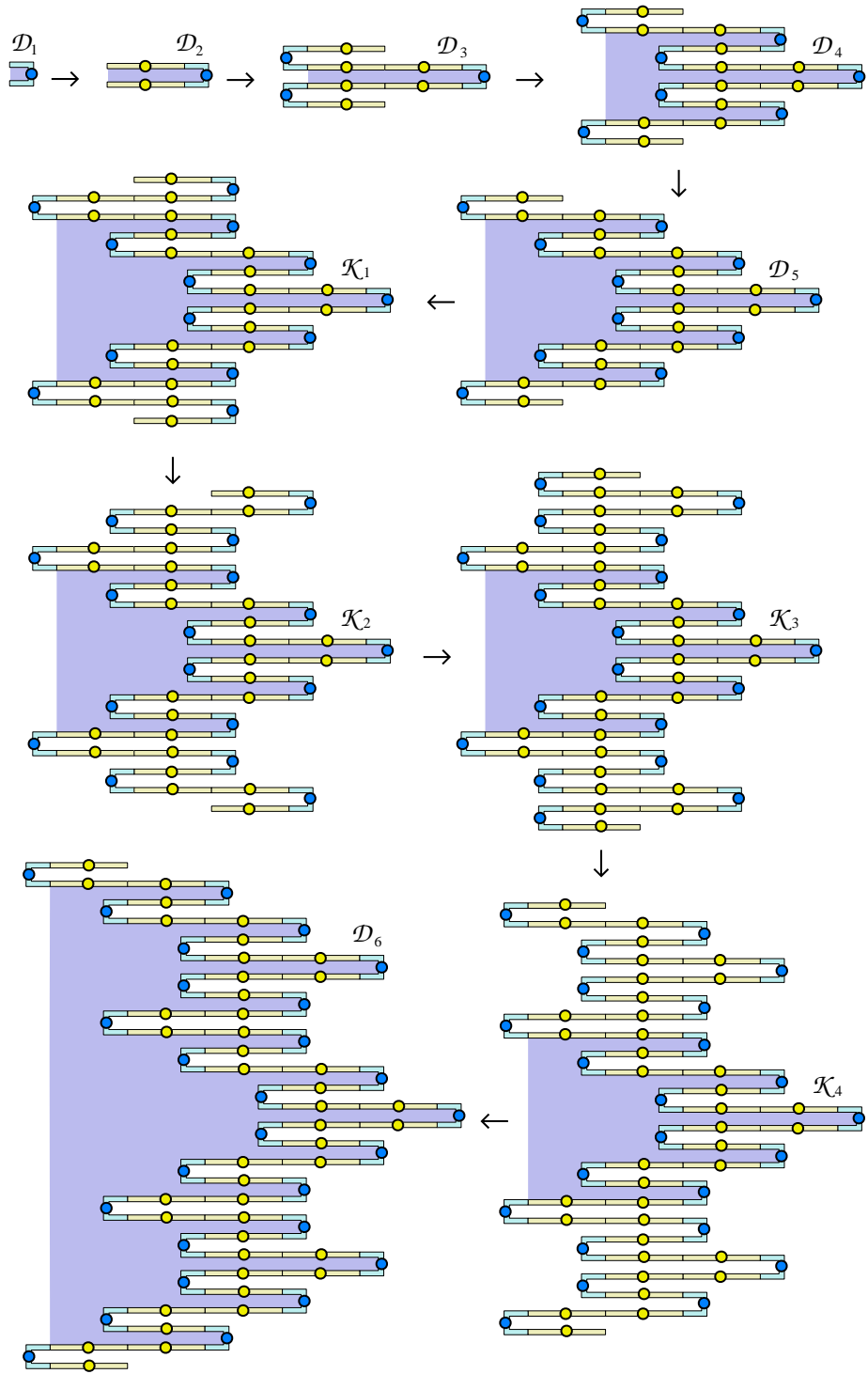

**Figure 2.13** : The tiling growth chart, moving outward from a central *ABA*, where it is reasonable to speculate the maximum deviation away from the zero deviation line occurs.

From figure 2.13 it is immediately obvious that this is, by far, a most efficient way of determining the possible deviation from zero structures that can occur. In producing figure 2.13, the tilings were first worked out in words. In practice, at each step, it was easiest to work out if the next factor to the right should be (*aba*) or (*baaba*) and then add the letters to the left that kept the required mirror symmetry. The bracketing of the factors when extending to the left is not always immediately obvious until a subsequent step has been taken. Desubstitution is used to ensure the legality of each step of the growth.

<div style="text-align:center">

*b*

*aba*

*a*)(*baaba*)(*aba*)

*a*)(*baaba*)(*baaba*)(*aba*)(*baaba*)

*aba*)(*aba*)(*baaba*)(*baaba*)(*aba*)(*baaba*)(*baaba*)

*a*)(*baaba*)(*aba*)(*baaba*)(*baaba*)(*aba*)(*baaba*)(*baaba*)(*aba*)

*aba*)(*aba*)(*baaba*)(*aba*)(*baaba*)(*baaba*)(*aba*)(*baaba*)(*baaba*)(*aba*)(*baaba*)

*a*)(*baaba*)(*aba*)(*baaba*)(*aba*)(*baaba*)(*baaba*)(*aba*)(*baaba*)(*baaba*)(*aba*)(*baaba*)(*aba*)

</div>

**Figure 2.14** : The first few stages of the tiling growth chart (of figure 2.13) moving outward from a central *ABA*, but in words rather than tiles. How the words extend was determined first, then the growth chart tiling drawn. This was because desubstitution needed to be repeatedly used to check which of the two possible extensions to the right was the legal one.

From how the growth chart has been constructed, or simply from looking at the tilings (figure 2.13) or the words (figure 2.14) it is clear that what is being "crawled along" with the extension to one end and its palindromic reflection to the other, is a piece of the fixed point of the iteration. This is not a surprise once it is noticed that,

$$\mathcal{F}_{3m} = \mathcal{F}_{3m-1}\mathcal{F}_{3m-2} \qquad m \in \mathbb{Z}, \; m \geq 1$$

$$= \mathcal{F}_{3m-2}\mathcal{F}_{3m-3}\mathcal{F}_{3m-2}.$$

This shows that a previous word of the form $\mathcal{F}_n$ for $n \equiv 3 \pmod 3$ is contained centrally in a subsequent word of that same form. In fact, because of the iterative nature of the process, it is contained centrally in **all** subsequent words of that same form. Figure 2.15 depicts some specific examples; $\mathcal{F}_6 = \mathcal{F}_4\mathcal{F}_3\mathcal{F}_4$, $\mathcal{F}_9 = \mathcal{F}_7\mathcal{F}_6\mathcal{F}_7$ and $\mathcal{F}_{12} = \mathcal{F}_{10}\mathcal{F}_9\mathcal{F}_{10}$. Colour and superposition are used to show the relationship between them. As figure 2.15 suggests, these separate results can be nested and we can write, for example, that,

$$\mathcal{F}_{12} = \mathcal{F}_{10}\mathcal{F}_7\mathcal{F}_4\mathcal{F}_3\mathcal{F}_4\mathcal{F}_7\mathcal{F}_{10}$$

This discussion readily generalises and is succinctly summarised by Lemma 2.3 as a formal result which will prove useful in the next section.

---

**Lemma 2.3 : Centralised Recursive Embedding of $\mathcal{F}_3$**

$$\mathcal{F}_{3m} = \mathcal{F}_{3m-2}\mathcal{F}_{3m-5}\ldots\mathcal{F}_7\mathcal{F}_4\left(\mathcal{F}_3\right)\mathcal{F}_4\mathcal{F}_7\ldots\mathcal{F}_{3m-5}\mathcal{F}_{3m-2} \qquad m \in \mathbb{Z}^+$$

*Proof*

As this is a straight forward proof by induction, it is left as an exercise for the interested reader.   □

---

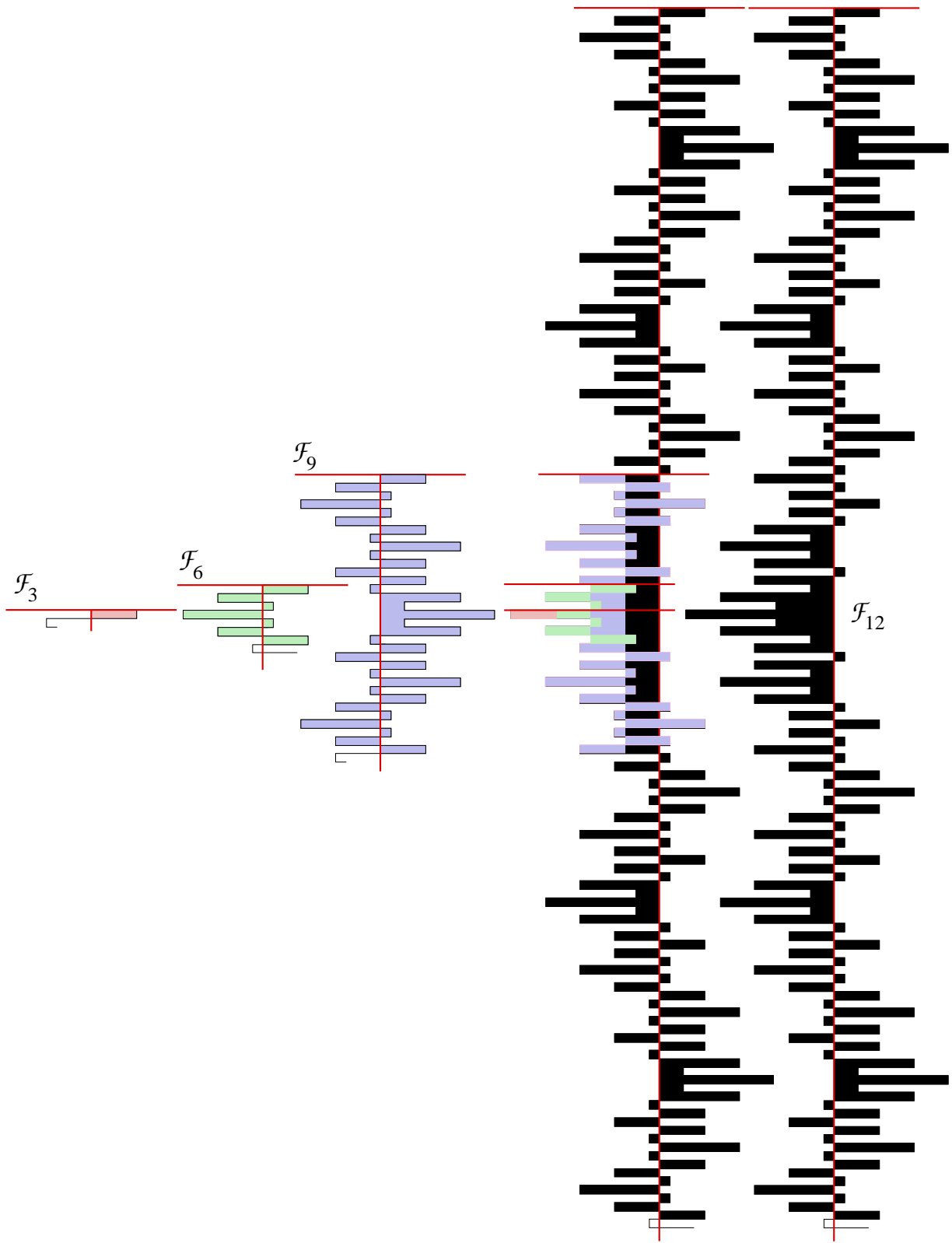

**Figure 2.15** : An illustration of how the shorter words appear at the centre of subsequent words when working only with words $\mathcal{F}_n$ where $n \equiv 0 \pmod{3}$. Sometimes the tiling path associated with a word needs reflecting before being overlayed. When such reflections are needed will be explained later. Colour is used to show how in $\mathcal{F}_{12}$ the previous words $\mathcal{F}_3$, $\mathcal{F}_6$ and $\mathcal{F}_9$ occur. One could say that $\mathcal{F}_3$ is embedded in $\mathcal{F}_6$ which is embedded in $\mathcal{F}_9$ which is embedded in $\mathcal{F}_{12}$. This figure was drawn by a computer running the LOGO programming language (Appendix A).

## 2.6  How Far From the Origin ?

For this final section of chapter 2, we return to a question raised earlier regarding how far from the origin the zigzag path can go. From a broad perspective, given that the Fibonacci word is aperiodic, it can be argued that, given any distance from the origin, the path will eventually exceed that distance. This is because being aperiodic means there can be no extended periods of repetition, no arbitrarily large periodic parts. This is seemingly at odds with the deviation from zero diagrams which have revealed a path having a nature that seems somewhat hesitant in moving away from the zero deviation line. The question is thus morphing into one about what mechanisms cause the wandering to, overall, have a maximum deviation away from the line that moves ever outward as the $n$ in $\mathcal{F}_n$ tends towards infinity. One key idea is to determine if a Fibonacci word represents a *direction reversal* or *direction sustain*. In other words, does it set up a subsequent piece of path to initially progress in the same direction or the opposite direction? Figure 2.16 gives two examples one with a direction sustain and one with a direction reversal.

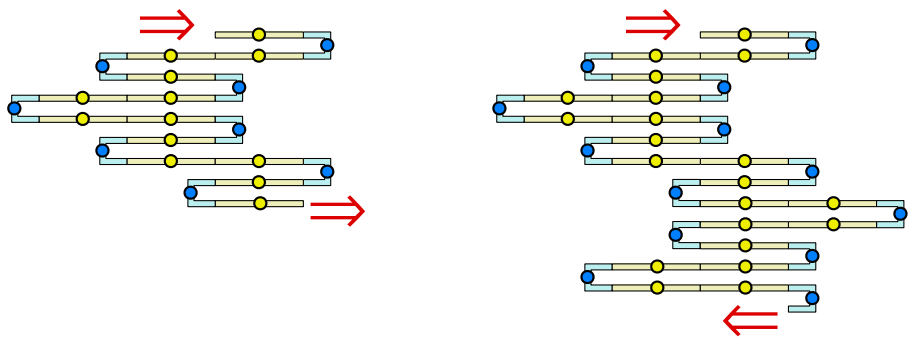

**Figure 2.16** :
To the left is $\mathcal{F}_6$ which sets up the next piece of path to flow in the same direction.
To the right is $\mathcal{F}_7$ which sets up the next piece of path to flow in the opposite direction.

Clearly, a word will represent a direction reversal if there are an odd number of the "about turn" instruction. That is, an odd number of letter $b$ in the word. It has already been observed (in table 1.1) that the number of letter $b$ in a word is in the Fibonacci number sequence and so, by considering this modulo 2, it is revealed if any given Fibonacci word represents a direction sustain or direction reversal. Table 2.2 shows the pattern that emerges.

| Word | $\mathcal{F}_1$ | $\mathcal{F}_2$ | $\mathcal{F}_3$ | $\mathcal{F}_4$ | $\mathcal{F}_5$ | $\mathcal{F}_6$ | $\mathcal{F}_7$ |
|---|---|---|---|---|---|---|---|
| N° of $b$ | 1 | 1 | 2 | 3 | 5 | 8 | 13 |
| N° of $b$ ( mod 2 ) | 1 | 1 | 0 | 1 | 1 | 0 | 1 |
| S or R | R | R | S | R | R | S | R |

**Table 2.2** : Fibonacci words are one of either direction sustain, *S*, or direction reversal, *R*.

As the number of $b$ in a word is the sum of the number of $b$ in the two previous words, the simple rules for combining odd and even numbers under addition guarantees that $\mathcal{F}_n$ is a direction sustain word when $n \equiv 0 \pmod{3}$ and a direction reversal word otherwise.

The second key idea is to determine the overall one-dimensional vector displacement, $\overrightarrow{\mathcal{F}_n}$, between the start and finish of the path associated with $\mathcal{F}_n$ relative to its initial direction. Recall, an *A* tile has length $\phi$, the golden ratio.

For a word $\mathcal{F}_n$ with $n \equiv 0 \pmod 3$, the associated tiling path for $\mathcal{W}_n$ has to both start and finish on the zero deviation line. This is because it's palindromic with mirror symmetry about a central letter *b*. The fact that words with $n \equiv 0 \pmod 2$ end in *ba* and words with $n \equiv 1 \pmod 2$ end in *ab* give us, in combination, that,

$$\overrightarrow{\mathcal{F}_n} = \begin{cases} +\phi, \text{ Sustain for } n \equiv 0 \pmod 6 \\ -\phi, \text{ Sustain for } n \equiv 3 \pmod 6 \end{cases} \text{ for } n \equiv 0 \pmod 3$$

Figure 2.17 shows the four possible situations that can arise.

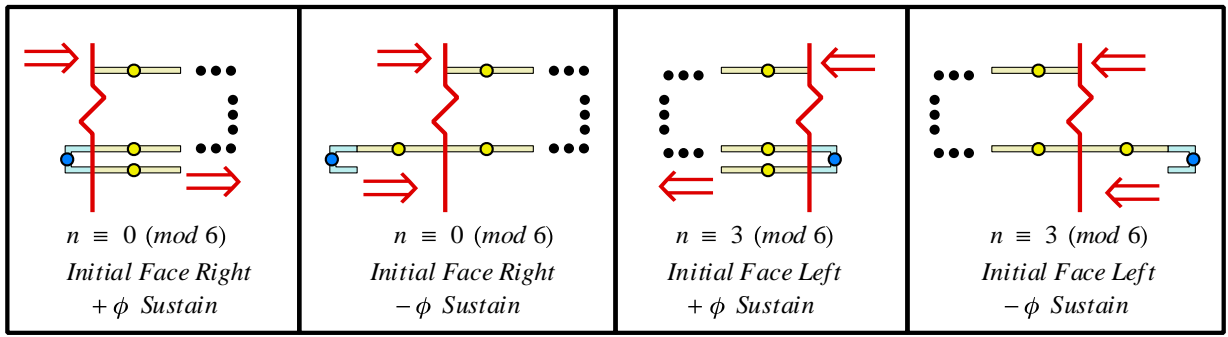

**Figure 2.17** : The four possible displacements relative to the initial direction for $n \equiv 0 \pmod 3$.

For a word $\mathcal{F}_n$ with $n \equiv 1 \pmod 3$, the associated tiling path for $\mathcal{W}_n$ has to both start and finish on the zero deviation line (because it's palindromic with half turn rotational symmetry, and has centre on the zero deviation line about a central empty word flanked on either side by an *a*). This time, drawing on the fact that words with $n \equiv 0 \pmod 2$ end in *ba* and words with $n \equiv 1 \pmod 2$ end in *ab* give us, in combination, that,

$$\overrightarrow{\mathcal{F}_n} = \begin{cases} +\phi, \text{ Reversal for } n \equiv 1 \pmod 6 \\ -\phi, \text{ Reversal for } n \equiv 4 \pmod 6 \end{cases} \text{ for } n \equiv 0 \pmod 3$$

Figure 2.18 shows the four possible situations that can arise.

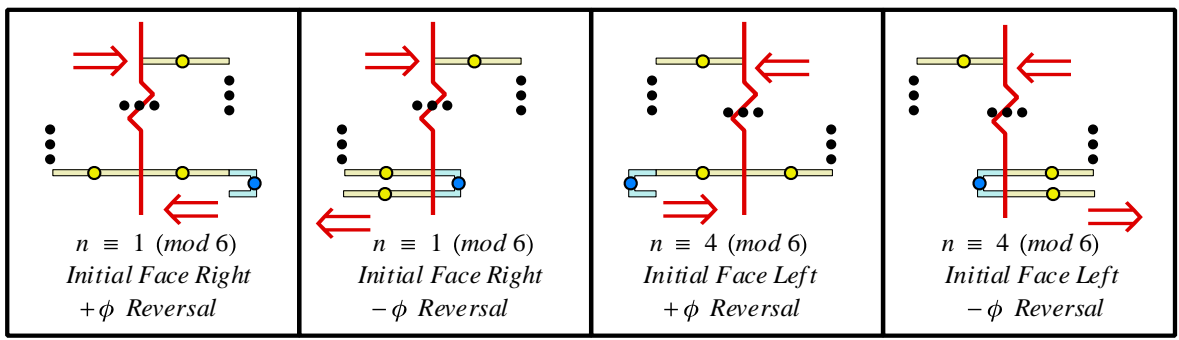

**Figure 2.18** : The four possible displacements relative to the initial direction for $n \equiv 1 \pmod 3$ .

All of the pieces are in place to work through an example that keeps track of a specific tiling's progress away from the zero deviation line by looking at a special decomposition of the word associated with that tiling.

The focus is upon words of the form $\mathcal{F}_n$ with $n \equiv 0 \pmod 3$ to which Lemma 2.3 on the centralised recursive embedding of $\mathcal{F}_3$ can be applied.

As an example consider $\mathcal{F}_{15}$ which, by Lemma 2.3, decomposes as,

$$\mathcal{F}_{15} = \mathcal{F}_{13}\, \mathcal{F}_{10}\, \mathcal{F}_{7}\, \mathcal{F}_{4}\, \mathcal{F}_{3}\, \mathcal{F}_{4}\, \mathcal{F}_{7}\, \mathcal{F}_{10}\, \mathcal{F}_{13}$$

Notice that all the words about the central $\mathcal{F}_3$ are of the form $n \equiv 1 \pmod 3$. Now to track how far from the zero deviation line the end of each of the tilings associated with the individual component words lie.

Table 2.3 presents the calculation.

| Component of $\mathcal{F}_{15}$ | $\mathcal{F}_{13}$ | $\mathcal{F}_{10}$ | $\mathcal{F}_{7}$ | $\mathcal{F}_{4}$ | $\mathcal{F}_{3}$ | $\mathcal{F}_{4}$ | $\mathcal{F}_{7}$ | $\mathcal{F}_{10}$ | $\mathcal{F}_{13}$ |
|---|---|---|---|---|---|---|---|---|---|
| Initial direction | $\Rightarrow$ | $\Leftarrow$ | $\Rightarrow$ | $\Leftarrow$ | $\Rightarrow$ | $\Rightarrow$ | $\Leftarrow$ | $\Rightarrow$ | $\Leftarrow$ |
| Relative displacement | $+\phi$ | $-\phi$ | $+\phi$ | $-\phi$ | $-\phi$ | $-\phi$ | $+\phi$ | $-\phi$ | $+\phi$ |
| Final direction | $\Leftarrow$ | $\Rightarrow$ | $\Leftarrow$ | $\Rightarrow$ | $\Rightarrow$ | $\Leftarrow$ | $\Rightarrow$ | $\Leftarrow$ | $\Rightarrow$ |
| Overall displacement | $+\phi$ | $+2\phi$ | $+3\phi$ | $+4\phi$ | $+3\phi$ | $+2\phi$ | $+\phi$ | $0$ | $-\phi$ |

**Table 2.3** : Tracking how far away from the zero deviation line the end of each component part lies. The relative displacements are relative to the initial direction.

As a check on the table 2.3 calculation, notice that, overall, $\overrightarrow{\mathcal{F}_{15}}$ is determined to be $(-\phi,\ \text{Sustain})$, which is as it should be. Such tracking calculations are at the heart of Theorem 2.2, the culminating highlight of this chapter.

---

**Theorem 2.2 : Infinite Wandering Of $\mathcal{F}_n$ From The Zero Deviation Line**

Given a distance, for convenience expressed as an integer multiple of $\phi$, the to and fro tiling path associated with the Fibonacci words $\mathcal{F}_n$ will eventually exceed that (one dimensional) distance away from the origin for sufficiently large enough $n$ where $n \equiv 0 \pmod 3$.

---

*Proof*

Two results to be used are, first, as noted above,

$$\overrightarrow{\mathcal{F}_n} = \begin{cases} +\phi, & \text{Sustain for } n \equiv 0 \pmod 6 \\ +\phi, & \text{Reversal for } n \equiv 1 \pmod 6 \\ -\phi, & \text{Sustain for } n \equiv 3 \pmod 6 \\ -\phi, & \text{Reversal for } n \equiv 4 \pmod 6 \end{cases}$$

and secondly, from Lemma 2.3,

$$\mathcal{F}_{3m} = \mathcal{F}_{3m-2}\, \mathcal{F}_{3m-5}\, \ldots\, \mathcal{F}_{7}\, \mathcal{F}_{4}\, (\mathcal{F}_3)\, \mathcal{F}_{4}\, \mathcal{F}_{7}\, \ldots\, \mathcal{F}_{3m-5}\, \mathcal{F}_{3m-2} \qquad m \in \mathbb{Z}^{+}$$

The proof uses these two results in each of two cases, one corresponding to wandering off to the right when $3m \equiv n \equiv 3 \pmod 6$ and the other corresponding to wandering off to the left when $3m \equiv n \equiv 0 \pmod 6$.

Case 1 : Wandering off to the right.
For $3m \equiv n \equiv 3 \pmod 6$ the generalised calculation will be as table 2.4 shows.

| Component of $\mathcal{F}_{3m}$ | $\mathcal{F}_{3m-2}$ | $\mathcal{F}_{3m-5}$ | ... | $\mathcal{F}_4$ | $\mathcal{F}_3$ | $\mathcal{F}_4$ | ... | $\mathcal{F}_{3m-5}$ | $\mathcal{F}_{3m-2}$ |
|---|---|---|---|---|---|---|---|---|---|
| Initial direction | ⇒ | ⇐ | ... | ⇐ | ⇒ | ⇒ | ... | ⇒ | ⇐ |
| Relative displacement | $+\phi$ | $-\phi$ | ... | $-\phi$ | $-\phi$ | $-\phi$ | ... | $-\phi$ | $+\phi$ |
| Final direction | ⇐ | ⇒ | ... | ⇒ | ⇒ | ⇐ | ... | ⇐ | ⇒ |
| Overall displacement | $+\phi$ | $+2\phi$ | ... | $(m-1)\phi$ | $(m-2)\phi$ | $(m-3)\phi$ | ... | $0$ | $-\phi$ |

**Table 2.4** : Generalised accountancy for how far away from the zero deviation line the end of each component part lies for $3m \equiv n \equiv 3 \pmod 6$.

Case 2 : Wandering off to the left.
For $3m \equiv n \equiv 0 \pmod 6$ the generalised calculation will be as table 2.5 shows.

| Component of $\mathcal{F}_{3m}$ | $\mathcal{F}_{3m-2}$ | $\mathcal{F}_{3m-5}$ | ... | $\mathcal{F}_4$ | $\mathcal{F}_3$ | $\mathcal{F}_4$ | ... | $\mathcal{F}_{3m-5}$ | $\mathcal{F}_{3m-2}$ |
|---|---|---|---|---|---|---|---|---|---|
| Initial direction | ⇒ | ⇐ | ... | ⇒ | ⇐ | ⇐ | ... | ⇒ | ⇐ |
| Relative displacement | $-\phi$ | $+\phi$ | ... | $-\phi$ | $-\phi$ | $-\phi$ | ... | $+\phi$ | $-\phi$ |
| Final direction | ⇐ | ⇒ | ... | ⇐ | ⇐ | ⇒ | ... | ⇐ | ⇒ |
| Overall displacement | $-\phi$ | $-\phi$ | ... | $(1-m)\phi$ | $(2-m)\phi$ | $(3-m)\phi$ | ... | $0$ | $+\phi$ |

**Table 2.5** : Generalised accountancy for how far away from the zero deviation line the end of each component part lies for $3m \equiv n \equiv 0 \pmod 6$.

Let $d\phi$ be an integer multiple of $\phi$.
This is then the fixed distance away from the origin that is to be exceeded.
Let $m = d + 1$.
Consider the word $\mathcal{F}_{3m}$, decomposed according to Lemma 2.3, in which the start of the tiling associated with the centrally embedded $\mathcal{F}_3$ will be at a distance $d\phi$ from the origin and orientated to project a further $\phi + 0.5$ units away from the origin. The tiling associated with $\mathcal{F}_{3m}$ thus exceeds the distance $d\phi$ from the origin.  □

Interestingly, Theorem 2.2 would yield the same conclusion it the length of the $A$ tiles where reduced from $\phi$ to 1.

### 2.7  Author's Notes on Chapter 2

Chapter 2 presents, as far as I know, some original mathematical research and has striven to both introduce an innovative visualisation of the Fibonacci word and prove new results in response to what was observed. The majority of the techniques involved are mainstream symbolic dynamical combinatorics on words manipulations. Thus, it was felt appropriate to see where the new ideas went, and show the techniques being applied in an unfamiliar setting.

Throughout this chapter there is a fascinating interplay of local rules giving rise to global behaviours and constraints. The Fibonacci word is often described in the literature as "an aperiodic word of minimal complexity". Even so, it was fascinating to realise that a single structure was emerging in spite of initially setting out with the expectation of finding and cataloguing a variety. In retrospect, I have seen this phenomenon elsewhere, when a carefully chosen symmetric piece of a tiling, a *patch,* provides an elegant way of producing, via the iteration of a substitution, a tiling of the plane, $\mathbb{R}^2$. The key aspect is that the patch is associated with a selected part of a fixed point. In this chapter that selected part is termed an embedded word. By way of a specific example of a similar situation occurring elsewhere, see figure 2.19, and for the full details of this "Chair Substitution" see [Rob99].

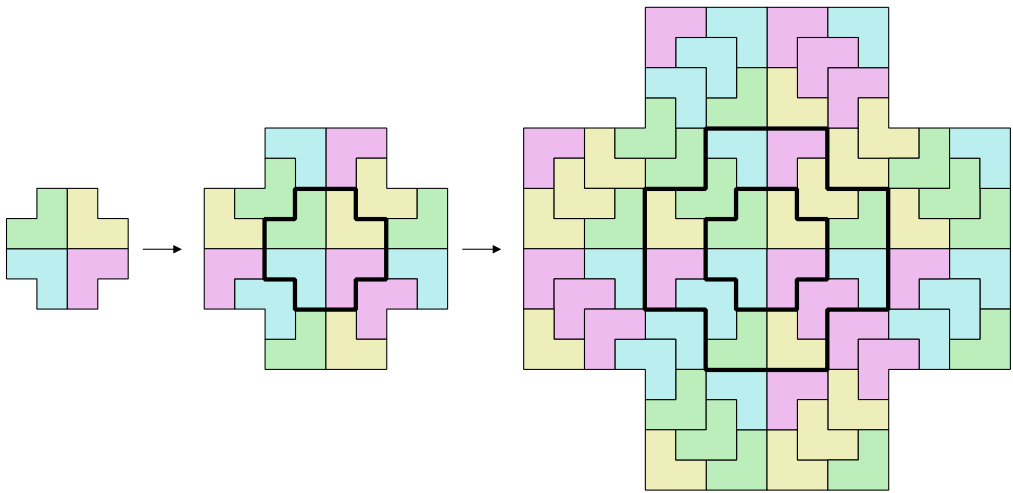

**Figure 2.19** : To the left is the initial patch associated with a selected embedded word from the fixed point of an iterated substitution. By iteration of this symmetrical start (ignoring the tile colourings) a particularly attractive demonstration of the substitutions structure is derived in much the same way as was done with the Fibonacci word in this chapter.

I originally set out in this chapter to explore a two dimensional drawing rule. However, the concept of a drawing rule is in itself so interesting that I started thinking about what the simplest possible drawing rule would be, even if it did not yield a two dimensional path. The to and fro drawing rule evolved to be amongst the most simple, and yet proved able to capture some of the symmetries of the Fibonacci word. The deviations from zero diagrams came about in response to discoveries made. I originally hit upon the idea of shading the area between the path and the artificial *y*-axis by accident but it was seized upon because of the strong visual impact of the result. The established practice is to put the control points at an end of a tile. I opted to placing it in the middle and then bent the *B* tile in half around it out of a desire to make more obvious, and then capture, the symmetry of the deviations from zero.

# Chapter 3

# Fractal In Nature

### 3.1  Overview

Previously we saw that the to and fro drawing rule on the Fibonacci word is a one-dimensional traversing, back and forth, along the line $\mathbb{R}$. In this chapter we continue to explore the nature of the Fibonacci word. A new drawing rule is applied that gives a two dimensional path across $\mathbb{R}^2$. It turns out to be a fractal.

### 3.2  Similarity Fractals

Benoit Mandelbrot's seminal work from 1982, "The Fractal Geometry of Nature", contains the example of the mathematically idealised fractal reproduced below in figure 3.1 [Man82, page 50]. Such fractals are created from two geometric objects that Mandelbrot terms an *initiator* and a *generator*. In figure 3.1 the initiator is the straight line segment, $A_1$, of unit length. The generator is the $A_2$ zigzag. The $A_2$ generator is composed of eight line segments each of length 0.25 units where the long double length segment is considered to be two singles. Each of these eight line segments can now be viewed as an initiator, each a quarter the length of the original initiator, $A_1$. The generator, scaled to fit, is applied to each and the result is $A_3$. Repeating the process, this time with a scaling factor of $0.25^2$ yields $A_4$. As the number of iterations tends to infinity, the fractal is the "curve" that results. The key idea is to have a generator that is constructed only from scaled and rotated copies of an initiator.

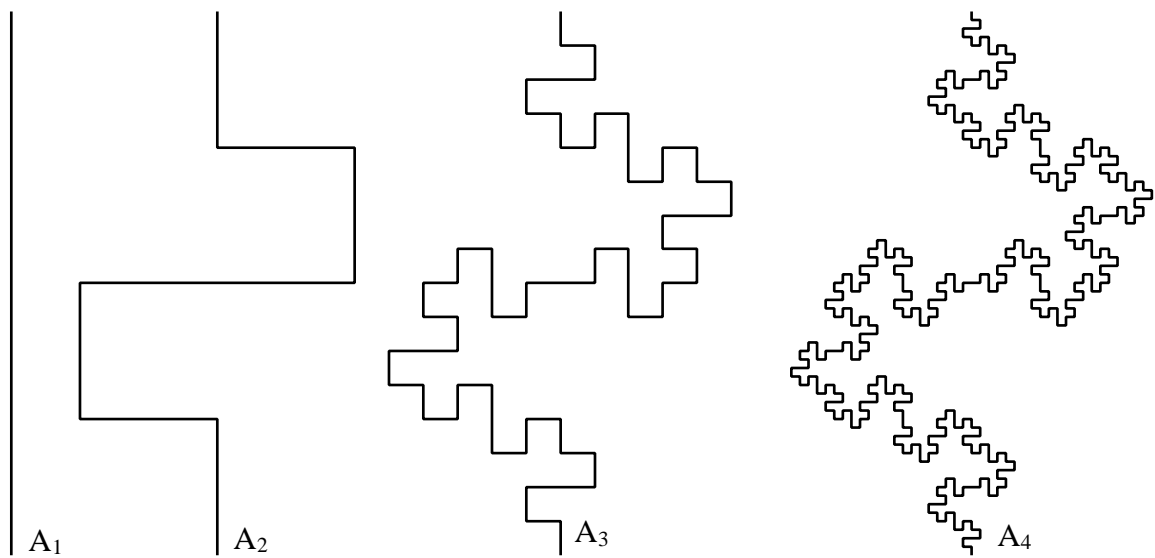

$A_1$            $A_2$            $A_3$            $A_4$

**Figure 3.1** : Geometric iteration giving rise to a fractal.

In general, a set made up of *m* copies of itself scaled by length scale factor *r* has a similarity dimension given by, $dim_{SYM} = -\dfrac{log\ m}{log\ r}$. For the figure 3.1 fractal this gives a value of 1.5 for the similarity dimension.

### 3.3 The Double Letter Drawing Rule

The introductory example was selected because of its passing resemblance to a meandering path. In striving to obtain a path of this nature from the Fibonacci word it would seem likely that, in addition to moving onward, turns to both the left and to the right will increase the likelihood of obtaining a path that does not overwrite itself excessively. The Fibonacci word's language contains three words of length two letters and this suggests a drawing rule such as that shown in table 3.1 where each possible two letter combination is assigned a different drawing instruction. Each instruction is assigned an overall length of 1 unit.

| Symbol | Action |
| --- | --- |
| *ab* | forward 0.5, turn right 90°, forward 0.5 |
| *aa* | forward 1 |
| *ba* | forward 0.5, turn left 90°, forward 0.5 |

**Table 3.1** : The double letter drawing rule.

The double letter drawing rule is only valid for Fibonacci words with an even number of letters. This is so the letters can be bracketed in pairs without a letter remaining unpartnered. This requirement causes the focus to fall upon words of the form $\mathcal{F}_n$ with $n \equiv 1 \pmod 3$. As examples satisfying these constraints let's consider $\mathcal{F}_4 = abaababa$ and $\mathcal{F}_7 = abaababaabaababaabaababaababaaba$.
These bracket in letter pairs as,

$$\mathcal{F}_4 = (ab)(aa)(ba)(ba)$$

and

$$\mathcal{F}_7 = (ab)(aa)(ba)(ba)(ab)(aa)(ba)(ba)(ab)(ab)(aa)(ba)(ab)(ab)(aa)(ba)(ab)$$

from which figure 3.2 is derived by applying the double letter drawing rule.

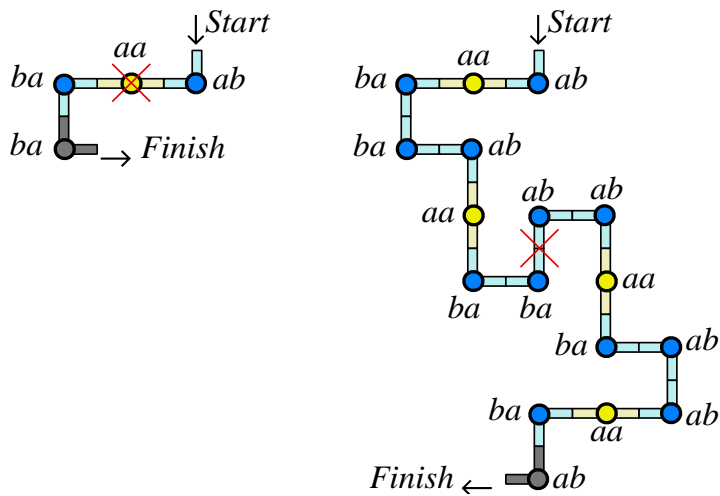

**Figure 3.2** : Left  : The tiling path for $\mathcal{F}_4$ under the double letter drawing rule.
         Right : The tiling path for $\mathcal{F}_7$ under the double letter drawing rule.
In both cases the path with the last tile removed has rotational symmetry about the point marked with a red cross. With that last tile removed, direction is preserved in the sense that the entrance movement is downward as is the exit movement.

Of course, from Lemma 2.1, we know that removing the last pair of letters from $\mathcal{F}_n$ leaves a palindromic word, $\mathcal{W}_n$. For the two above examples this gives,

$$\mathcal{W}_4 = (ab)\,(aa)\,(ba)$$

and

$$\mathcal{W}_7 = (ab)\,(aa)\,(ba)\,(ba)\,(ab)\,(aa)\,(ba)\,(ba)\,(ab)\,(ab)\,(aa)\,(ba)\,(ab)\,(ab)\,(aa)\,(ba)$$

where the middle section of each word is highlighted in red. In fact there can be no middle section other than these two possibilities because of the illegality of $bb$ and the requirement that $\mathcal{W}_n$ be palindromic. Before moving on from figure 3.2 there are a couple of observations to extract. Firstly, notice that the last tile of each of the two paths is greyed out and in both cases the coloured path that remains has half turn rotational symmetry about the point marked with a red cross. This red cross corresponds to either the tile associated with ...$(aa)$... or the tiles associated with ...$(ba)(ab)$... at the middle of $\mathcal{W}_n$. Secondly, notice that the direction along a tile is reversed in its image under the half turn rotational symmetry. To further explore and explain this symmetry, table 3.2 shows all possible orientations of the tiles drawn by the drawing rule along with the image of each under a half turn rotation with the direction "flipped". This flip is to keep the direction along the tiled path correct. This direction element cannot be ignored as the tile for a right turn $(ab)$, if passed through the "wrong way" would be seen as a left turn $(ba)$, for example. Fundamentally, the fact that a word is read from left to right needs to be preserved in the tiling path that it gives rise to.

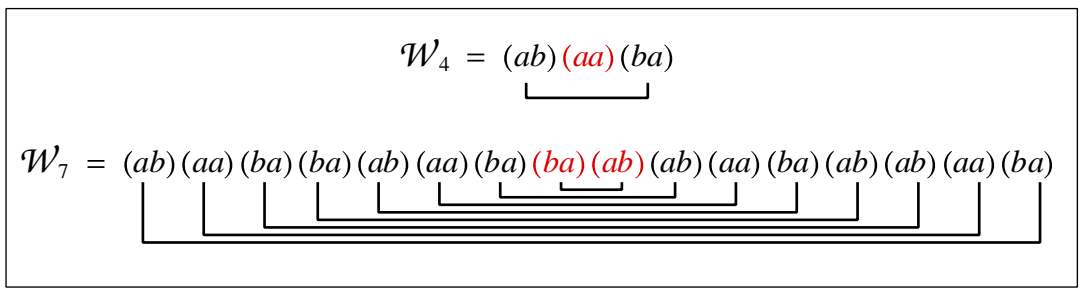

**Table 3.2** : All possible orientations of tiles by the double letter drawing rule. A tile selected from the upper row of the table has image under a half turn rotation given directly underneath, (and vice versa). Once rotated, a tile has the direction of the path through it reversed.

In table 3.2 it is noted that in all cases, when the direction is also considered, $(ab)$ has image $(ba)$, $(aa)$ is its own image, and $(ba)$ has image $(ab)$. The same matching is evident in $\mathcal{W}_4$ and $\mathcal{W}_7$ due to the palindromic nature of these words. Figure 3.3 demonstrates this;

$$\mathcal{W}_4 = (ab)\,(aa)\,(ba)$$

$$\mathcal{W}_7 = (ab)\,(aa)\,(ba)\,(ba)\,(ab)\,(aa)\,(ba)\,(ba)\,(ab)\,(ab)\,(aa)\,(ba)\,(ab)\,(ab)\,(aa)\,(ba)$$

**Figure 3.3** : Matching of $(ab)$ with $(ba)$, of $(aa)$ with itself, and of $(ba)$ with $(ab)$ in $\mathcal{W}_4$ and $\mathcal{W}_7$.

**Lemma 3.1 : Half Turn Rotational Symmetry of $\mathcal{W}_n$ when $n \equiv 1$ (mod 3)**

For $n \equiv 1$ (mod 3), $\mathcal{W}_n$ under the double letter drawing rule gives rise to a path that has half turn rotational symmetry. Furthermore, the entrance direction is the same as the exit direction.

*Proof*

Case 1 : $n \equiv 1$ (mod 6).

The number of letters in $\mathcal{F}_n$ will be of the form $4k + 2$ for some non-negative integer value $k$. With the last two (rightmost) letters removed what remains, $\mathcal{W}_n$, is a palindrome (by Lemma 2.1) and of length $4k$. As the subword $bb$ is not legal, at the middle of this palindrome must be the letter configuration ...$(ba)(ab)$... which corresponds geometrically under the double letter drawing rule to one of the tiling configurations shown in figure 3.3;

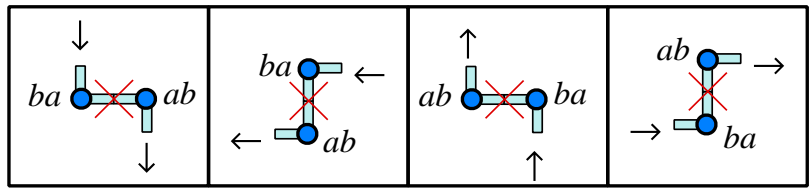

**Figure 3.3** : The four possible orientations of the ...$(ba)(ab)$... tiles at the middle of the tiling path associated with $\mathcal{W}_n$ with $n \equiv 1$ (mod 6).

Case 2 : $n \equiv 4$ (mod 6).

The number of letters in $\mathcal{F}_n$ will be of the form $4k$ for some non-negative integer value $k$. With the last two (rightmost) letters removed what remains, $\mathcal{W}_n$, is a palindrome (by Lemma 2.1) and of length $4k - 2$. As the subword $bb$ is not legal, at the middle of this palindrome must be the letter configuration ...$(aa)$... which corresponds geometrically under the double letter drawing rule to one of the tiling configurations shown in figure 3.4;

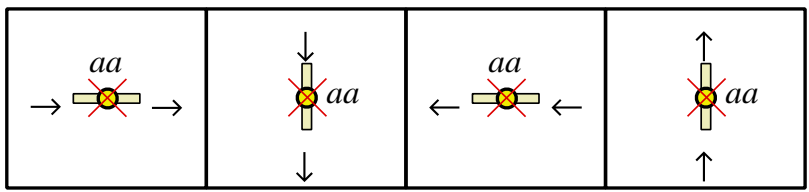

**Figure 3.4** : The four possible orientations of the ...$(aa)$... tile at the middle of the tiling path associated with $\mathcal{W}_n$ with $n \equiv 4$ (mod 6).

In figures 3.3 and 3.4 the arrows indicate the direction in which the tile is traversed when it is considered as a piece of path. Notice that all of the possible configurations in these figures have half turn rotational symmetry about the centre marked with a red cross and also that all have an exit direction the same as their entrance direction. Cases 1 and 2 now combine to give a basis for an inductive proof.

For the assumption step, we assume that there is some existing piece of tiled path that has half turn rotational symmetry about its middle and with an entrance and exit in the same direction. This is represented diagrammatically as shown in figure 3.5 and, without loss of generality, it will be assumed that the orientation and the direction through the tiling path is as shown.

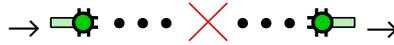

**Figure 3.5** : The assumption step in diagrammatic form representing a piece of path. The green half tile represents any of three possible tiles associated with (*ab*), (*aa*) or (*ba*). It is assumed that this piece of path has half turn rotational symmetry about the red cross at the middle of the piece. The orientation of this path and the entrance and exit directions are as shown, without loss of generality.

The piece of path in figure 3.5 will extend to the right in one of three possible ways, either through the addition of the tile associated with *ab,* or *aa* or *ba*. From the palindromic nature of the Fibonacci words, these three possible extensions to the right must match with extensions of *ba*, *aa*, or *ab* respectively to the left. Diagrammatically, the three possibilities are shown in figure 3.6.

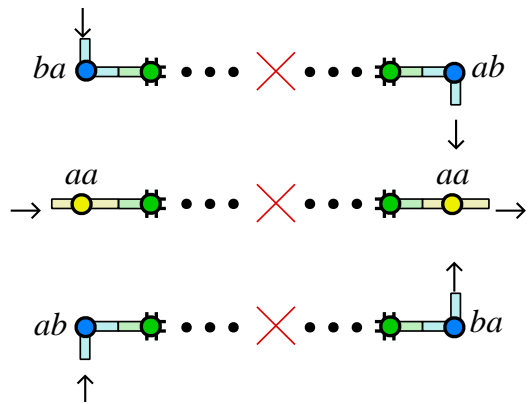

**Figure 3.6** : The inductive step, showing that all possible extensions of the path by one tile to the right cause a tile to be added to the left that retains the half turn rotational symmetry about the red cross and gives the entrance and exit directions that continue to be the same.

In all cases, after the inductive step, the tiling path continues to have half turn rotational symmetry about the red cross and the entrance and exit directions continue to be the same. The principle of induction can now be invoked to complete the proof.    □

### 3.4  Is The Path A Fractal ?

By way of further illustrating Lemma 3.1, figure 3.7 on the next page shows simplified paths (where the individual tiles are implied rather than shown explicitly) for the two examples considered previously, $\mathcal{F}_4$ and $\mathcal{F}_7$, alongside the next two Fibonacci words for which the *n* in $\mathcal{F}_n$ satisfies $n \equiv 1 \pmod 3$, $\mathcal{F}_{10}$ and $\mathcal{F}_{13}$. These images strengthen the suggestion that the drawing rule is giving rise to a path that is fractal in nature. The aim of this chapter is to prove that this is so. Several of the ideas used in this section are inspired by the 2009 paper "The Fibonacci Word Fractal" by Alexis Monnerot-Dumaine [Mon09] where they are applied to a different but related fractal.

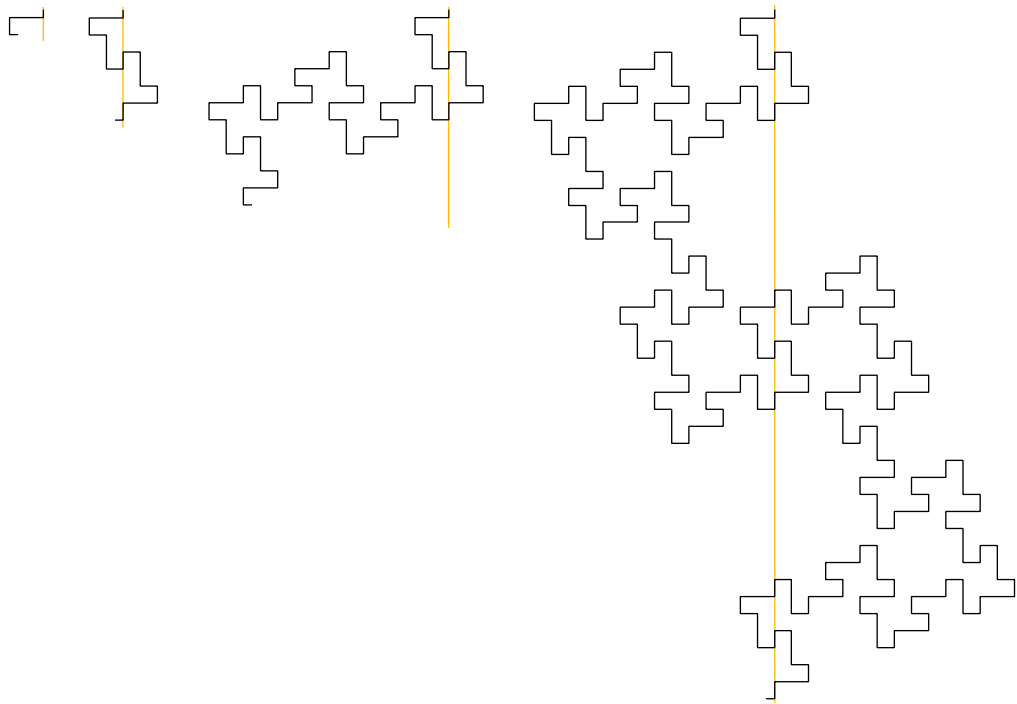

**Figure 3.7** : From left to right, the tiling paths associated with $\mathcal{F}_4$, $\mathcal{F}_7$, $\mathcal{F}_{10}$ and $\mathcal{F}_{13}$ under the double letter drawing rule. The paths start at the top of the diagram and finish at the bottom. With the very last tile removed these paths represent $\mathcal{W}_4$, $\mathcal{W}_7$, $\mathcal{W}_{10}$ and $\mathcal{W}_{13}$.

Theorem 3.1 is now stated; a crucial result underpinning our subsequent work.

---

**Theorem 3.1 : The Fractal Structure of $\mathcal{W}_n$** (**Monnerot-Dumaine** [Mon09])

Fibonacci words, $\mathcal{F}_n$, with the last two letters removed are denoted $\mathcal{W}_n$ and, from Lemma 2.1, are palindromic, $n \geq 1$. Let the two letters removed be $uv$ where $uv$ is one of the two possible endings; either $ab$ (with $u = a$ and $v = b$ when $n$ is odd) or $ba$ (with $u = b$ and $v = a$ when $n$ is even).

Then, for $n \geq 6$,

$$\mathcal{W}_n = \mathcal{W}_{n-3}\,(vu)\,\mathcal{W}_{n-3}\,(vu)\,\mathcal{W}_{n-6}\,(uv)\,\mathcal{W}_{n-3}\,(uv)\,\mathcal{W}_{n-3}$$

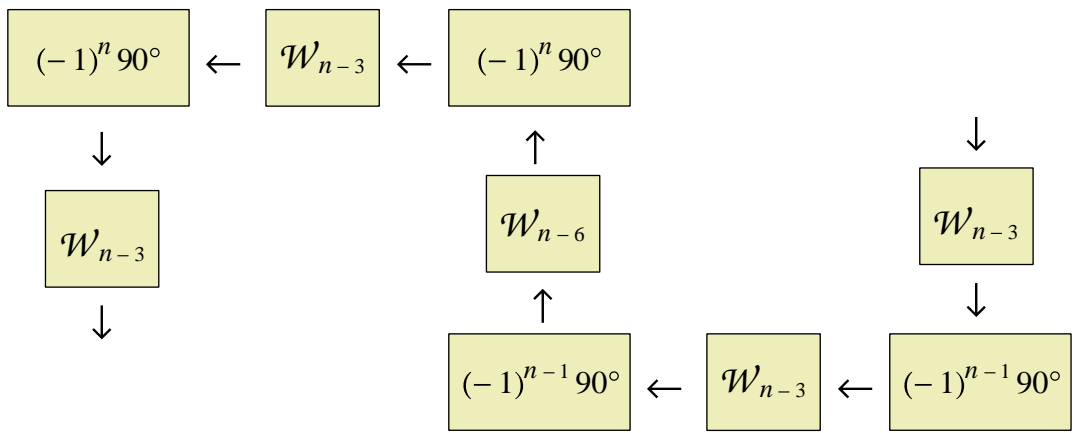

In particular this is true when $n \equiv 1 \pmod{3}$ corresponding to Fibonacci words that are valid under the double letter drawing rule.

---

Intuitively, a fractal is generally regarded as an endless, repeating, (or "almost" repeating) pattern within which can be found echos of its parts at different scales. Indeed, the whole maybe found within itself at a reduced scale. Theorem 3.1 shows that the Fibonacci word under the double letter drawing rule is giving rise to an entity of this nature. The proof of Theorem 3.1 will require that a few lemmas be established first and these will be presented shortly. First, however, figure 3.8 provides an illustrative example for Theorem 3.1 where the tiled path associated with $\mathcal{W}_{10}$ is shown broken down into its component parts as given by $\mathcal{W}_{10} = \mathcal{W}_7\,(ab)\,\mathcal{W}_7\,(ab)\,\mathcal{W}_4\,(ba)\,\mathcal{W}_7\,(ba)\,\mathcal{W}_7$.

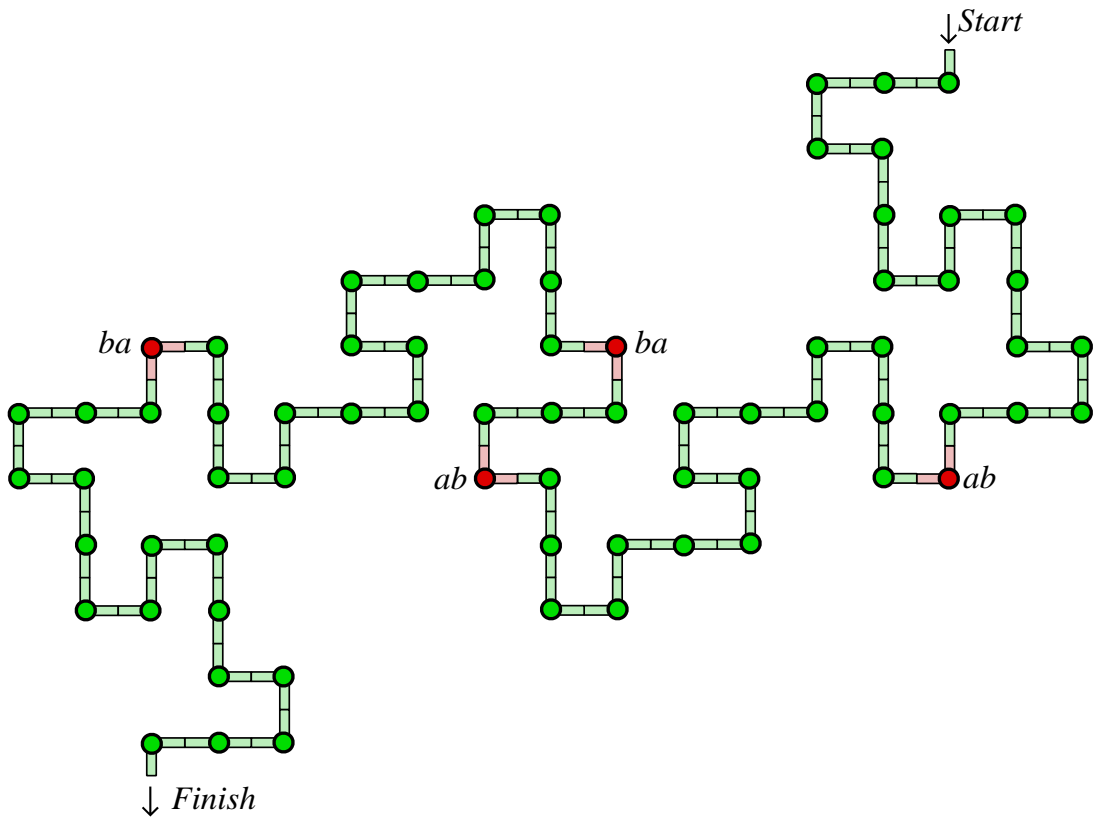

**Figure 3.8** : Illustrating Theorem 3.1 for $\mathcal{W}_{10} = \mathcal{W}_7\,(ab)\,\mathcal{W}_7\,(ab)\,\mathcal{W}_4\,(ba)\,\mathcal{W}_7\,(ba)\,\mathcal{W}_7$.

---

**Lemma 3.2 : End Letter Pair Reversal Rule   (Monnerot-Dumaine [Mon09])**

Let $\mathcal{F}_n = \mathcal{W}_n(uv)$ where $u$ and $v$ represent the penultimate and last letters respectively of $\mathcal{F}_n$. Let $\mathcal{T}_n = \mathcal{W}_n(vu)$ where the two end letters have swapped places. Then,

$$\mathcal{T}_{n+1} = \mathcal{F}_n\,\mathcal{T}_{n-1}$$

---

*Proof*

The Fibonacci words, by definition, have the relationship,

$\quad\mathcal{F}_{n+1} = \mathcal{F}_n\,\mathcal{F}_{n-1}\qquad$ for $n \geq 1$, $\mathcal{F}_0 = a$, $\mathcal{F}_1 = ab$

$\mathcal{W}_{n+1}(vu) = \mathcal{F}_n\,\mathcal{W}_{n-1}(vu)\qquad \mathcal{F}_{n+1},\,\mathcal{F}_{n-1}$ have the same two letter ending

$\mathcal{W}_{n+1}(uv) = \mathcal{F}_n\,\mathcal{W}_{n-1}(uv)\qquad$ reverse the two letter ending on both sides

$\quad\mathcal{T}_{n+1} = \mathcal{F}_n\,\mathcal{T}_{n-1}\qquad\qquad\qquad\qquad\qquad\qquad\qquad\square$

**Lemma 3.3 : "Almost" Commutativity**  (**Monnerot-Dumaine** [Mon09])

Let $\mathcal{F}_n = \mathcal{W}_n(uv)$ where $u$ and $v$ represent the penultimate and last letters respectively of $\mathcal{F}_n$. Let $\mathcal{T}_n = \mathcal{W}_n(vu)$ where the two end letters have swapped places. Then,

$$\mathcal{F}_n = \mathcal{F}_{n-2}\,\mathcal{T}_{n-1} \quad \text{and} \quad \mathcal{T}_n = \mathcal{F}_{n-2}\,\mathcal{F}_{n-1}$$

*Proof*

From $\mathcal{F}_n = \mathcal{W}_n(uv)$ note that, for example, $\mathcal{F}_{n+1} = \mathcal{W}_{n+1}(vu)$ because the two letter endings of the Fibonacci words alternate between *ab* and *ba*. From $\mathcal{T}_n = \mathcal{W}_n(vu)$ note that, for example, $\mathcal{T}_{n+1} = \mathcal{W}_{n+1}(uv)$ for the same reason. From the definition of the Fibonacci substitution acting repeatedly on $\mathcal{F}_0 = a$, we have that $\mathcal{F}_n$ for $1 \leqslant n \leqslant 4$ are as follows along with the corresponding $\mathcal{T}_n$,

- $\mathcal{F}_1 = ab, \ \mathcal{T}_1 = ba$
- $\mathcal{F}_2 = aba, \ \mathcal{T}_2 = aab$
- $\mathcal{F}_3 = abaab, \ \mathcal{T}_3 = ababa$
- $\mathcal{F}_4 = abaababa, \ \mathcal{T}_4 = abaabaab$

By way of establishing a basis for a proof by induction consider the case $n = 3$. Then, $\mathcal{F}_{3-2}\,\mathcal{T}_{3-1} = \mathcal{F}_1\,\mathcal{T}_2 = (ab)(aab) = \mathcal{F}_3$ which matches the lemma's first claim, and also $\mathcal{F}_{3-2}\,\mathcal{F}_{3-1} = \mathcal{F}_1\,\mathcal{F}_2 = (ab)(aba) = \mathcal{T}_3$ which matches the second.

It is now assumed the result is true for $n = k$ for some $k \in \mathbb{Z}, k \geqslant 4$.

In other words, $\mathcal{F}_k = \mathcal{F}_{k-2}\,\mathcal{T}_{k-1}$ and $\mathcal{T}_k = \mathcal{F}_{k-2}\,\mathcal{F}_{k-1}$ by assumption.

For the inductive step, consider $n = k + 1$, in which case,

$$\begin{aligned}
\mathcal{F}_{k+1} &= \mathcal{F}_k\,\mathcal{F}_{k-1} && \text{(By definition)} \\
&= \mathcal{F}_{k-1}\,\mathcal{F}_{k-2}\,\mathcal{F}_{k-1} \\
&= \mathcal{F}_{k-1}\,\mathcal{T}_k && \text{(From the assumptive step)}
\end{aligned}$$

which is the first result with $n$ replaced with $k + 1$.

Also,

$$\begin{aligned}
\mathcal{F}_{k-1}\,\mathcal{F}_k &= \mathcal{F}_{k-1}\,\mathcal{F}_{k-2}\,\mathcal{T}_{k-1} && \text{(From the assumptive step)} \\
&= \mathcal{F}_k\,\mathcal{T}_{k-1} && \text{(By definition)} \\
&= \mathcal{T}_{k+1} && \text{(By Lemma 3.1)}
\end{aligned}$$

which is the second result with $n$ replaced with $k + 1$,

Induction now gives the required pair of results. □

Lemma 3.3 is telling us that if we take the concatenation of two consecutive Fibonacci words $\mathcal{F}_{n-1}$ with $\mathcal{F}_{n-2}$ (which is equal to $\mathcal{F}_n$) and concatenate them the other way round as $\mathcal{F}_{n-2}$ with $\mathcal{F}_{n-1}$ then the resulting word is very nearly still $\mathcal{F}_n$. In fact, the lemma tells us that this resulting word is $\mathcal{F}_n$ but with the last two letters exchanged. One way of remembering this is to think of two consecutive Fibonacci words being "almost commutative".

The preparations are now done, and we are all set to prove Theorem 3.1.

**Theorem 3.1 : The Fractal Structure of $\mathcal{W}_n$  (Monnerot-Dumaine** [Mon09]**)**
Fibonacci words, $\mathcal{F}_n$, with the last two letters removed are denoted $\mathcal{W}_n$ and, from Lemma 2.1, are palindromic, $n \geqslant 1$. Let the two letters removed be $uv$ where $uv$ is one of the two possible endings; either $ab$ (with $u = a$ and $v = b$ when $n$ is odd) or $ba$ (with $u = b$ and $v = a$ when $n$ is even).
Then, for $n \geqslant 6$,

$$\mathcal{W}_n = \mathcal{W}_{n-3}\,(vu)\,\mathcal{W}_{n-3}\,(vu)\,\mathcal{W}_{n-6}\,(uv)\,\mathcal{W}_{n-3}\,(uv)\,\mathcal{W}_{n-3}$$

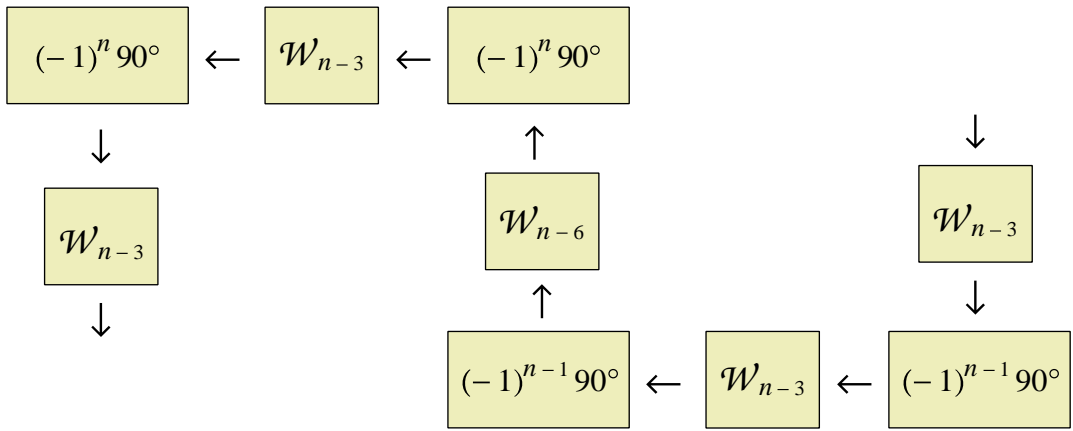

In particular this is true when $n \equiv 1 \pmod 3$ corresponding to Fibonacci words that are valid under the double letter drawing rule.

*Proof*

$$\begin{aligned}
\mathcal{F}_n &= \mathcal{F}_{n-1}\,\mathcal{F}_{n-2} & \text{(By definition)} \\
&= \mathcal{F}_{n-2}\,\mathcal{F}_{n-3}\,\mathcal{F}_{n-3}\,\mathcal{F}_{n-4} \\
&= \mathcal{F}_{n-3}\,\mathcal{F}_{n-4}\,\mathcal{F}_{n-4}\,\mathcal{F}_{n-5}\,\mathcal{F}_{n-3}\,\mathcal{F}_{n-4} \\
&= \mathcal{F}_{n-3}\,\mathcal{F}_{n-4}\,\mathcal{F}_{n-5}\,\mathcal{F}_{n-6}\,\mathcal{F}_{n-5}\,\mathcal{F}_{n-4}\,\mathcal{F}_{n-5}\,\mathcal{F}_{n-4} \\
&= \mathcal{F}_{n-3}\,\mathcal{F}_{n-3}\,\mathcal{F}_{n-6}\,\mathcal{T}_{n-3}\,\mathcal{T}_{n-3} & \text{(By Lemma 3.3)} \\
&= \mathcal{W}_{n-3}\,(vu)\,\mathcal{W}_{n-3}\,(vu)\,\mathcal{F}_{n-6}\,\mathcal{W}_{n-3}\,(uv)\,\mathcal{W}_{n-3}\,(uv)
\end{aligned}$$

Thus, $\mathcal{W}_n = \mathcal{W}_{n-3}\,(vu)\,\mathcal{W}_{n-3}\,(vu)\,\mathcal{W}_{n-6}\,(uv)\,\mathcal{W}_{n-3}\,(uv)\,\mathcal{W}_{n-3}$  □

### 3.5  A Vector Relationship Between Tiling Paths

Theorem 3.1 proves that the double letter drawing rule has given rise to a fractal. It is also the key to exploring other aspects of the meandering path. One striking feature of figure 3.7 is that the exit of the tiling path associated with $\mathcal{W}_7$ and $\mathcal{W}_{13}$ is directly vertically below the entrance. Figure 3.9 and figure 3.10 extend the sequence of fractals presented in figure 3.7 and show that $\mathcal{W}_{19}$ is the next word with a tiling path that has this property.

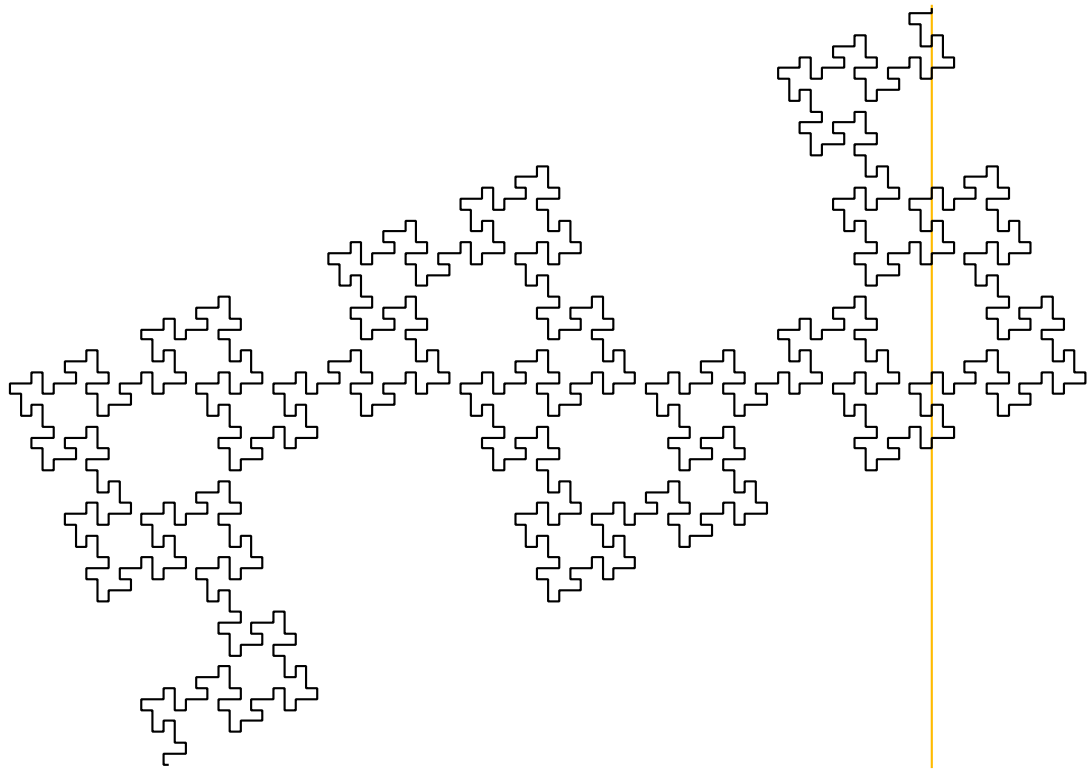

**Figure 3.9** : The tiling path associated with $\mathcal{F}_{16}$ under the double letter drawing rule. The path starts top and finishes bottom. $\mathcal{W}_{16}$ has the same path but with the very last tile removed.

Let $\overrightarrow{\mathcal{W}_n}$ be the vector with its tail at the start of the tiling path for $\mathcal{W}_n$ and its head at the finish. Table 3.3 gives $\overrightarrow{\mathcal{W}_n}$ for some values of $n$.

| $\mathcal{W}_n$ | N° tiles | $\overrightarrow{\mathcal{W}_n}$ | $\mathcal{W}_n$ | N° tiles | $\overrightarrow{\mathcal{W}_n}$ |
|---|---|---|---|---|---|
| $\mathcal{W}_4$ | 3 | $\begin{pmatrix} -2 \\ -1 \end{pmatrix}$ | $\mathcal{W}_{13}$ | 304 | $\begin{pmatrix} 0 \\ -40 \end{pmatrix}$ |
| $\mathcal{W}_7$ | 16 | $\begin{pmatrix} 0 \\ -6 \end{pmatrix}$ | $\mathcal{W}_{16}$ | 1291 | $\begin{pmatrix} -70 \\ -69 \end{pmatrix}$ |
| $\mathcal{W}_{10}$ | 71 | $\begin{pmatrix} -12 \\ -11 \end{pmatrix}$ | $\mathcal{W}_{19}$ | 5472 | $\begin{pmatrix} 0 \\ -238 \end{pmatrix}$ |

**Table 3.3** : Observed values of $\overrightarrow{\mathcal{W}_n}$ from careful counting on figures 3.7, 3.9 and 3.10.

Two previous vectors can be used to determine a subsequent vector in table 3.3 by making use of Theorem 3.1 and taking care over how the 90° turns flip the vector. For example, for $n \equiv 1 \pmod 6$, $\mathcal{W}_{13}$ is obtained from $\mathcal{W}_{10}$ and $\mathcal{W}_7$ by,

$$\overrightarrow{\mathcal{W}_{13}} = \begin{pmatrix} -12 \\ -11 \end{pmatrix} + \begin{pmatrix} \frac{1}{2} \\ -\frac{1}{2} \end{pmatrix} + \begin{pmatrix} 11 \\ -12 \end{pmatrix} + \begin{pmatrix} \frac{1}{2} \\ \frac{1}{2} \end{pmatrix} + \begin{pmatrix} 0 \\ 6 \end{pmatrix} + \begin{pmatrix} \frac{1}{2} \\ \frac{1}{2} \end{pmatrix} + \begin{pmatrix} 11 \\ -12 \end{pmatrix} + \begin{pmatrix} \frac{1}{2} \\ -\frac{1}{2} \end{pmatrix} + \begin{pmatrix} -12 \\ -11 \end{pmatrix}$$

$$= \begin{pmatrix} 0 \\ -40 \end{pmatrix}$$

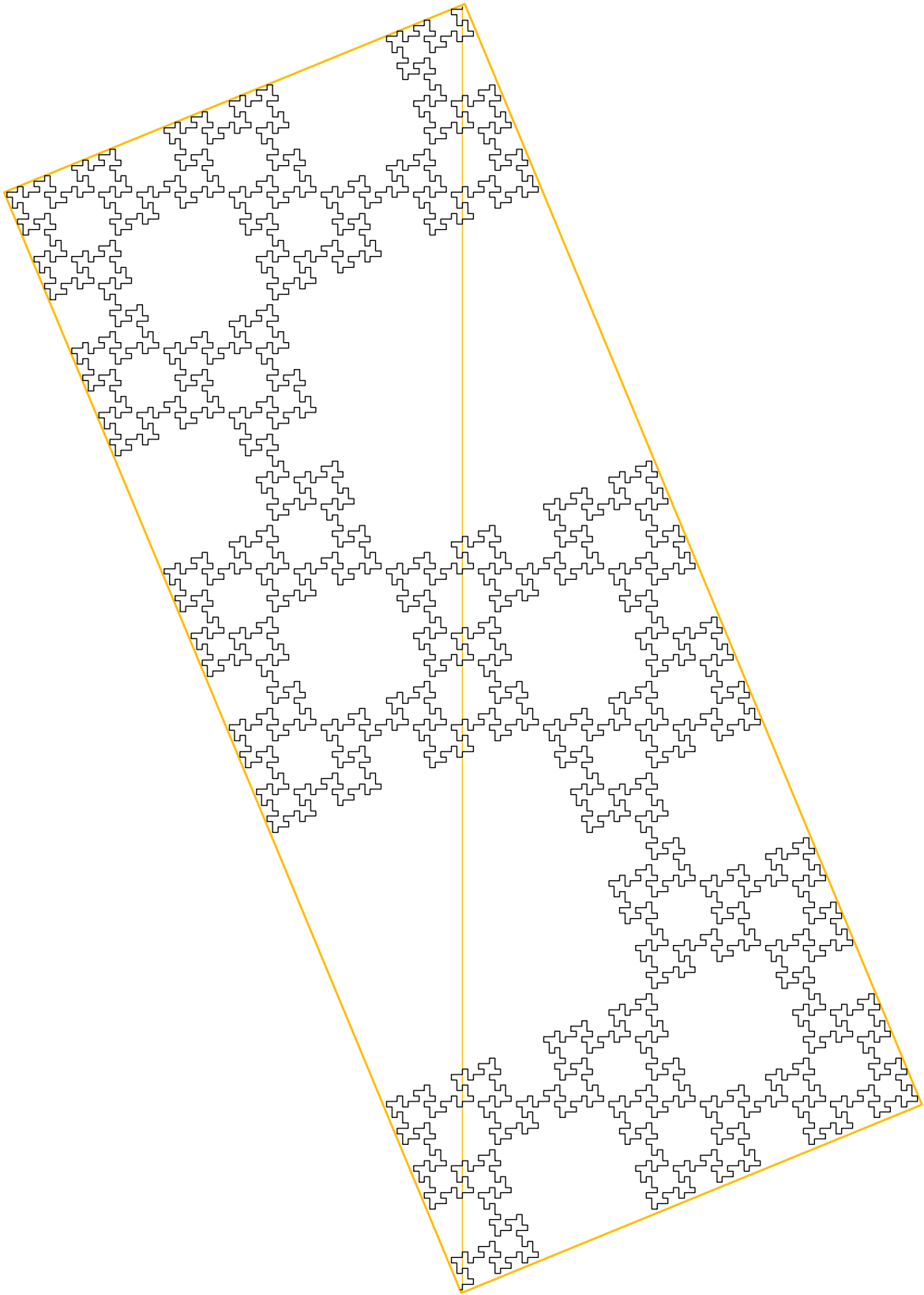

**Figure 3.10** : The tiling path associated with $\mathcal{F}_{19}$ under the double letter drawing rule. The path starts at the top of the diagram and finishes at the bottom. A bounding rectangle has been fitted around the fractal with a vertical diagonal showing that the entrance is directly above the exit of the tiling path. With the very last tile removed this is also the tiling path for $\mathcal{W}_{19}$.

As a further example, for $n \equiv 4 \pmod 6$, $\mathcal{W}_{10}$ is obtained from $\mathcal{W}_7$ and $\mathcal{W}_4$ by,

$$\overrightarrow{\mathcal{W}_{10}} = \begin{pmatrix} 0 \\ -6 \end{pmatrix} + \begin{pmatrix} -\frac{1}{2} \\ -\frac{1}{2} \end{pmatrix} + \begin{pmatrix} -6 \\ 0 \end{pmatrix} + \begin{pmatrix} -\frac{1}{2} \\ \frac{1}{2} \end{pmatrix} + \begin{pmatrix} 2 \\ 1 \end{pmatrix} + \begin{pmatrix} -\frac{1}{2} \\ \frac{1}{2} \end{pmatrix} + \begin{pmatrix} -6 \\ 0 \end{pmatrix} + \begin{pmatrix} -\frac{1}{2} \\ -\frac{1}{2} \end{pmatrix} + \begin{pmatrix} 0 \\ -6 \end{pmatrix}$$

$$= \begin{pmatrix} -12 \\ -11 \end{pmatrix}$$

Lemma 3.4 provides general formulae for such calculations.

---

**Lemma 3.4 : Iterative Formulae for $\overrightarrow{\mathcal{W}_n}$**

Let $\overrightarrow{\mathcal{W}_{n-3}} = \begin{pmatrix} p \\ q \end{pmatrix}$ and $\overrightarrow{\mathcal{W}_{n-6}} = \begin{pmatrix} r \\ s \end{pmatrix}$ where $p$, $q$, $r$ and $s$ are integer constants and $n \equiv 1 \pmod 3$ such than the vectors represent displacements between the entrance and exit of the associated tiling paths under the double letter drawing rule on the Fibonacci word. Then,

$$\overrightarrow{\mathcal{W}_n} = \begin{cases} \begin{pmatrix} 2p - 2q - r + 2 \\ 2p + 2q - s \end{pmatrix} & \text{for } n \equiv 1 \pmod 6 \\ \begin{pmatrix} 2p + 2q - r - 2 \\ -2p + 2q - s \end{pmatrix} & \text{for } n \equiv 4 \pmod 6 \end{cases}$$

---

*Proof*
From Theorem 3.1, taking care over the orientation of the vectors implied by the 90° turns, and with $\overrightarrow{\mathcal{W}_{n-3}}$ and $\overrightarrow{\mathcal{W}_{n-6}}$ as defined above,

Case 1 : For $n \equiv 1 \pmod 6$,

$$\overrightarrow{\mathcal{W}_n} = \begin{pmatrix} p \\ q \end{pmatrix} + \begin{pmatrix} \frac{1}{2} \\ -\frac{1}{2} \end{pmatrix} + \begin{pmatrix} -q \\ p \end{pmatrix} + \begin{pmatrix} \frac{1}{2} \\ \frac{1}{2} \end{pmatrix} + \begin{pmatrix} -r \\ -s \end{pmatrix} + \begin{pmatrix} \frac{1}{2} \\ \frac{1}{2} \end{pmatrix} + \begin{pmatrix} -q \\ p \end{pmatrix} + \begin{pmatrix} \frac{1}{2} \\ -\frac{1}{2} \end{pmatrix} + \begin{pmatrix} p \\ q \end{pmatrix}$$

$$= \begin{pmatrix} 2p - 2q - r + 2 \\ 2p + 2q - s \end{pmatrix} \text{ as claimed.}$$

Case 2 : For $n \equiv 4 \pmod 6$,

$$\overrightarrow{\mathcal{W}_n} = \begin{pmatrix} p \\ q \end{pmatrix} + \begin{pmatrix} -\frac{1}{2} \\ -\frac{1}{2} \end{pmatrix} + \begin{pmatrix} q \\ -p \end{pmatrix} + \begin{pmatrix} -\frac{1}{2} \\ \frac{1}{2} \end{pmatrix} + \begin{pmatrix} -r \\ -s \end{pmatrix} + \begin{pmatrix} -\frac{1}{2} \\ \frac{1}{2} \end{pmatrix} + \begin{pmatrix} q \\ -p \end{pmatrix} + \begin{pmatrix} -\frac{1}{2} \\ -\frac{1}{2} \end{pmatrix} + \begin{pmatrix} p \\ q \end{pmatrix}$$

$$= \begin{pmatrix} 2p + 2q - r - 2 \\ -2p + 2q - s \end{pmatrix} \text{ as claimed.}$$

□

Previously, in table 3.3, there appeared to be an alternation between vectors of the form $\overrightarrow{W_n} = \begin{pmatrix} 0 \\ w \end{pmatrix}$ when $n \equiv 1$ (mod 6) and $\overrightarrow{W_n} = \begin{pmatrix} x \\ x+1 \end{pmatrix}$ when $n \equiv 4$ (mod 6). Our next result, Lemma 3.5, shows that this alternation remains a feature as $n$ tends towards infinity.

---

**Lemma 3.5 : Alternating Nature of $\overrightarrow{W_n}$**

For all $n \equiv 1$ (mod 3) the vector $\overrightarrow{W_n}$ is of the form $\begin{pmatrix} 0 \\ w \end{pmatrix}$ when $n \equiv 1$ (mod 6) and of the form $\begin{pmatrix} x \\ x+1 \end{pmatrix}$ when $n \equiv 4$ (mod 6). In other words $\overrightarrow{W_n}$ alternates between these two forms repeatedly as $n$ tends to infinity.

---

*Proof*

Case 1 : For $n \equiv 1$ (mod 6),
$\overrightarrow{W_{n-3}}$ is of the form $\begin{pmatrix} x \\ x+1 \end{pmatrix}$ and $\overrightarrow{W_{n-6}}$ is of the form $\begin{pmatrix} 0 \\ w \end{pmatrix}$.

From Lemma 3.4, $\overrightarrow{W_n} = \begin{pmatrix} 2p - 2q - r + 2 \\ 2p + 2q - s \end{pmatrix}$ with $p = x$, $q = x+1$, $r = 0$, $s = w$

$$\overrightarrow{W_n} = \begin{pmatrix} 0 \\ 4x + 2 - w \end{pmatrix}$$

This is of the required form $\begin{pmatrix} 0 \\ w' \end{pmatrix}$ where $w' = 4x + 2 - w$.

Case 2 : For $n \equiv 4$ (mod 6),
$\overrightarrow{W_{n-3}}$ is of the form $\begin{pmatrix} 0 \\ w \end{pmatrix}$ and $\overrightarrow{W_{n-6}}$ is of the form $\begin{pmatrix} x \\ x+1 \end{pmatrix}$.

From Lemma 3.4, $\overrightarrow{W_n} = \begin{pmatrix} 2p + 2q - r - 2 \\ -2p + 2q - s \end{pmatrix}$ with $p = 0$, $q = w$, $r = x$, $s = x+1$

$$\overrightarrow{W_n} = \begin{pmatrix} (2w - x - 2) \\ (2w - x - 2) + 1 \end{pmatrix}$$

This is of the required form $\begin{pmatrix} x' \\ x'+1 \end{pmatrix}$ where $x' = 2w - x - 2$.

Given that $\overrightarrow{W_4}$ is of the form $\begin{pmatrix} x \\ x+1 \end{pmatrix}$ with $x = -2$ followed by $\overrightarrow{W_7}$ being of the form $\begin{pmatrix} 0 \\ w \end{pmatrix}$ with $w = -6$, we have an initial alternating basis to which an inductive argument is then applied.

Alternate use of the case 1 and case 2 results, as appropriate, proves that, as $n$ tends to infinity, the claimed alternation of $\overrightarrow{W_n}$ between the two forms of vector continues to hold. □

### 3.6 Determining Fractal Dimension

An important measure attached to a fractal is its dimension. There are several different measures of dimension in common usage and there is a skill in selecting which is most appropriate for a given fractal. When dealing with relatively uncomplicated fractals which have component parts (for example, line segments) that do not excessively overlap, the generic term "fractal dimension" is sometimes used, especially if the different measures of dimension give the same numerical result. For our fractal the dimension can be determined using a modified version of the formula previously given in section 3.2;

$$dim_{SYM} = \frac{\log m}{\log s} \quad \left(\text{where } s \text{ is } \frac{1}{r} \text{ for the } r \text{ of the earlier formula}\right).$$

In this, $m$ is the ratio between $\mathcal{F}_{n-3}$ and $\mathcal{F}_n$ (in the limit) at which the number of line segments is increasing and $s$ is the expansive length scale factor (in the limit) at which (for example) the corresponding diagonal of a rectangular box that is bounding the fractal is increasing. This rectangular bounding box becomes more intuitively obvious as the number of iterations increases and is shown explicitly in figure 3.10 as an orange outline along with the diagonal corresponding to $\overrightarrow{\mathcal{W}}_{13}$. For our fractal, the similarity dimension will have the same value as that obtained from more sophisticated measures. In Kenneth Falconer's, *Fractal Geometry,* [Fal03], for example, box dimension, $dim_B$, and Hausdorff dimension, $dim_H$, will take on the same value. Indeed, for the example presented earlier in figure 3.1, $dim_B = dim_H = dim_{SYM} = 1.5$.

---

**Proposition 3.1 : Value of $m$**

For the fractal associated with the double letter drawing rule, in the limit, the value of $m$ in $dim_{SYM} = \frac{\log m}{\log s}$ is $\phi^3$, where $\phi = \frac{1+\sqrt{5}}{2}$ and $m$ is the limiting ratio at which the number of line segments between $\mathcal{F}_{n-3}$ and $\mathcal{F}_n$ is increasing.

---

*Proof*

The number of letters in the Fibonacci word $\mathcal{F}_n$ (see table 1.1) is given by the Fibonacci number $f(n)$ which, by definition, is given by,

$$f(n) = f(n-1) + f(n-2) \text{ for } n \geq 2 \text{ with } f(0) = 1, f(1) = 2$$

In consequence, $\dfrac{f(n)}{f(n-1)} = \dfrac{f(n-1) + f(n-2)}{f(n-1)} = 1 + \dfrac{f(n-2)}{f(n-1)}$

Taking limits, $\lim\limits_{n \to \infty} \dfrac{f(n)}{f(n-1)} = \lim\limits_{n \to \infty} 1 + \lim\limits_{n \to \infty} \dfrac{f(n-2)}{f(n-1)}$.

Let $x = \lim\limits_{n \to \infty} \dfrac{f(n)}{f(n-1)}$ in which case the above becomes $x = 1 + \dfrac{1}{x}$.

Solving gives the expansive value of $x$ sought. It is, $\phi = \dfrac{1+\sqrt{5}}{2}$.

Hence, $m = \lim\limits_{n \to \infty} \dfrac{f(n)}{f(n-3)} = \phi^3$.  □

Note that, in Proposition 3.1 it does not matter that it is pairs of letters that give rise to the number of line segments because it is a ratio that is being found. Also note that, because we are taking limits as $n$ tends to infinity, we can be relaxed about the one tile (two letter) difference between $\mathcal{F}_n$ and $\mathcal{W}_n$.

---

**Proposition 3.2 : Value of $s$**

For the fractal associated with the double letter drawing rule, in the limit, the value of $s$ in $dim_{SYM} = \dfrac{\log m}{\log s}$ is $1 + \sqrt{2}$ where $s$ is the expansive length scale factor in the limit between $\overrightarrow{\mathcal{W}_{n-3}}$ and $\overrightarrow{\mathcal{W}_n}$.

---

*Proof*

Consider the two consecutive vectors $\overrightarrow{\mathcal{W}_{n-6}} = \begin{pmatrix} x \\ x+1 \end{pmatrix}$ and $\overrightarrow{\mathcal{W}_{n-3}} = \begin{pmatrix} 0 \\ w \end{pmatrix}$ in which case $\overrightarrow{\mathcal{W}_n} = \begin{pmatrix} (2w - x - 2) \\ (2w - x - 2) + 1 \end{pmatrix}$ by Lemma 3.4.

As we are dealing with an expansive iterative process, as the number of iterations increases, $w$ and $x$ become large with $\overrightarrow{\mathcal{W}_{n-6}} \to \begin{pmatrix} x \\ x \end{pmatrix}$ and $\overrightarrow{\mathcal{W}_n} \to \begin{pmatrix} 2w - x \\ 2w - x \end{pmatrix}$. The lengths of these three vectors are tending towards the following; $\left|\overrightarrow{\mathcal{W}_{n-6}}\right| \to \sqrt{2}\, x$, $\left|\overrightarrow{\mathcal{W}_{n-3}}\right| \to w$ and $\left|\overrightarrow{\mathcal{W}_n}\right| \to \sqrt{2}\,(2w - x)$. In the limit, the expansive length scale factor is given by,

$\dfrac{w}{\sqrt{2}\, x} = \dfrac{\sqrt{2}\,(2w - x)}{w}$ from which, $w^2 - 4wx + 2x^2 = 0$.

Solving this quadratic equation gives, $w = \sqrt{2}\, x\,(\sqrt{2} \pm 1)$.

For an expansive value, it follows that, $s = \dfrac{w}{\sqrt{2}\, x} = 1 + \sqrt{2}$   □

---

The chapter now concludes with the important result presented as Theorem 3.2;

---

**Theorem 3.2 : Value of Hausdorff Dimension**

The Hausdorff dimension of the fractal associated with the double letter drawing rule is $dim_H = \dfrac{\log \phi^3}{\log(1 + \sqrt{2})} = 1.6379$

---

*Proof*

Using the values of $m$ and $s$ from Propositions 3.1 and 3.2 respectively,

$$dim_H = dim_{SYM} = = \dfrac{\log m}{\log s} = \dfrac{\log \phi^3}{\log(1 + \sqrt{2})} = 1.6379 \quad \square$$

### 3.7 Author's Notes on Chapter 3

The fractal known in the literature as "The Fibonacci Word Fractal" [Wik21b] arises from what is termed "the odd-even drawing rule" directly applied to the Fibonacci word, assigning one instruction to one letter. The instructions to draw the fractal tell the drawing pen to, starting before the first letter of $\mathcal{F}$,

- take the next letter and draw a unit line segment forward
- if the letter is an *a* then,
  - turn 90° left if the *a* is in an even position within $\mathcal{F}$
  - turn 90° right if the *a* is in an odd position within $\mathcal{F}$
- repeat the process indefinitely

Figure 3.11 gives an example of the path produced.

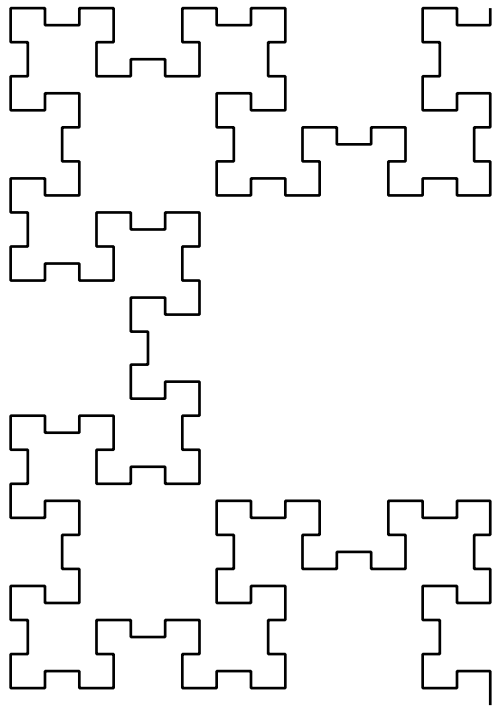

**Figure 3.11** : The meandering path produced by the odd-even drawing rule applied to the Fibonacci word, $\mathcal{F}_{12}$. The path starts at the top of the diagram and finishes at the bottom.

The odd-even drawing rule is not an obvious drawing rule to apply because of the need to know for each letter *a* if it is in an odd or even numbered location. The advantage of the rule, however, is that it draws a fractal that is more straight forward to analyse. To me, taking the letters in pairs seemed more natural. I was curious to know if the "alternative" fractal obtained was amenable to a relatively straight forward mathematical analysis; it was a relief to find a way of determining the Hausdorff dimension! In some respects the two fractals are different. For example, the odd-even fractal has a limiting ratio of width : height of $1 : \sqrt{2}$ whereas the double letter fractal's ratio is $1 : 1 + \sqrt{2}$ (The interested reader may like to confirm that result). And yet they have exactly the same value for their Hausdorff dimension.

# Chapter 4

# Future Directions

### 4.1 Generalised To and Fro

Inevitably, in a work of less than 12,000 words, we have but sampled some flavoursome fruit from the low hanging branches of a vast tree. The ideas from chapter 2 have obvious extensions into two dimensions even without opting for a completely different drawing rule like that of chapter 3. Table 4.1 gives a generalised to and fro drawing rule. With $\varphi = 180°$ it is as before but other values of $\varphi$ are waiting to be explored. For example, with $\varphi = 108°$ the attractive motif shown left in figure 4.1 is obtained form $\mathcal{F}_{12}$. Increasing the $n$ in $\mathcal{F}_n$ to 18 gives a massive amount of overwriting but reveals the path shown on the right of the same figure with a steady overall bend. Increasing $n$ further will cause it to eventually circle back on itself but is an annulus the eventual form? Chasing this idea a little further, figure 4.2 shows a similar circular path when $\varphi = 136°$ for $\mathcal{F}_{13}$ (left) and $\mathcal{F}_{18}$ (right). Computer explorations suggest there may be a special "fire hose" angle around $\varphi = 137.4°$ that projects the path forward overall. The name arises because finding this critical angle is akin to holding onto the tap end of a hosepipe. The hose is trying to snake up and down and wants to swing back on itself spraying water from its free end in all directions but forward. Finding such "fire hose" angles (if they exists) using mathematics may prove interesting.

| Symbol | Action |
| --- | --- |
| $a$ | forward $\phi$ |
| $b$ | forward 0.5, turn $\varphi°$, forward 0.5 |

**Table 4.1** : The generalised to and fro drawing rule making use of an angle $\varphi$ instead of 180°.

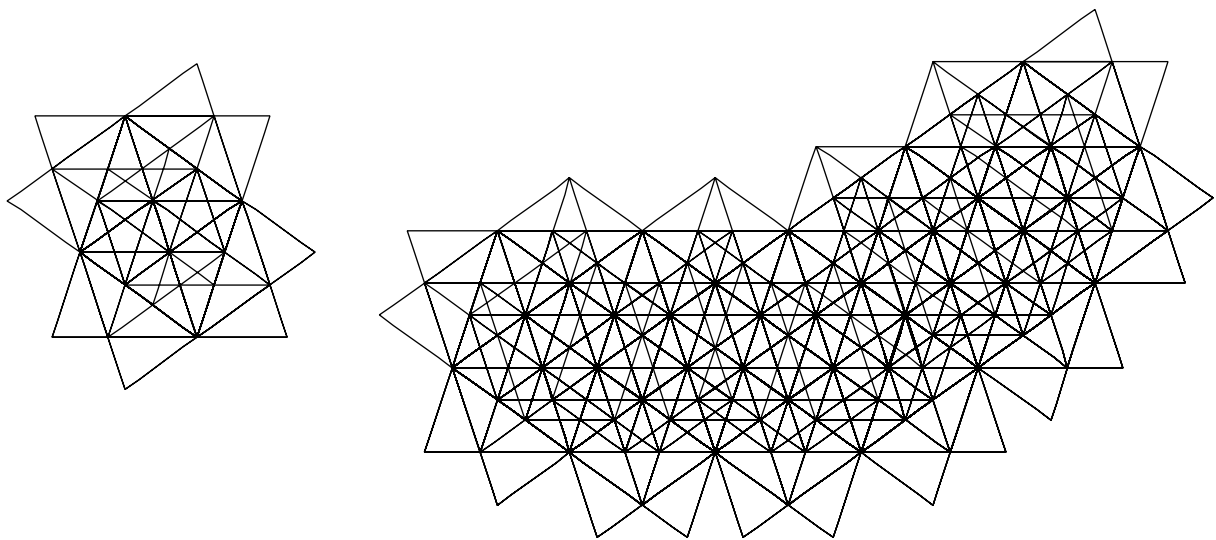

**Figure 4.1** : Left : The tiling path when the $\varphi = 108°$ drawing rule is applied to $\mathcal{F}_{12}$.
Right : This path has an overall bend as the $n$ in $\mathcal{F}_n$ is increased.

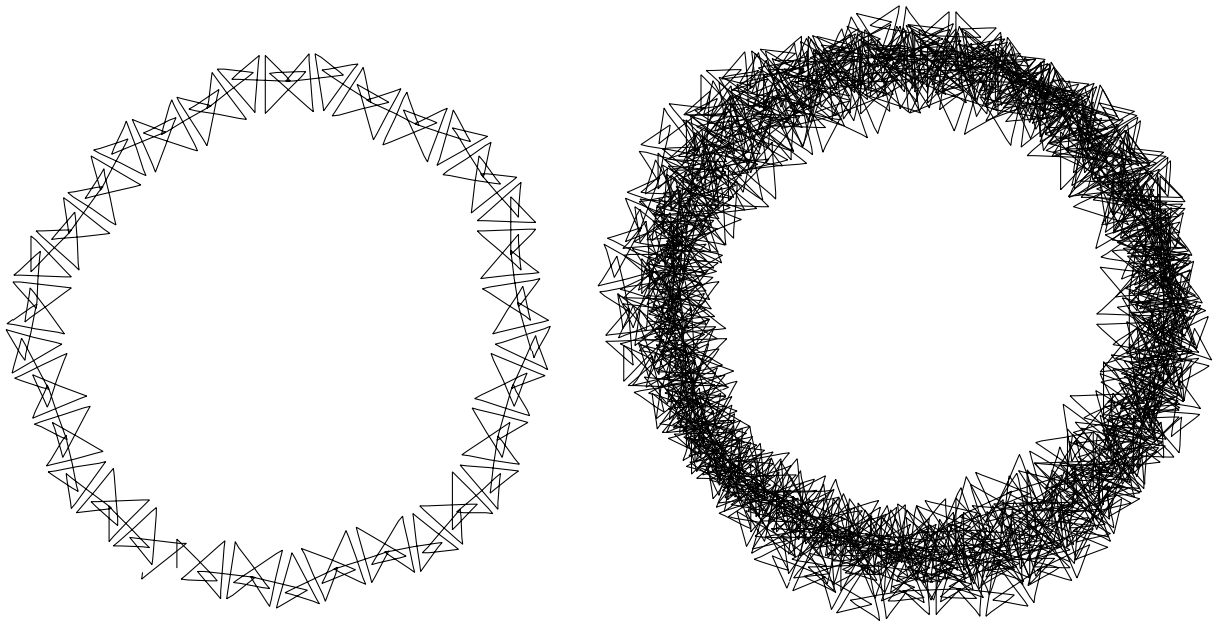

**Figure 4.2** : Left : The tiling path when the $\varphi = 136°$ drawing rule is applied to $\mathcal{F}_{13}$.
Right : The same rule applied to $\mathcal{F}_{18}$.

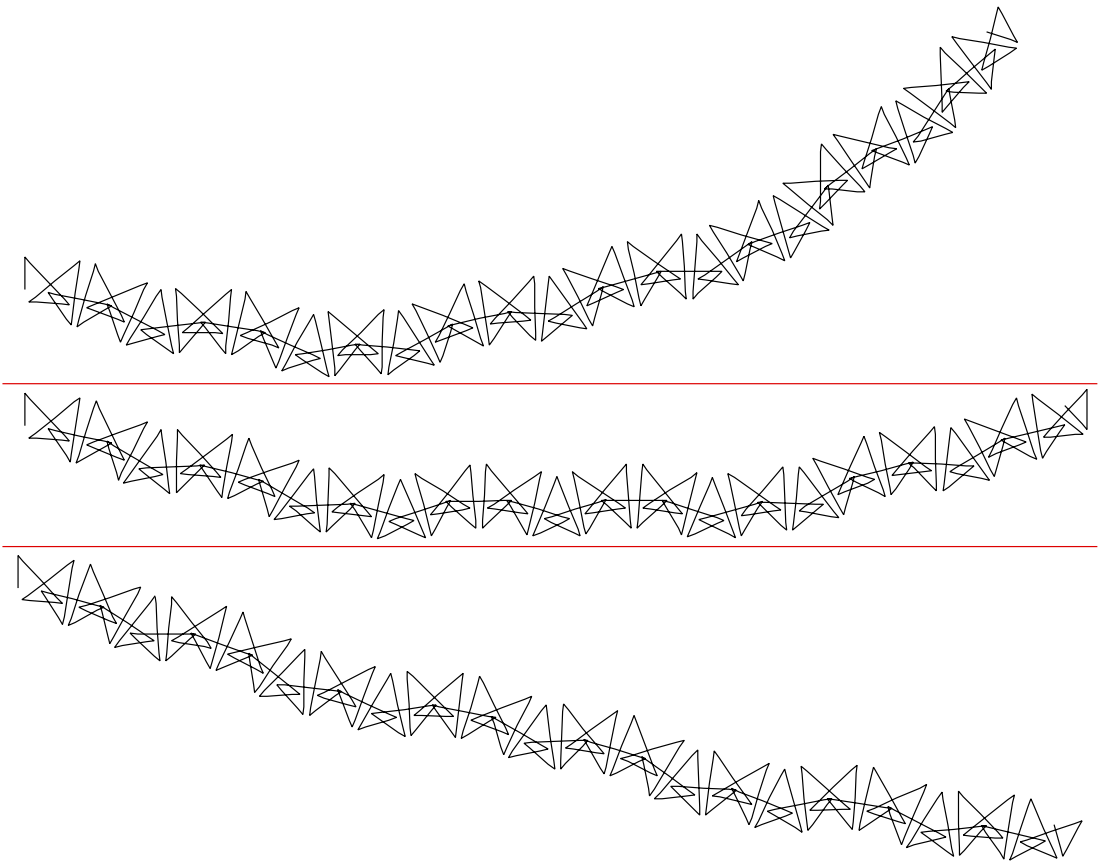

**Figure 4.3** : Using a computer to explore the possibility that there is a value of $\varphi$ in the generalised to and fro style drawing rule that, overall, projects the tiling path in a straight direction. Here the word $\mathcal{F}_{12}$ is suggesting that the angle sought is greater than $\varphi = 137.0°$ (upper diagram) and also greater than $\varphi = 137.2°$ (middle diagram) both of which gives a tiling path that seems to overall have an anticlockwise bend, whereas $\varphi = 137.4°$ gives an almost straight overall projection of the path (lower diagram). A striking feature of the path is that it seems to be composed from two motifs thus giving a visualisation of a factorisation of $\mathcal{F}$.

### 4.2 Self Avoiding Walks

In the Monnerot-Dumaine paper [Mon09], proving that the Fibonacci fractal does not intersect itself is identified as an open problem. There is a substantial body of work on self avoiding walks which may provide some pointers. By way of an initial look at this problem, Figure 4.4 shows a fractal with a straight line segment as initiator, $A_1$, of 1 unit length. The generator is the $A_2$ dogleg. The substitution, $\Omega$, underpinning the iterative generation of this fractal is;

$$\Omega : \begin{cases} F \rightarrow F L F R F R F L F \\ L \rightarrow L \\ R \rightarrow R \end{cases}$$

and the iterative fixed point of interest is $\Omega^n(F)$ as $n \rightarrow \infty$.
The drawing rule is as given in table 4.2.

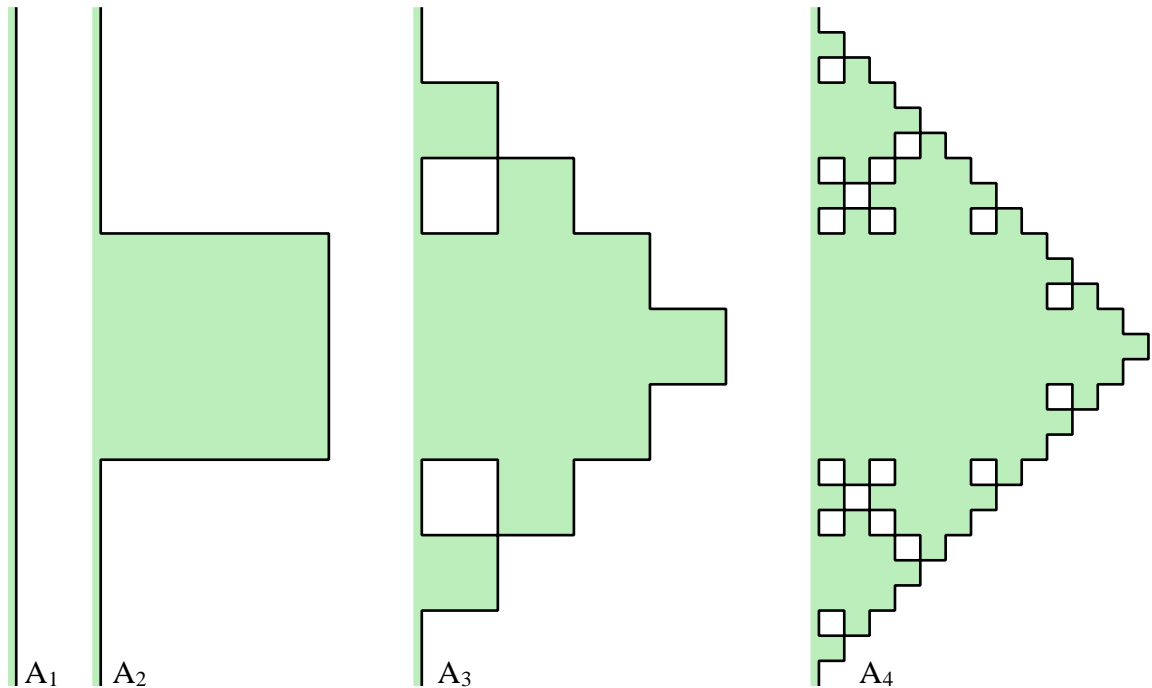

**Figure 4.4** : An example of a self intersecting fractal.

| Symbol | Action |
|---|---|
| $F$ | forward 1 |
| $L$ | turn left (clockwise) 90° |
| $R$ | turn right (anticlockwise) 90° |

**Table 4.2** : The drawing rule that acts on the substitution $\Omega^\infty(F)$ to give the fractal of figure 4.4.

The infinite word $\Omega^\infty(F)$ begins $F L F R F R \textcolor{red}{F L F L F L F} R F R F L F R$ ... but the embedded word highlighted in red is a loop that self intersects. This proves that this fractal is not a self avoiding walk. In general, determining if an infinite walk is self avoiding is difficult but, possibly, the almost periodic nature of the Fibonacci word means this is a problem that can be resolved in this special case. It would be interesting to work on this.

### 4.3 Letter Frequency Analysis

Tied in with the problem of determining if the double letter fractal is a self avoiding walk or not, would seem to be the issue of analysing the frequency with which the three letter pairs *ab*, *aa* and *ba* occur. Without going into the detail, it turns out that *ab* and *ba* have a relative frequency of about 38% with *aa* claiming the remaining 24%. However, it's not just these frequencies that are important; as we read through $\mathcal{F}$ from left to right, how much can the count of one letter pair get ahead or behind the others? There would also seem to be merit in looking into the possible factorisations of $\mathcal{F}$ and the frequencies of letters (or double letters) within and between factors. Theorem 2.1 gave a proof that one such factorisation exists; there are many others.

### 4.4 Periodic Approximations to the Aperiodic

Figure 4.5 shows a fractal with a straight line segment as initiator, $A_1$, of 1 unit length. The generator is the $A_2$ meander. This was inspired by the tiling path associated with $\mathcal{F}_7$ under the double letter drawing rule. This fractal can be thought of as a periodic approximation to our aperiodic fractal from $\mathcal{F}$ under the double letter drawing rule. This approximation has 17 segments and a scaling factor of one seventh giving a fractal dimension of 1.46. A sequence of such approximating a fractals based on $\mathcal{F}_{13}$, $\mathcal{F}_{19}$, $\mathcal{F}_{25}$, ... would seem to provide ever better periodic approximations. For example, $\mathcal{F}_{31}$ has 1762289 segments and a scaling factor of $8119^{-1}$ giving a fractal dimension of 1.60, which is getting closer to the 1.64 of the double letter fractal. Are there situations where it would be easier to work with these periodic approximations rather than the aperiodic real fractal ?

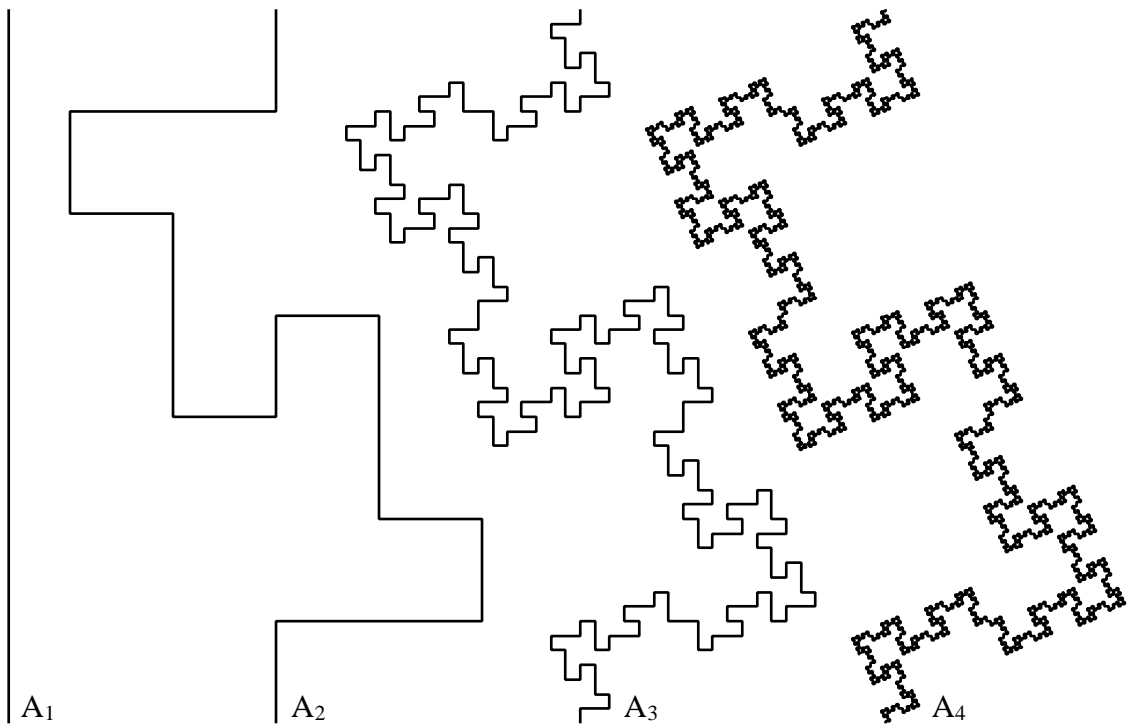

$A_1$   $A_2$   $A_3$   $A_4$

**Figure 4.5** : With a straight line initiator, using an $\mathcal{F}_7$ inspired generator to obtain a periodic fractal to approximate the aperiodic fractal obtained for the Fibonacci word, $\mathcal{F}$, under the double letter drawing rule. This approximation has fractal dimension 1.46 rather than the 1.64 of $\mathcal{F}$ but increasingly better approximations are obtained by using as inspiration $\mathcal{F}_{13}$, $\mathcal{F}_{19}$, $\mathcal{F}_{25}$, $\mathcal{F}_{31}$, ...

### 4.5  Author's Notes on Chapter 4

A major part of the attraction of working with the Fibonacci word is the abundance of patterns giving mathematical footholds and suggesting possible avenues to explore. In the established literature the Fibonacci substitution and the resulting word has been generalised in several different ways and it would certainly be interesting to extend this work to embrace them. Another possible future project would be to add an exercise to the end of each chapter by way of expanding and exploring further the content therein.

In researching this topic I became very aware that there is a large body of mathematics, far more sophisticated than that covered here, that analyses the Fibonacci word in ever more subtle ways, and from an ever more general viewpoint. However, the desire to take this work in that direction has been tempered for now by the fact that the dissertation deadline is upon me; it is time to stop!

# Appendix A1 :
## Remembering Professor Uwe Grimm

It was Uwe who introduced me to that most fascinating of mathematical toys, the infinite Fibonacci word. It's a deceptively simple substitution, $\theta$, on an alphabet of only two letters, $\mathcal{A}\{a,b\}$, defined by $a \to ab$ and $b \to a$. It gives us the finite Fibonacci words, $\mathcal{F}_n = \theta^n(a)$. The first few are; $\mathcal{F}_0 = a$, $\mathcal{F}_1 = ab$, $\mathcal{F}_2 = aba$, $\mathcal{F}_3 = abaab$ and so on. Throw away the last couple of letters on any given word and what's left is a palindrome. As example $\mathcal{F}_4 = abaababa$ which, without its rightmost two letters, is $abaaba$. This palindromic nature along with the remarkable concatenation property that $\mathcal{F}_n = \mathcal{F}_{n-1}\mathcal{F}_{n-2}$ for $n \geq 3$, guarantees that the Fibonacci words abound with symmetries. As $n \to \infty$ the infinite Fibonacci word emerges as a fixed point of the iteration.

Uwe took pleasure in finding geometric visualisations to complement his algebraic researches. These were often stunningly beautiful creations that non-mathematicians could marvel over. When he died, I had just begun studying the Open University's M840 graduate course, Aperiodic Tilings and Symbolic Dynamics. In the course topic guide, (co-authored with Reem Yassawi), he showed how, via the Fibonacci word's substitution incidence matrix, the left eigenvector gave rise to an aperiodic tiling of $\mathbb{R}^+$.

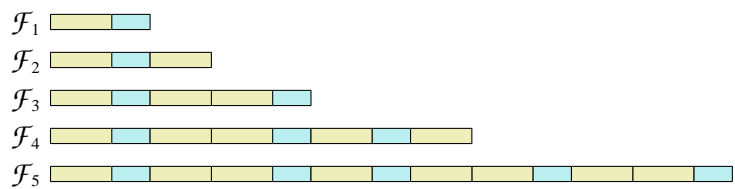

For my dissertation I ended up running with this idea, under the watchful eye of Dan Rust who kindly stepped in to supervise Uwe's orphaned students and keep the course going. I took Uwe's tiled path and twisted it back and forth, often with it tiling over itself, and looked at the properties of the resulting figures. The twisting was via a drawing rule that took each letter of the Fidonacci word in turn and used it as an instruction to say how the next tile should be placed. Some attractive images resulted. To give a flavour of what can occur the adjacent image is for $\mathcal{F}_{22}$ under the following drawing rule,

| Symbol | Action |
|--------|--------|
| $a$ | forward $\phi$   (The golden ratio, about 1.618) |
| $b$ | forward 0.5, turn 108°, forward 0.5 |

I like to thing that my visualisation of $\mathcal{F}_{22}$ is in the spirit of the mathematics that inspired it; Uwe's mathematics. And that he would approve.

<div style="text-align:right">Martin Hansen, June 2022</div>

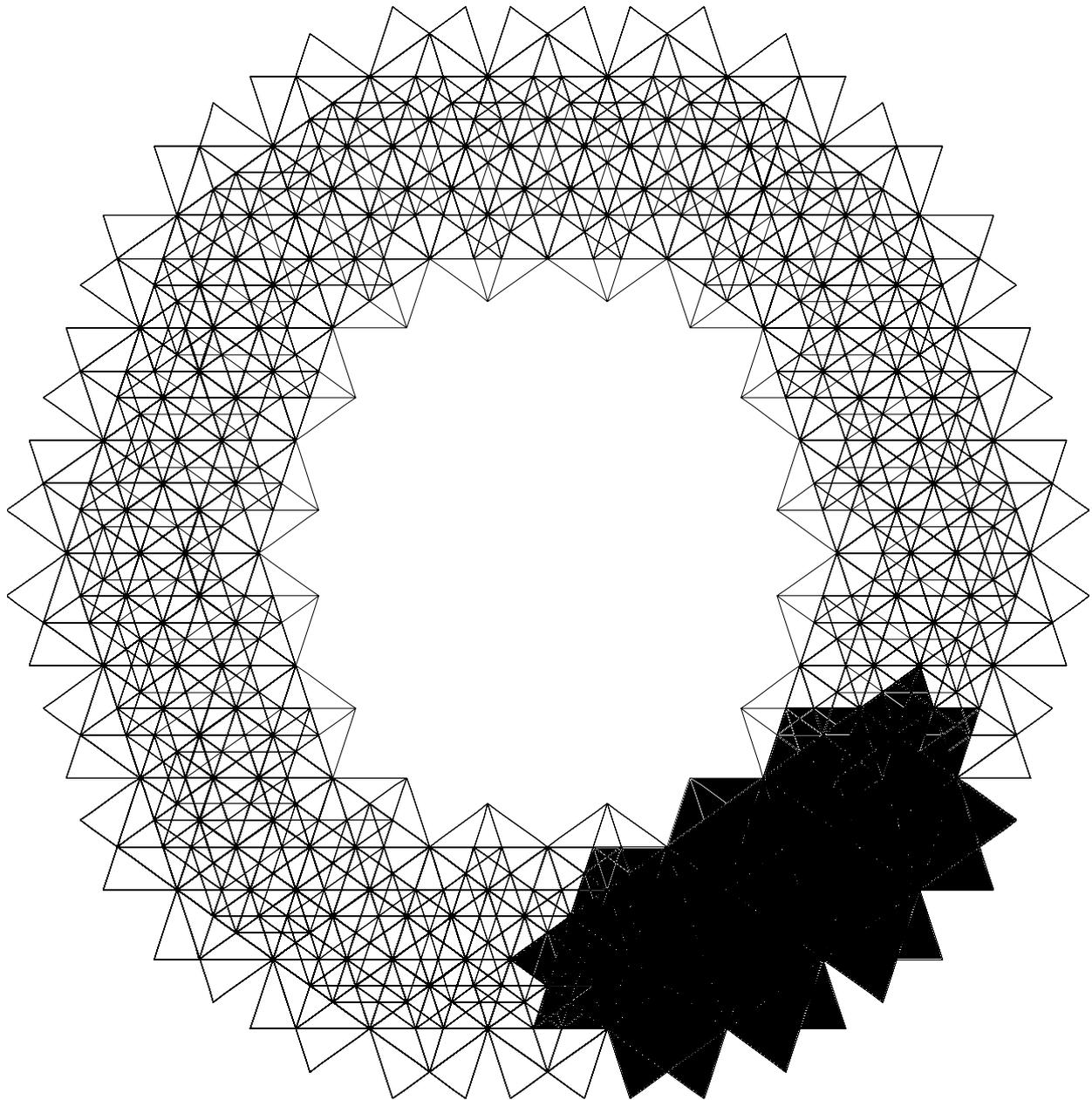

"Remembering Professor Uwe Grimm" was produced as a poster for *The Aperiodic*, a conference organised in memory of Uwe Grimm at The Open University in Milton Keynes in June 2022.

It was subsequently published in the December 2022 edition of M500, the mathematics magazine of The Open University.

# Appendix A2 : Programs

### A2.1  To and Fro Drawings : Samples of Simplified Logo Program Code

| | |
|---|---|
| to Fib01right | : Program name, this will draw $\mathcal{F}_1 = ab$ |
| forward 1.618 | : Length of tile *A*, the golden ratio to three decimal places |
| forward 0.5 | : First half length of tile *B* |
| right 90 | |
| forward 0.4 | : To make an artificial step downward at the centre of tile *B* |
| right 90 | |
| forward 0.5 | : Second half length of tile *B* |
| end | |
| | |
| to Fib01left | : Similar to Fib01right with right 90 replaced with left 90 |
| | : Whether moving to the right or the left the line needs to |
| | : zigzag downwards |
| | |
| to Fib02right | : Program name, this will draw $\mathcal{F}_2 = aba$ |
| Fib01right | : Call program titled Fib01right |
| forward 1.618 | : The extra tile *A* on the end of $\mathcal{F}_2$ having just drawn $\mathcal{F}_1$ |
| end | |
| | |
| to Fib02left | : The moving left version of Fib02right |
| | |
| to Fib03right | : Drawing $\mathcal{F}_3$ using the relationship $\mathcal{F}_3 = \mathcal{F}_2 \mathcal{F}_1$ |
| Fib02right | : Call program titled Fib02right |
| Fib01left | : Call program titled Fib01left |
| end | |
| | |
| to Fib03left | : Similar to Fib03right with Fib02right replaced with Fib02left |
| | :              and with Fib01left replaced with Fib01right |

| Number Wonder Mathematics Software | |
|---|---|
| Title | Logo Programming Language |
| Author | Martin Hansen |
| License | The Open University |
| Version | 1.0 |
| Dated | 1st of March 2022 |

**Figure A2.1** : A LOGO programming language interpreter written to allow short "turtle graphics" programs to be written that implement and explore various drawing rules. Many of this dissertation's diagrams were drawn using this software.